\RequirePackage{fix-cm}
\documentclass[smallextended]{svjour3}       
\smartqed  
\usepackage[cmex10]{amsmath}
\usepackage{amssymb} 

\usepackage{algorithmic}
\usepackage[plain]{algorithm}

\usepackage{array}

\usepackage{url}

\DeclareMathOperator*{\argmin}{arg\,min}

\DeclareMathOperator{\diag}{diag}
\def\mplus{\mathrel{%
  \ooalign{\raise.29ex\hbox{$\scriptscriptstyle\mathbf{+}$}\cr}}}
\usepackage{color}

\usepackage{comment}

\usepackage{graphicx}
\usepackage{pgf,tikz}

 \usepackage{caption}
 \usepackage[caption=false]{subfig}

\usepackage{mathtools}
\usepackage{booktabs}

\begin{document}

\title{Performance of First- and Second-Order Methods for $\ell_1$-Regularized Least Squares Problems
}


\author{Kimon Fountoulakis         \and
        Jacek Gondzio 
}


\institute{Kimon Fountoulakis \at
              School of Mathematics and Maxwell Institute, The University of Edinburgh, James Clerk Maxwell Building, The King's Buildings, Peter Guthrie Tait Road, Edinburgh EH9 3JZ, Scotland UK\\
              Tel.: +44 131 650 5083\\
              \email{K.Fountoulakis@sms.ed.ac.uk}           
           \and
	      Jacek Gondzio \at
              School of Mathematics and Maxwell Institute, The University of Edinburgh, James Clerk Maxwell Building, The King's Buildings, Peter Guthrie Tait Road, Edinburgh EH9 3JZ, Scotland UK\\
              Tel.: +44 131 650 8574\\
              Fax: +44 131 650 6553\\
              \email{J.Gondzio@ed.ac.uk}           
}

\date{Received: date / Accepted: date}

\maketitle

\begin{abstract}
We study the performance of first- and second-order optimization methods 
for $\ell_1$-regularized sparse least-squares problems as the conditioning 
of the problem changes and the dimensions of the problem increase up to one trillion.
A rigorously defined generator is presented which allows control of the dimensions, the conditioning and the sparsity of the problem. 
The generator has very low memory requirements and scales well with the dimensions of the problem. 

\keywords{$\ell_1$-regularised least-squares \and First-order methods \and Second-order methods \and Sparse least squares instance generator \and Ill-conditioned problems}
\end{abstract}

\section{Introduction}\label{sec:introducation}
\begin{figure}%
\centering
\subfloat[$\ell_1$ regularized]{%
\label{fig0a}%
\includegraphics[scale=0.35]{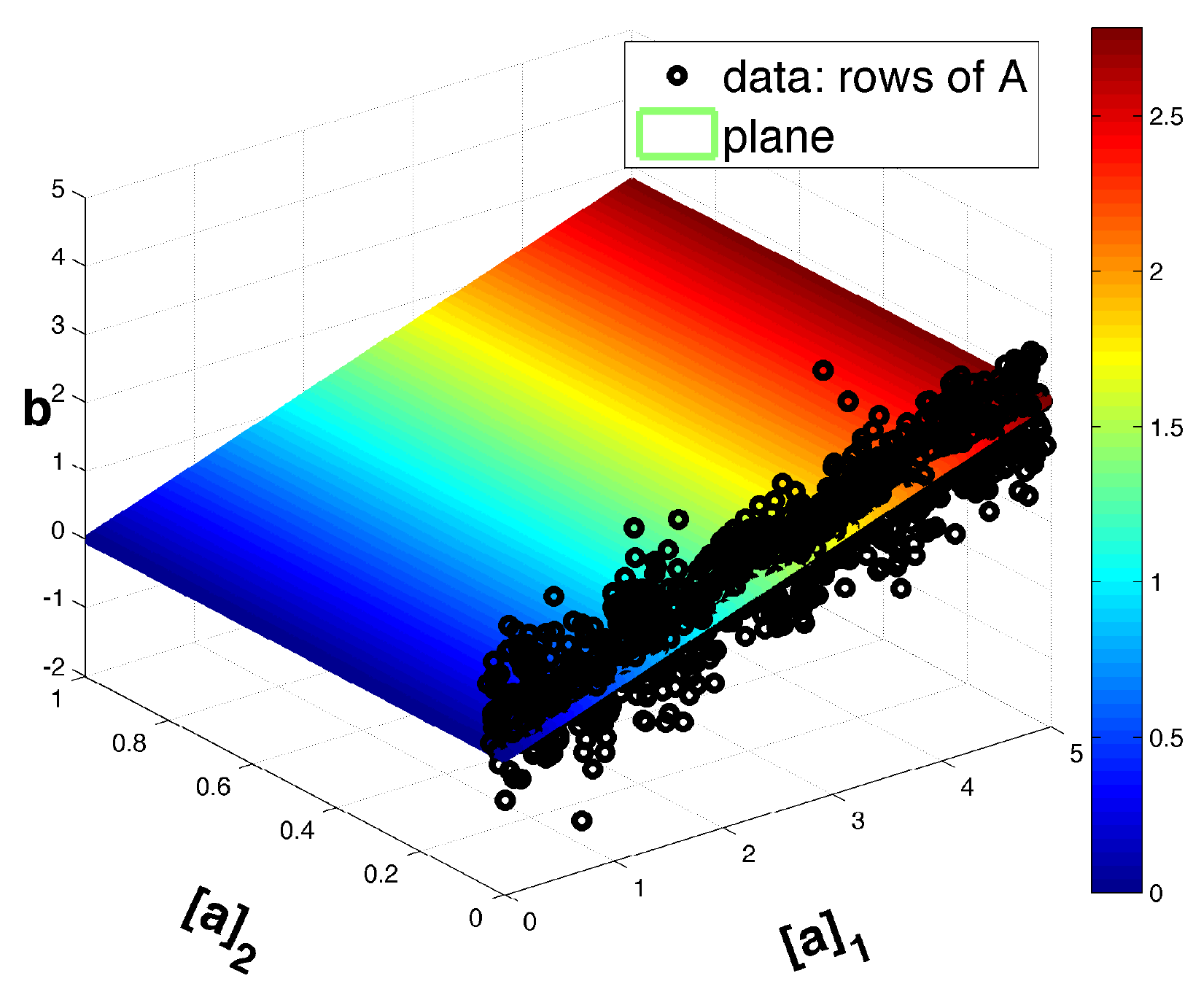}}
\hspace{0.1cm}
\subfloat[$\ell_2$ regularized]{%
\label{fig0b}%
\includegraphics[scale=0.35]{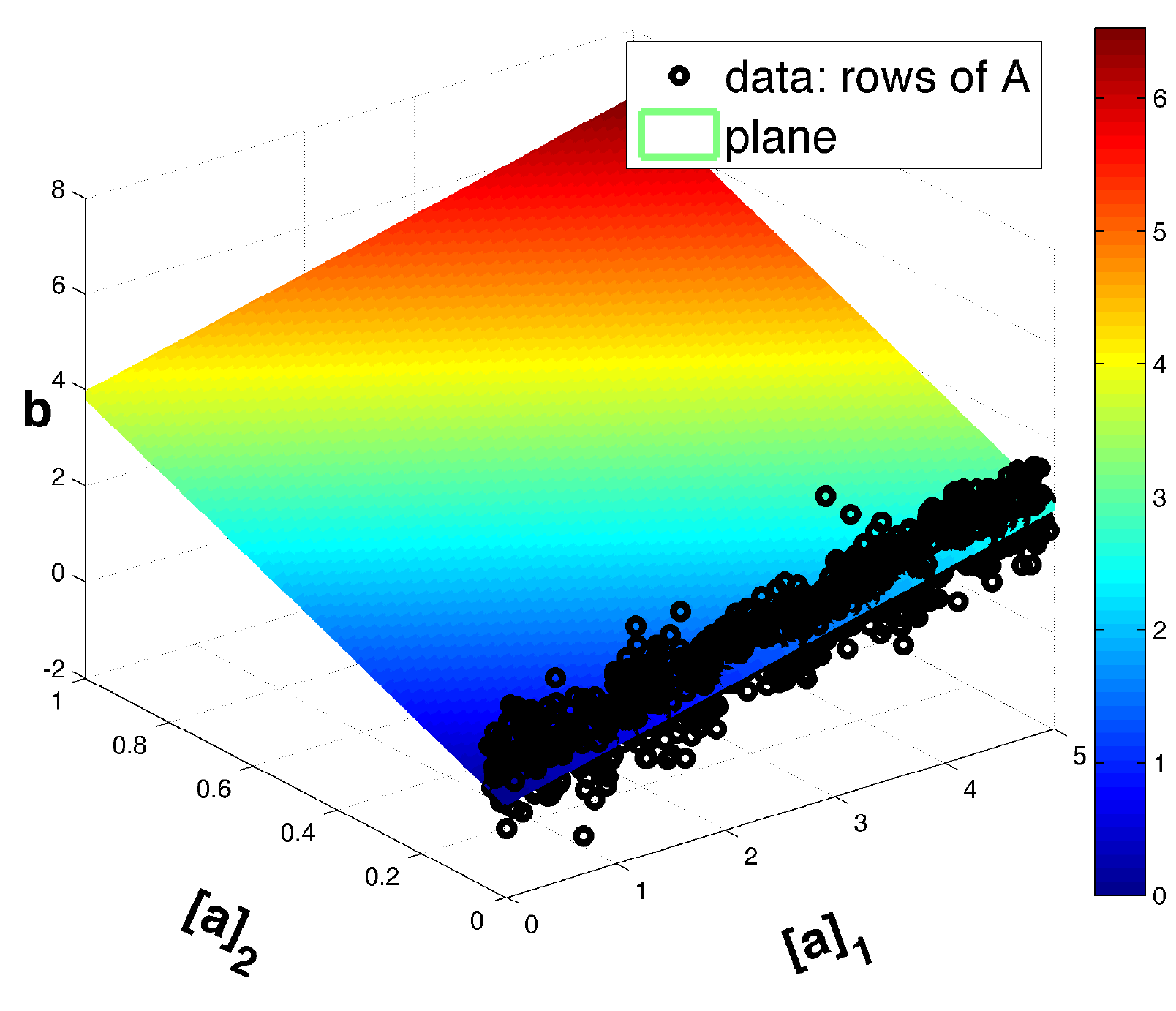}}\\
\caption{Demonstration of the purpose of the $\ell_1$ norm for data fitting problems.}
\label{fig0}%
\end{figure}

We consider the problem
\begin{equation}\label{prob1}
\mbox{minimize} \  f_\tau(x):=\tau \|x\|_1 + \frac{1}{2}\|Ax-b\|^2_2,
\end{equation}
where $x\in \mathbb{R}^n$, $\|\cdot\|_1$ denotes the $\ell_1$-norm, $\|\cdot\|_2$ denotes the Euclidean norm, $\tau>0$, 
$A\in\mathbb{R}^{m\times n}$ and $b\in\mathbb{R}^m$.
An application that is formulated as in \eqref{prob1} is sparse data fitting, where the aim is to approximate $n$-dimensional sampled points (rows of matrix $A$)
using a linear function, which depends on less than $n$ variables, i.e., its slope is a sparse vector. 
Let us assume that we sample $m$ data points $(a_i,b_i)$, where $a_i\in\mathbb{R}^{n}$ and $b_i\in\mathbb{R}$ $\forall i=1,2,\cdots,m$.
We assume linear dependence of $b_i$ on $a_i$:
$$
b_i = a_i^\intercal x + e_i \quad \forall i=1,2,\cdots,m,
$$
where $e_i$ is an error term due to the sampling process being innacurate.
Depending on the application some statistical information is assumed about vector $e$.
In matrix form the previous relationship is:
\begin{equation}\label{linearfit}
b = A x + e,
\end{equation}
where $A\in\mathbb{R}^{m \times n}$ is a matrix with $a_i$'s as its rows and $b\in\mathbb{R}^m$ is a vector with $b_i$'s
as its components. The goal is to find a sparse vector $x$ (with many zero components) such that the error $\|Ax - b\|_2$ is minimized.
To find $x$ one can solve problem \eqref{prob1}. 
The purpose of the $\ell_1$ norm in \eqref{prob1} is to promote sparsity in the optimal solution \cite{SparsityInducing}.
An example that demonstrates the purpose of the $\ell_1$ norm is presented in Figure \ref{fig0}. Figure \ref{fig0} shows a two dimensional 
instance where $n=2$, $m=1000$ and matrix $A$ is full-rank. Notice that the data points $a_i$ $\forall i=1,2,\cdots,m$
have large variations with respect to feature $[a_i]_1$ $\forall i$, where $[\cdot]_j$ is the $j$th component of the input vector, while there is only a small variation with respect to feature $[a_i]_2$ $\forall i$.
This property is captured when problem \eqref{prob1} is solved with $\tau=30$. The fitted plane in Figure \ref{fig0a}
depends only on the first feature $[a]_1$, while the second feature $[a]_2$ is ignored because $[x^*]_2=0$, where $x^*$ is the optimal solution of \eqref{prob1}. This can be observed through the level sets of the plane shown with the colored map; for each value of $[a]_1$ the level sets remain constant for all values 
of $[a]_2$. On the contrary, this is not the case when one solves a simple least squares problem ($\tau=0$ in \eqref{prob1}). Observe in Figure \ref{fig0a} that the fitted plane depends on both features $[a]_1$ and $[a]_2$.

A variety of sparse data fitting applications originate from the fields of signal processing and statistics.
Five representative examples are briefly described below.
\begin{itemize}\setlength\itemsep{0.001em}
\item[-] Magnetic Resonance Imaging (MRI): A medical imaging tool used to 
scan the anatomy and the physiology of a body \cite{sparseMRI}.
\item[-] Image inpainting: A technique 
for reconstructing degraded parts of an image \cite{inpainting}.
\item[-]  Image deblurring: Image processing tool
for removing the blurriness of a photo caused by natural phenomena, such as motion \cite{deblurringimages}.  
\item[-] Genome-Wide Association study (GWA): DNA comparison 
between two groups of people (with/without a disease) in order to investigate factors that a disease 
depends on \cite{gwacs}.
\item[-]  Estimation of global temperature based on historic data \cite{datafitting}.
\end{itemize}

Data fitting problems frequently require the
\textit{analysis of large scale data sets}, i.e., gigabytes or terabytes of data. 
In order to address large scale problems 
there has been a resurgence in methods with computationally inexpensive iterations. 
For example many first-order methods were recovered and refined, such as coordinate descent \cite{gmlnet,HsiehChang,petermartin,tsengblkcoo,nesterovhuge,tsengyun,wrightaccel,tonglange}, 
alternating direction method of multipliers \cite{distributedadmm,admmtutorial,goldsteinadmm,convergenceadmm,admmtv}, proximal first-order methods \cite{fista,ista,proximalmethods} and 
first-order smoothing methods \cite{IEEEhowto:Nesta,convexTemplates,nesterovSmooth}. The previous are just few representative examples, the list is too long for a complete demonstration, 
many other examples can be found in \cite{SparsityInducing,introtv}. 
Often the goal of modern first-order methods is to reduce the computational complexity per iteration,
while preserving the theoretical worst case iteration complexity of classic first-order methods \cite{booknesterov}.
Many modern first order methods meet the previous goal. For instance, coordinate descent methods can have
up to $n$ times less computational complexity per iteration \cite{peterbigdata,petermartin}.

First-order methods have been very successful in various scientific fields, such as support vector machine \cite{svmcomparison}, compressed sensing \cite{IEEEhowto:DonohoCompSens}, image processing \cite{ista} and data fitting \cite{datafitting}. 
Several new first-order type approaches have recently been proposed for various imaging problems in the special issue edited by M. Bertero, V. Ruggiero and L. Zanni \cite{coap-vol54}.
However, even for the simple unconstrained problems that arise in the previous fields there exist more challenging instances. Since first-order methods do not capture sufficient second-order information, their performance might degrade 
unless the problems are well conditioned \cite{2ndpaperstrongly}. 
On the other hand, the second-order methods capture the curvature of the objective function sufficiently well, 
but by consensus they are usually applied only on medium scale problems
or when high precision accuracy is required. In particular, it is frequently claimed \cite{fista,IEEEhowto:Nesta,convexTemplates,haleyin,shwartzTewari}
that the second-order methods do not scale favourably as the dimensions of the problem increase because of their high computational complexity
per iteration.
Such claims are based on an assumption that a \textit{full} second-order information has to be used.
However,
there is evidence \cite{2ndpaperstrongly,IEEEhowto:Jacekmf} that for non-trivial problems, \textit{inexact} second-order methods
can be very efficient.

In this paper we will exhaustively study the performance of first- and second-order methods. We will perform numerical experiments for large-scale problems with sizes up 
to one trillion of variables. We
will examine conditions under which certain methods are favoured or not. We hope that by the end of this paper the reader will have a clear view about the 
performance of first- and second-order methods. 

Another contribution of the paper is the development of a rigorously defined instance generator for problems of the form of \eqref{prob1}.
The most important feature of the generator is that it scales well with the size of the problem and can inexpensively create instances where the user controls the sparsity and the conditioning of the problem.
For example see Subsection \ref{subsec:bigproblem}, where an instance of one trillion variables is created using the proposed generator. 
We believe that the flexibility of the proposed generator will cover the need for generation of various good test problems. 

This paper is organised as follows. In Section \ref{sec:briefmethods} we briefly discuss the structure 
of first- and second-order methods. In Section \ref{sec:gen} we give the details of the instance generator. In Section \ref{subsec:conA} we provide examples for constructing 
matrix $A$. In Section \ref{subset:cond}, we present some measures of the conditioning of problem \eqref{prob1}.
These measures will be used to examine the performance of the methods in the numerical experiments. 
In Section \ref{sec:contructsol} we discuss how the optimal solution of the problem is selected. 
In Section \ref{sec:existinggen} we briefly describe known problem generators and explain how our propositions add value to the existing approaches.
In Section \ref{sec:num} we present the practical performance of first- and second-order methods
as the conditioning and the size of the problems vary. Finally, in Section \ref{sec:con} we give our conclusions.

\section{Brief discussion on first- and second-order methods}\label{sec:briefmethods}

We are concerned with the performance of unconstrained optimization methods which have the following intuitive setting. 
At every iteration a convex function $Q_\tau(y; x): \mathbb{R}^n \to \mathbb{R}$ is created that locally approximates $f_\tau$ at a given point $x$.
Then, function $Q_\tau$ is minimized to obtain the next point. An example that covers the previous setting is the Generic Algorithmic Framework (GFrame) which is given below. Details
of GFrame for each method used in this paper are presented in Section \ref{sec:num}.

\begin{algorithm}
\floatname{algorithm}{Algorithm}
\renewcommand{\thealgorithm}{}
\caption{Generic Framework (GFrame)}
\begin{algorithmic}[1]
\STATE Initialize $x_0\in\mathbb{R}^n$ and $y_0\in\mathbb{R}^n$\\
For $k=0,1,2,\ldots $ until some termination criteria are satisfied \\ 
\STATE \hspace{0.5cm} Create a convex function $Q_\tau(y;y_k)$ that approximates $f_\tau$ in a neighbourhood of $y_k$ 
\STATE \hspace{0.5cm} Approximately (or exactly) solve the subproblem 
\begin{equation}
\label{bd3}
 x_{k+1}\approx \argmin_{y} \  Q(y; y_k)
\end{equation}
\STATE \hspace{0.5cm} Find a step-size $\alpha>0$ based on some criteria and set
$$
y_{k+1} = x_k + \alpha ({x}_{k+1} - x_k)
$$
\\
end-for
\STATE Return approximate solution $x_{k+1}$
\end{algorithmic}
\end{algorithm} 

Loosely speaking, close to the optimal solution of problem \eqref{prob1}, the better the approximation $Q_\tau$ of $f_\tau$ at any point $x$ the fewer 
iterations are required to solve \eqref{prob1}.
On the other hand, the practical performance of such methods is a trade-off between careful incorporation of the curvature of $f_\tau$, i.e. 
second-order derivative information in $Q_\tau$ and the cost of solving 
subproblem \eqref{bd3} in GFrame. 

\begin{figure}%
\centering
\subfloat[Separable quadratic]{%
\label{fig_sepquad}%
\includegraphics[scale=0.35]{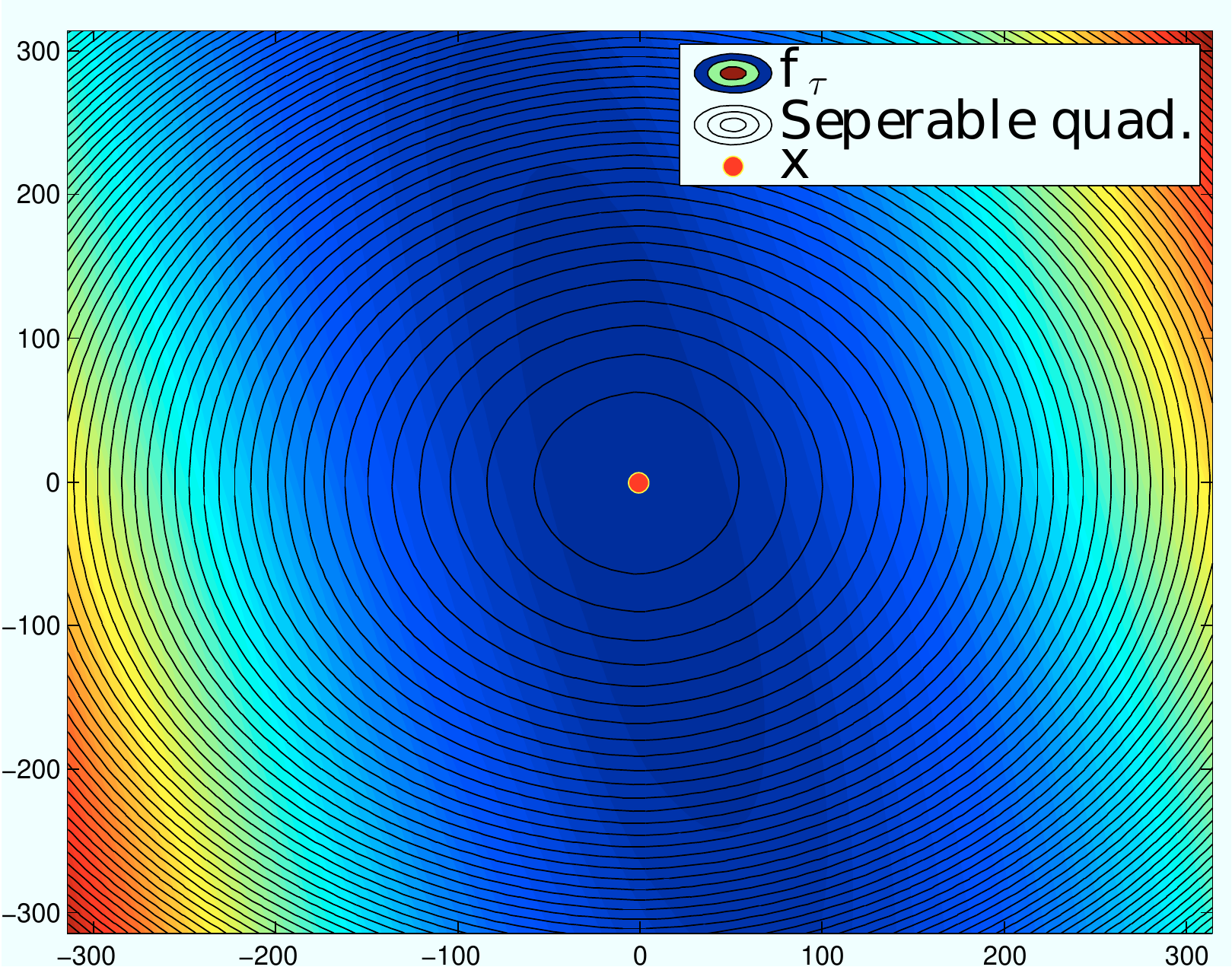}}
\hspace{0.1cm}
\subfloat[Non separable quadratic]{%
\label{fig_nonsepquad}%
\includegraphics[scale=0.35]{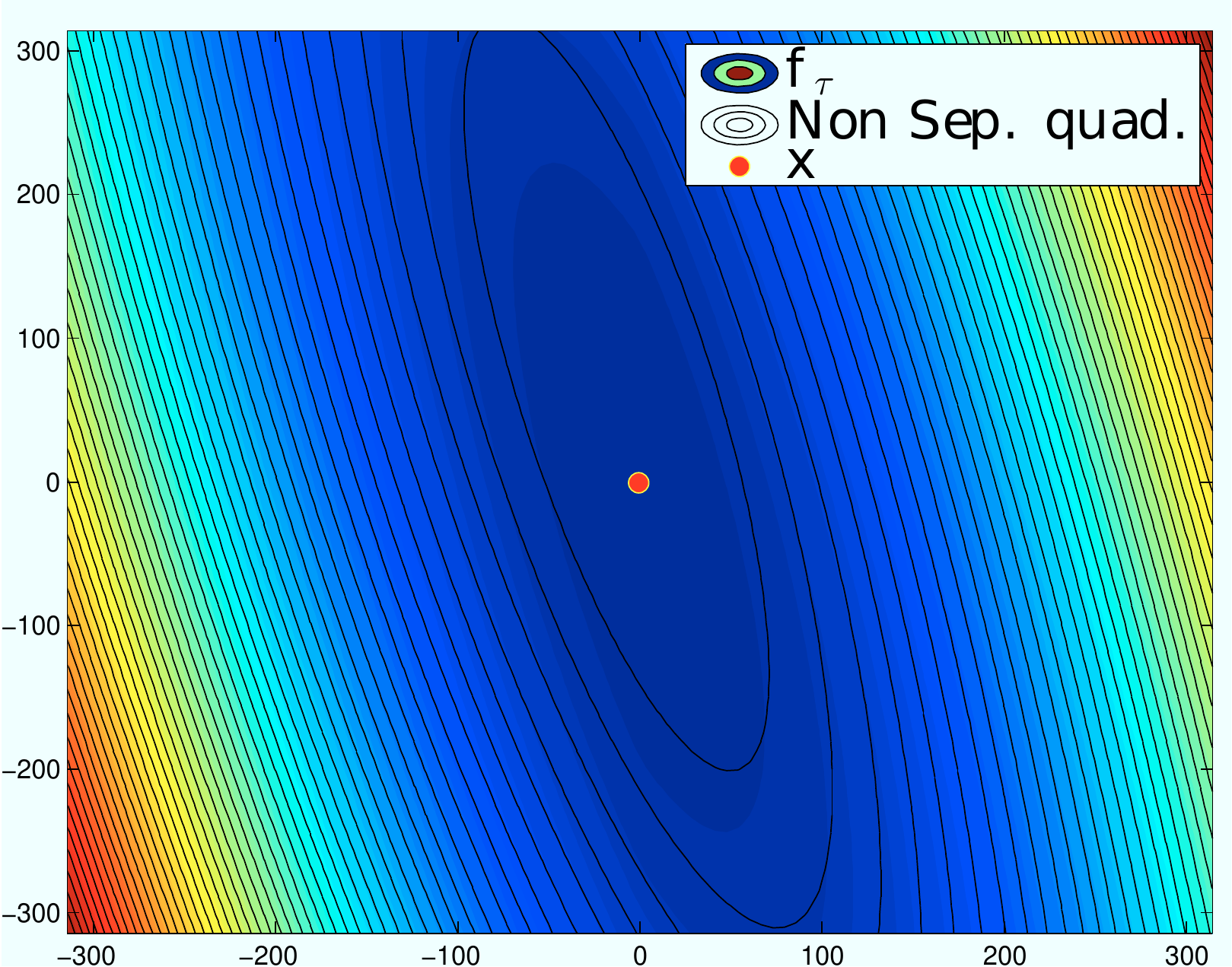}}\\
\caption{Demonstration of the contours of two different types of function $Q_\tau$, which locally approximate function $f_\tau$ at point $x$. In the left subfigure, function $Q_\tau$ is a simple separable quadratic for which, $H = \lambda_{max}(A^\intercal A)I_n$
in \eqref{bd1},
that is frequently used in first-order methods. In the right subfigure, function $Q_\tau$ is a non separable quadratic \eqref{bd4} which is used in some of the second-order methods. }
\label{fig_quad}%
\end{figure}

Discussion on two examples of $Q_\tau$ which consider different trade-off follows. 
First, let us fix the structure of $Q_\tau$ for problem \eqref{prob1} to be
\begin{equation}
\label{bd1}
Q_\tau(y;x) := \tau \|y\|_1 + \frac{1}{2}\|Ax-b\|^2 + (A^\intercal(Ax-b))^\intercal (y-x) + \frac{1}{2} (y-x)^\intercal H (y-x),
\end{equation}
where $H \in\mathbb{R}^{n\times n}$ is a positive definite matrix.
Notice that the decision of creating $Q_\tau$ has been reduced to a decision of selecting $H$.
Ideally, matrix $H$ should be chosen such that
it represents curvature information of $1/2\|Ax-b\|^2$ at point $x$, i.e. matrix $H$ should have similar spectral decomposition to $A^\intercal A$. 
Let $B(x):=\{v \in \mathbb{R}^n \ | \ \|v - x\|_2^2 \le 1\}$ be a unit ball centered at $x$. Then, $H$ should be selected in an optimal way: 
\begin{align}
\label{bd2}
\min_{H\succeq 0} \int_{B}|f_\tau(y) - Q_\tau(y;x)| dB.
\end{align}
The previous problem simply states that $H$ should minimize the sum of the absolute values of the residual $f_\tau - Q_\tau$ over $B$.
Using twice the fundamental theorem of calculus on $f_\tau$ from $x$ to $y$ we have that \eqref{bd2} is equivalent to
\begin{align}
\label{bd5}
\min_{H\succeq 0} \frac{1}{2} \int_{B} \Big| (y-x)^\intercal(A^\intercal A - H)(y-x)\Big| dB.
\end{align}
It is trivial to see that the best possible $H$ is simply $H=A^\intercal A$. However, this makes every subproblem \eqref{bd3}
as difficult to be minimized as the original problem \eqref{prob1}. One has to reevaluate the trade-off between a matrix $H$ that sufficiently well represents curvature information of $1/2\|Ax-b\|^2$
at a point $x$ compared to a simple matrix $H$ that is not as good approximation but offers an inexpensive solution of subproblem \eqref{bd3}.
An example can be obtained by setting $H$ to be a positively scaled identity, 
which gives a solution to problem \eqref{bd2} $H = \lambda_{max}(A^\intercal A)I_n$, where $\lambda_{max}(\cdot)$ denotes the largest eigenvalue of the input matrix
and $I_n$ is the identity matrix of size $n\times n$. The contours of such a function $Q_\tau$ compared to those of function $f_\tau$ are presented in Subfigure \ref{fig_sepquad}.
Notice that the curvature information of function $f_\tau$ is lost, this is because nearly all spectral properties of $A^\intercal A$ are lost.
However, for such a function $Q_\tau$ the subproblem \eqref{bd3} has an inexpensive closed form solution known as iterative shrinkage-thresholding \cite{ista,sparseMRI}. The computational complexity per iteration is so low that one hopes that this
will compensate for the losses of curvature information. Such methods, are called first-order methods and have been shown to be efficient for some large scale problems of the form of \eqref{prob1} \cite{fista}.

Another approach of constructing $Q_\tau$ involves the approximation of the $\ell_1$-norm with the pseudo-Huber function 
\begin{equation}\label{pseudoHuber}
 \psi_\mu(x) = \sum_{i=1}^n \Big((\mu^2+{x_i^2})^{\frac{1}{2}} - \mu\Big),
\end{equation}
where $\mu>0$ is an approximation parameter. 
This approach is frequently used by methods that aim in using at every iteration full information from the Hessian matrix $A^\intercal A$, see for example
\cite{2ndpaperstrongly,IEEEhowto:Jacekmf}.
Using \eqref{pseudoHuber}, problem \eqref{prob1} is replaced with
\begin{equation}\label{prob2}
 \mbox{minimize} \  f_{\tau}^{\mu}(x) := \tau \psi_{\mu}(x) + \frac{1}{2}\|Ax-b\|^2.
\end{equation}
The smaller $\mu$ is the better the approximation of problem \eqref{prob2} to \eqref{prob1}. 
The advantage is that $f_\tau^\mu$ in \eqref{prob2} is a smooth function which has derivatives of all degrees. 
Hence, smoothing will allow access to second-order information and essential curvature information will be exploited. 
However, for very small $\mu$
certain problems arise for optimization methods of the form of GFrame, see \cite{2ndpaperstrongly}. For example, 
the optimal solution of \eqref{prob1} is expected to have many zero components, on the other hand, the optimal solution 
of \eqref{prob2} is expected to have many nearly zero components. However, for small $\mu$ one can expect to obtain
a good approximation of the optimal solution of \eqref{prob1} by solving \eqref{prob2}.
For the smooth problem \eqref{prob2}, the convex approximation $Q_\tau$ at $x$ is:
\begin{equation}
\label{bd4}
Q_\tau(y;x) := f_{\tau}^{\mu}(x) + \nabla f_{\tau}^{\mu}(x)^\intercal (y-x) + \frac{1}{2} (y-x)^\intercal \nabla^2 f_{\tau}^{\mu}(x) (y-x).
\end{equation} 
The contours of such a function $Q_\tau$ compared to function $f_\tau$ are presented in Subfigure \ref{fig_nonsepquad}.
Notice that $Q_\tau$ captures the curvature information of function $f_\tau$.
However, minimizing
the subproblem \eqref{bd3} might be a more expensive operation. Therefore, we rely on an approximate solution of \eqref{bd3} using some iterative method which 
requires only simple matrix-vector product operations with matrices $A$ and $A^\intercal$.
In other words we use only an approximate second-order information. 
It is frequently claimed \cite{fista,IEEEhowto:Nesta,convexTemplates,haleyin,shwartzTewari}
that second-order methods do not scale favourably with the dimensions of the problem because of the more costly 
task of solving approximately the subproblems in \eqref{bd3}, instead of having an inexpensive closed form solution.
Such claims are based on an assumption that \textit{full} second-order information has to be used when solving subproblem \eqref{bd3}.
Clearly, this is not necessary: an \textit{approximate} second-order information suffices. 
Studies in \cite{2ndpaperstrongly,IEEEhowto:Jacekmf} provided theoretical background as well as the preliminary computational results 
to illustrate the issue. In this paper, we provide rich computational evidence which demonstrates that second-order methods can be very efficient.

\section{Instance Generator}\label{sec:gen}
In this section we discuss an instance generator for \eqref{prob1} 
for the cases $m\ge n$ and $m<n$. The generator is inspired by the one presented in Section $6$ of \cite{nesterovgen}. The advantage of our modified 
version is that it allows to control the properties of matrix $A$ and the optimal solution $x^*$ of \eqref{prob1}. For example, the sparsity of matrix $A$, its spectral decomposition, the sparsity and the norm of $x^*$,
since $A$ and $x^*$ are defined by the user.

Throughout the paper we will denote the $i^{th}$ component of a vector, by the name of the vector with subscript $i$. Whilst, the $i^{th}$ column of a matrix
is denoted by the name of the matrix with subscript $i$.

\subsection{Instance Generator for $m\ge n$}
\label{subsec:mgn}

Given $\tau>0$, $A\in\mathbb{R}^{m\times n}$ and $x^*\in\mathbb{R}^n$ the generator returns a vector $b\in\mathbb{R}^m$ such that 
$x^*:= \argmin_x f_\tau(x)$.
For simplicity we assume that the given matrix $A$ has rank $n$. 
The generator is described in Procedure IGen below.
\begin{algorithm}
\floatname{algorithm}{Procedure}
\renewcommand{\thealgorithm}{}
\caption{Instance Generator (IGen)}
\begin{algorithmic}[1]
\STATE Initialize $\tau>0$, $A\in\mathbb{R}^{m\times n}$ with $m\ge n$ and rank $n$, $x^*\in\mathbb{R}^n$
\STATE Construct $g\in\mathbb{R}^n$ such that $g\in \partial \|x^*\|_1$:
\hspace{0.5cm}
\begin{equation}
\label{bd6}
    g_i \in
\begin{cases}
    \{1\},& \text{if } x^*_i > 0\\
    \{-1\},& \text{if } x^*_i < 0\\
    [-1,1],& \text{if } x^*_i = 0
\end{cases}
\quad
\forall i=1,2,\cdots,n
\end{equation}
\STATE Set $e=\tau A(A^\intercal A)^{-1}g$
\STATE Return $b = Ax^* + e$
\end{algorithmic}
\end{algorithm} 

In procedure IGen, given $\tau$, $A$ and $x^*$ we are aiming in finding a vector $b$ such that $x^*$ satisfies the optimality conditions of problem \eqref{prob1}
\begin{equation*}
A^\intercal(Ax^* - b) \in -\tau \partial \|x^*\|_1,
\end{equation*}
where $\partial \|x\|_1=[-1,1]^n$ is the subdifferential of the $\ell_1$-norm at point $x$. By fixing a subradient $g\in\partial \|x^*\|_1$ as defined in \eqref{bd6}
and setting $e=b - Ax^*$, the previous optimality conditions can be written as
\begin{equation}
\label{bd7}
A^\intercal e  = \tau g.
\end{equation}
The solution to the underdetermined system \eqref{bd7} is set to $e=\tau A(A^\intercal A)^{-1}g$ and then we simply obtain $b=Ax^* + e$; Steps $3$ and $4$ in IGen, respectively.
Notice that for a general matrix $A$, Step $3$ of IGen can be very expensive. Fortunately, using elementary linear transformations, such as Givens rotations, 
we can iteratively construct a sparse matrix $A$ with a known singular value decomposition and guarantee that the inversion of matrix $A^\intercal A$ in Step $3$ of IGen is trivial.
We provide a more detailed argument in Section \ref{subsec:conA}.

\subsection{Instance Generator for $m < n$}\label{app:igen}
In this subsection we extend the instance generator that was proposed in Subsection \ref{subsec:mgn} 
to the case of matrix $A\in\mathbb{R}^{m\times n}$ with more columns than rows, i.e. $m<n$. 
Given $\tau>0$, $B\in\mathbb{R}^{m\times m}$, $N\in\mathbb{R}^{m\times n-m}$ and  $x^*\in\mathbb{R}^n$ the generator returns a vector $b\in\mathbb{R}^m$ and a matrix $A\in\mathbb{R}^{m\times n}$ such that 
$x^*:= \argmin_x f_\tau(x)$.

For this generator we need
to discuss first some restrictions on matrix $A$ and the optimal solution $x^*$.
Let 
\begin{equation}\label{bd13}
S:=\{i \in \{1,2,\cdots,n\} \ | \ x_i^* \neq 0 \}
\end{equation}
with $|S|=s$ and $A_S\in\mathbb{R}^{m\times s}$ be a collection of columns from matrix $A$ which correspond
to indices in $S$. Matrix $A_S$ must have rank $s$ otherwise problem \eqref{prob1} is not well-defined. To see this, 
let $\mbox{sign}(x^*_S)\in\mathbb{R}^s$ be the sign function applied component-wise to $x^*_S$, where $x^*_S$ is a vector with components of $x^*$ that correspond to indices in $S$.
Then problem \eqref{prob1} reduces to the following:
\begin{equation}
\label{bd9}
\mbox{minimize} \ \tau \mbox{sign}(x^*_S)^\intercal x_s + \frac{1}{2}\|A_Sx_s-b\|^2_2,
\end{equation}
where $x_s\in\mathbb{R}^s$. The first-order stationary points of problem \eqref{bd9} satisfy
$$
A_S^\intercal A_S x_s = -\tau \mbox{sign}(x^*_S) + A_S^\intercal b.
$$ 
If $\mbox{rank}(A_S)< s$, the previous linear system does not have a unique solution and problem \eqref{prob1} does not have a unique minimizer. 
Having this restriction in mind, let us now present the instance generator for $m < n$ in Procedure IGen2 below.
\begin{algorithm}[H]
\floatname{algorithm}{Procedure}
\renewcommand{\thealgorithm}{}
\caption{Instance Generator 2 (IGen2)}
\begin{algorithmic}[1]
\STATE Initialize $\tau>0$, $B\in\mathbb{R}^{m\times m}$ with rank $m$, $N\in\mathbb{R}^{m\times n-m}$, $x^*\in\mathbb{R}^n$ with
$S:= \{1,2,\cdots,s\}$ and $s\le m$
\STATE Construct $g\in\mathbb{R}^m$ such that $g\in\partial \|x^*(1,2,\cdots,m)\|_1$:
\hspace{0.5cm}
\begin{equation}
\label{bd11}
    g_i \in
\begin{cases}
    \{1\},& \text{if } x^*_i > 0\\
    \{-1\},& \text{if } x^*_i < 0\\
    [-1,1],& \text{if } x^*_i = 0
\end{cases}
\quad
\forall i=1,2,\cdots,m
\end{equation}
\STATE Set $e=\tau B^{-\intercal}g$
\STATE Construct matrix $\tilde{N}\in\mathbb{R}^{m\times n-m}$ with the following loop:\\
For $k=1,2,\ldots, n-m $ \\
\hspace{0.5cm} $$\tilde{N}_i = \frac{\xi \tau}{ | N_i^\intercal e| } N_i, \ \mbox{where} \ \xi \ \mbox{is a random variable in} \ [-1,1]$$  \\
end-for
\STATE Return $A = [B,\tilde{N}]$ and $b = Ax^* + e$
\end{algorithmic}
\end{algorithm} 

In IGen2, given $\tau$, $B$, $N$ and $x^*$ we are aiming in finding a vector $b$ and a matrix $\tilde{N}$ such that for $A = [B,\tilde{N}]$, $x^*$ satisfies the optimality conditions of problem \eqref{prob1}
\begin{equation*}
A^\intercal(Ax^* - b) \in -\tau \partial \|x^*\|_1.
\end{equation*}
Without loss of generality it is assumed that all nonzero components of $x^*$ correspond to indices in $S=\{1,2,\cdots,s\}$. 
By fixing a partial subradient $g\in\partial \|x^*(1,2,\cdots,m)\|_1$ as in \eqref{bd11}, where $x^*(1,2,\cdots,m)\in\mathbb{R}^m$
is a vector which consists of the first $m$ components of $x^*$, 
and defining a vector $e=b-Ax^*$, the previous optimality conditions can be written as:
\begin{equation}
\label{bd10}
e  = \tau B^{-\intercal} g \quad \mbox{and} \quad \tilde{N}^\intercal e \in \tau [-1,1]^{n-m}.
\end{equation}
It is easy to check that by defining $\tilde{N}$ as in Step $4$ of IGen2 conditions \eqref{bd10} are satisfied. 
Finally, we obtain $b = Ax^* + e$.

Similarly to IGen in Subsection \eqref{subsec:mgn}, for Step $3$ in IGen2 we have to perform a matrix inversion, which generally can be an
expensive operation. However, in the next section we discuss techniques how this matrix inversion can be executed using a sequence of elementary 
orthogonal transformations.

\section{Construction of matrix $A$}\label{subsec:conA}
In this subsection we provide a paradigm on how matrix $A$ can be inexpensively constructed such that its singular value decomposition is known
and its sparsity is controlled. We examine the case of instance generator IGen where $m\ge n$. The paradigm can be 
easily extended to the case of IGen2, where $m<n$.

Let $\Sigma \in\mathbb{R}^{m\times n}$ be a rectangular matrix with the singular values $\sigma_1,\sigma_2,\cdots,\sigma_n$ on its diagonal and zeros elsewhere:
\[
\Sigma=
\left[\begin{array}{c}
\diag(\sigma_1,\sigma_2,\ldots,\sigma_n) \\
\hline
O_{m-n\times n} \\
\end{array}\right],
\]
where $O_{m-n\times n}\in\mathbb{R}^{m-n\times n}$ is a matrix of zeros, and let $G(i,j,\theta)\in\mathbb{R}^{n \times n}$ be a Givens rotation matrix, which rotates plane $i$-$j$ by an angle $\theta$:
\begin{equation*}
G(i,j,\theta) = 
\begin{bmatrix}   1   & \cdots &    0       & \cdots &    0      & \cdots &    0   \\
                      \vdots & \ddots & \vdots  &           & \vdots &            & \vdots \\
                         0      & \cdots &    c      & \cdots &    -s     & \cdots &    0   \\
                      \vdots &            & \vdots & \ddots & \vdots &            & \vdots \\
                         0      & \cdots &   s       & \cdots &    c      & \cdots &    0   \\
                      \vdots &            & \vdots &            & \vdots & \ddots & \vdots \\
                         0     & \cdots &    0       & \cdots &    0      & \cdots &    1
       \end{bmatrix},
\end{equation*}
where $i,j\in\{1,2,\cdots,n\}$, $c = \cos \theta$ and $s=\sin \theta$.
Given a sequence of Givens rotations $\{G(i_k,j_k,\theta_k)\}_{k=1}^K$ we define the following composition of them:
$$
G = G(i_1,j_1,\theta_1)G(i_2,j_2,\theta_2)\cdots G(i_K,j_K,\theta_K).
$$
Similarly, let $\tilde{G}(l,p,\vartheta)\in\mathbb{R}^{m \times m}$ be a Givens rotation matrix where $l,p\in\{1,2,\cdots,m\}$ and
$$
\tilde{G} = \tilde{G}(l_1,p_1,\vartheta_1)\tilde{G}(l_2,p_2,\vartheta_2)\cdots \tilde{G}(l_{\tilde{K}},p_{\tilde{K}},\vartheta_{\tilde{K}})
$$
be a composition of $\tilde{K}$ Givens rotations.
Using $G$ and $\tilde G$ we define matrix $A$ as
\begin{equation}
\label{bd8}
A = (P_1 \tilde{G} P_2) \Sigma G^\intercal,
\end{equation}
where $P_1,P_2\in\mathbb{R}^{m\times m}$ are permutation matrices.  
Since the matrices $P_1\tilde{G}P_2$ and $G$ are orthonormal it is clear that the left singular vectors of matrix $A$ are the columns of $P_1\tilde{G}P_2$, 
$\Sigma$ is the matrix of singular values and the right singular vectors are the columns of $G$. Hence, in Step $3$ of IGen we simply set 
$
(A^\intercal A)^{-1} = G(\Sigma^\intercal \Sigma)^{-1} G^\intercal,
$
which means that Step $3$ in IGen costs two matrix-vector products with $G$ and a diagonal scaling with $(\Sigma^\intercal \Sigma)^{-1}$.
Moreover, the sparsity of matrix $A$ is controlled by the numbers $K$ and $\tilde{K}$ of Givens rotations, the type, i.e. $(i,j,\theta)$ and $(l,p,\vartheta)$, and the order of Givens rotations.
Also, notice that the sparsity of matrix $A^\intercal A$ is controlled only by matrix $G$.
Examples are given in Subsection \ref{subsec:constA}.

It is important to mention that other settings of matrix $A$ in \eqref{bd8} could be used, for example different combinations of permutation matrices and Givens rotations.
The setting chosen in \eqref{bd8} is flexible, it allows for an inexpensive construction of matrix $A$ and makes the control of the singular value 
decomposition and the sparsity of matrices $A$ and $A^\intercal A$ easy. 

Notice that matrix $A$ does not have to be calculated and stored. 
In particular, in case that the method which is applied to solve problem \eqref{prob1} requires only matrix-vector product operations using matrices $A$ and $A^\intercal$, one can simply consider matrix $A$ as an operator. It is only required to predefine the triplets $(i_k,j_k,\theta_k)$ $\forall k=1,2,\cdots,K$ for matrix $G$, the triplets $(l_k,p_k,\theta_k)$ $\forall k=1,2,\cdots,\tilde{K}$ for matrix $\tilde{G}$ and the permutation matrices $P_1$ and $P_2$. 
The previous implies that the generator is inexpensive in terms of memory requirements.
Examples of matrix-vector product operations with matrices $A$ and $A^\intercal$ in case of \eqref{bd8} are given below in Algorithms MvPA and MvPAt, respectively.

\begin{algorithm}[H]
\floatname{algorithm}{Algorithm}
\renewcommand{\thealgorithm}{}
\caption{Matrix-vector product with $A$ (MvPA)}
\begin{algorithmic}[1]
\STATE Given a matrix $A$ defined as in \eqref{bd8} and an input vector $x\in\mathbb{R}^n$, do\\
\STATE Set $y_0 = x$ \\ 
For $k=1,2,\ldots, K $
\STATE \hspace{0.5cm} $y_k = G_k^\intercal y_{k-1}$ \\ 
end-for
\STATE Set $\tilde{y}_0 = P_2\Sigma y_K$
\STATE 
For $k=1,2,\ldots, \tilde{K} $
\STATE \hspace{0.5cm} $\tilde{y}_k = \tilde{G}_{\tilde{K}- k+1} \tilde{y}_{k-1}$ \\ 
end-for
\STATE Return $P_1\tilde{y}_{\tilde{K}}$
\end{algorithmic}
\end{algorithm} 

\begin{algorithm}[H]
\floatname{algorithm}{Algorithm}
\renewcommand{\thealgorithm}{}
\caption{Matrix-vector product with $A^\intercal$ (MvPAt)}
\begin{algorithmic}[1]
\STATE Given a matrix $A$ defined as in \eqref{bd8} and input vector $y\in\mathbb{R}^m$, do\\
\STATE Set $\tilde{x}_0 = P_1^\intercal y$\\
For $k=1,2,\ldots, \tilde{K} $
\STATE \hspace{0.5cm} $\tilde{x}_k = \tilde{G}_k^\intercal \tilde{x}_{k-1}$ \\ 
end-for
\STATE Set $x_0 = \Sigma^\intercal P_2^\intercal \tilde{x}_{\tilde{K}}$\\
For $k=1,2,\ldots, K $ \\
\STATE \hspace{0.5cm} $x_k = G_{K - k + 1} x_{k-1}$\\
end-for
\STATE Return $x_K$
\end{algorithmic}
\end{algorithm} 

\subsection{An example using Givens rotation}\label{subsec:constA}
Let us assume that $m,n$ are divisible by two and $m \ge n$. 
Given the singular values matrix $\Sigma$ and rotation angles $\theta$ and $\vartheta$, we construct matrix $A$ as
$$
A = (P\tilde{G}P)\Sigma G^\intercal ,
$$
where $P$ is a random permutation of the identity matrix, $G$ is a composition of $n/2$ Givens rotations:
$$
G = G(i_1,j_1,\theta)G(i_2,j_2,\theta)\cdots ,G(i_k,j_k,\theta),\cdots,G(i_{n/2},j_{n/2},\theta)
$$
with 
$$
i_k = 2k-1, \quad j_k = 2k \quad \mbox{for } \ k=1,2,3,\cdots,\frac{n}{2} 
$$
and $\tilde G$ is a composition of $m/2$ Givens rotations:
$$
\tilde{G} = \tilde{G}(l_1,p_1,\vartheta)\tilde{G}(l_2,p_2,\vartheta)\cdots ,\tilde{G}(l_k,p_k,\vartheta),\cdots,\tilde{G}(l_{m/2},p_{m/2},\vartheta)
$$
with 
$$
l_k = 2k-1, \quad p_k = 2k \quad \mbox{for } \ k=1,2,3,\cdots,\frac{m}{2}. 
$$
Notice that the angle $\theta$ is the same for all Givens rotations $G_k$, this means that the total memory requirement for matrix $G$ is low. In particular, it consists only of the storage of a $2\times 2$ rotation matrix.
Similarly, the memory requirement for matrix $\tilde{G}$ is also low.

\subsection{Control of sparsity of matrix $A$ and $A^\intercal A$}\label{subsec:constA_2}
We now present examples in which we demonstrate how sparsity of matrix $A$ can be controlled 
through Givens rotations. 

\begin{figure}
\centering
	\subfloat[$A=(P\tilde{G}P)\Sigma G^\intercal $]{%
			\label{spy_A_1}%
		\includegraphics[scale=0.66]{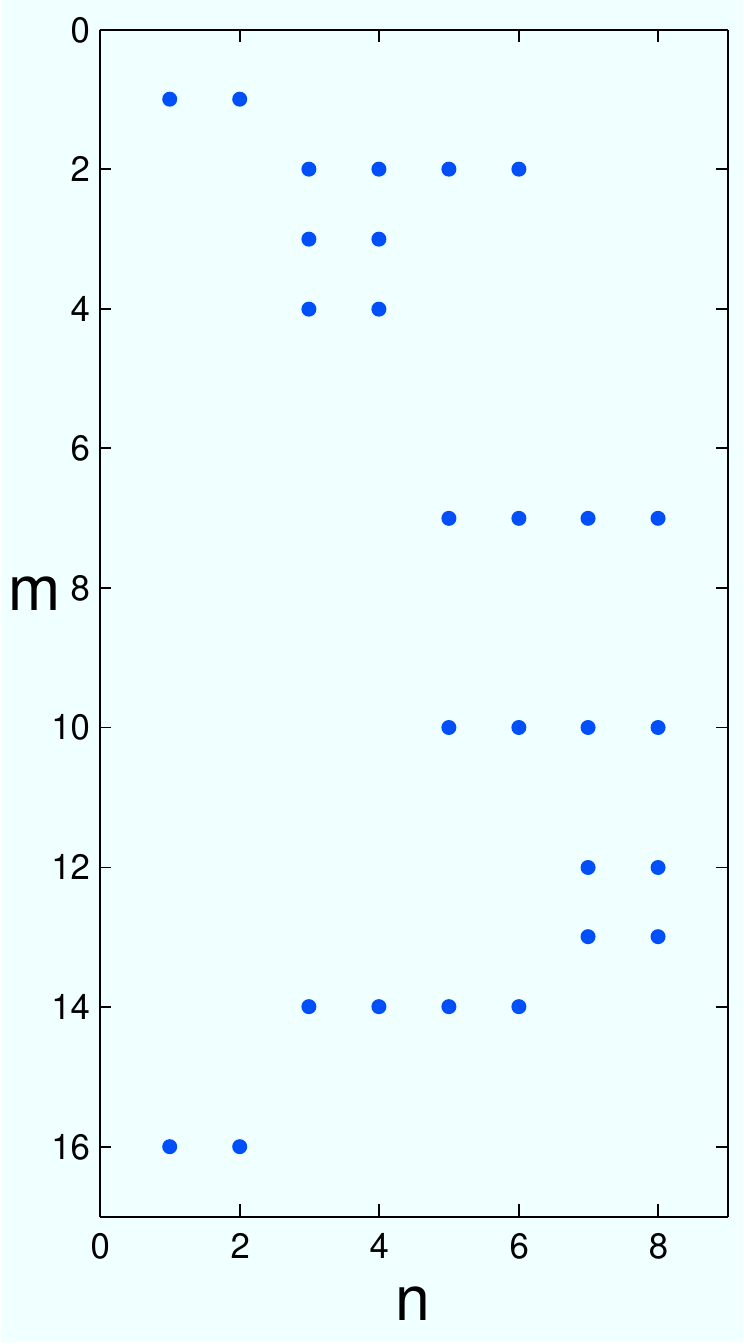}}		
         \quad
	\subfloat[$A=(P\tilde{G}_2\tilde{G}P)\Sigma (G_2G)^\intercal $]{%
		\label{spy_A_2}%
		\includegraphics[scale=0.66]{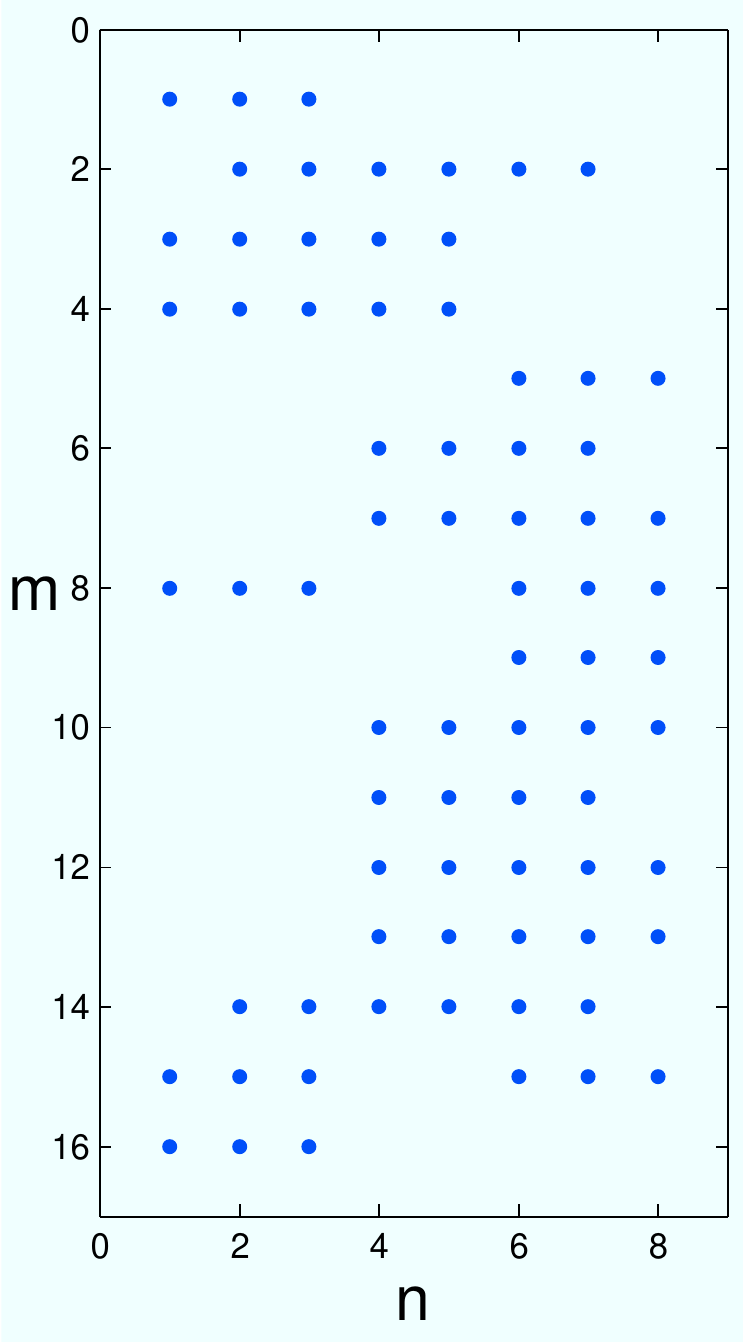}}
        \\
	\subfloat[$A=(P\tilde{G}\tilde{G}_2\tilde{G}P)\Sigma (GG_2G)^\intercal $]{%
		\label{spy_A_3}%
		\includegraphics[scale=0.66]{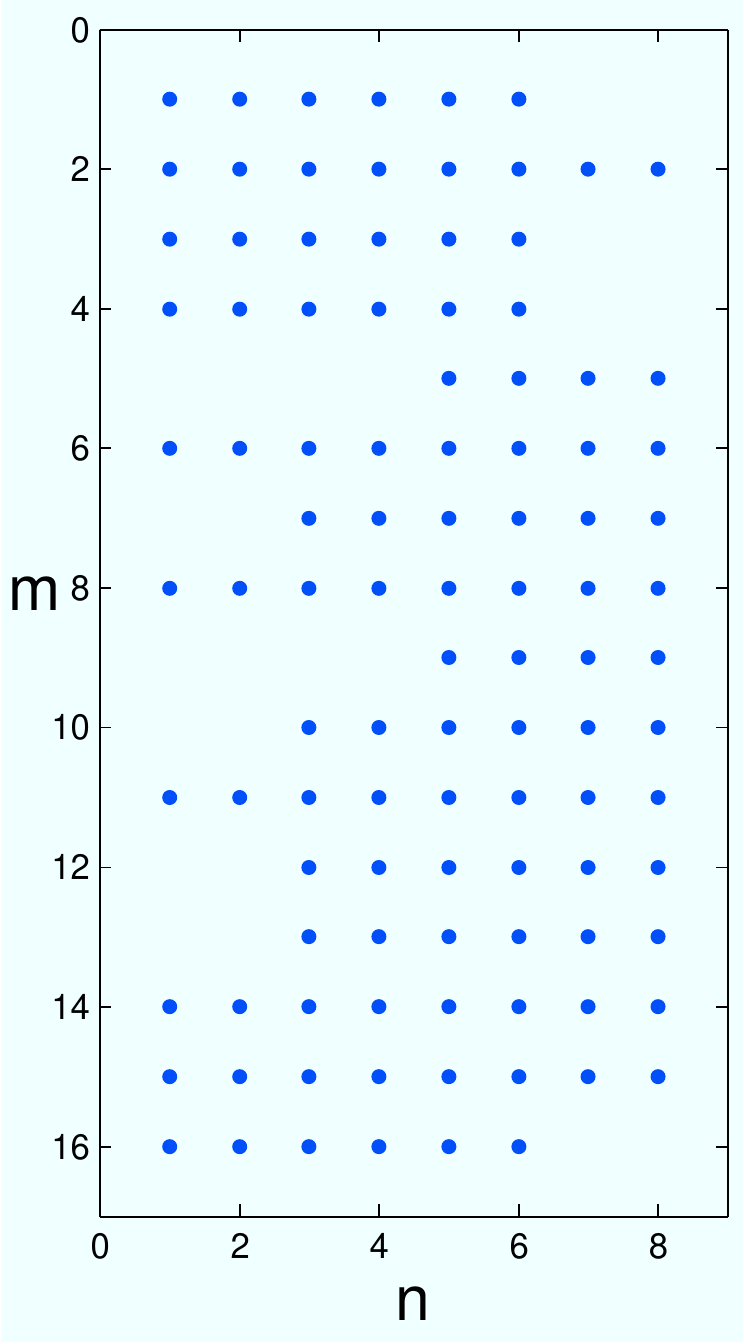}}
         \quad
	\subfloat[$A=(P\tilde{G}_2\tilde{G}\tilde{G}_2\tilde{G}P)\Sigma (G_2GG_2G)^\intercal $]{%
		\label{spy_A_4}%
		\includegraphics[scale=0.66]{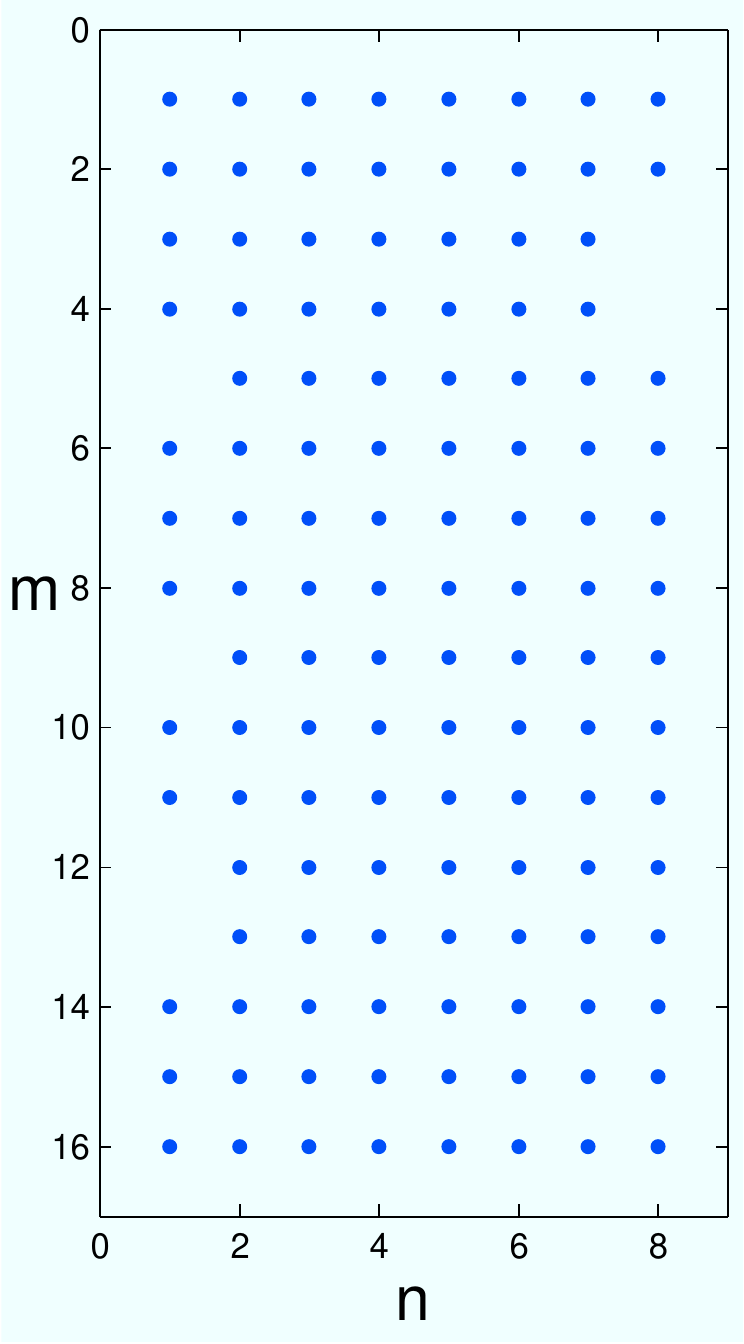}}
	\caption{Sparsity pattern of four examples of matrix $A$, the Givens rotations $G$ and $G_2$ are explained in Subsections \ref{subsec:constA} and \ref{subsec:constA_2}.}
	\label{figGivensrot}%
\end{figure}

\begin{figure}
\centering
	\subfloat[$A^\intercal A$, $A=(P\tilde GP)\Sigma G^\intercal $]{%
			\label{spy_AtA_1}%
		\includegraphics[scale=0.44]{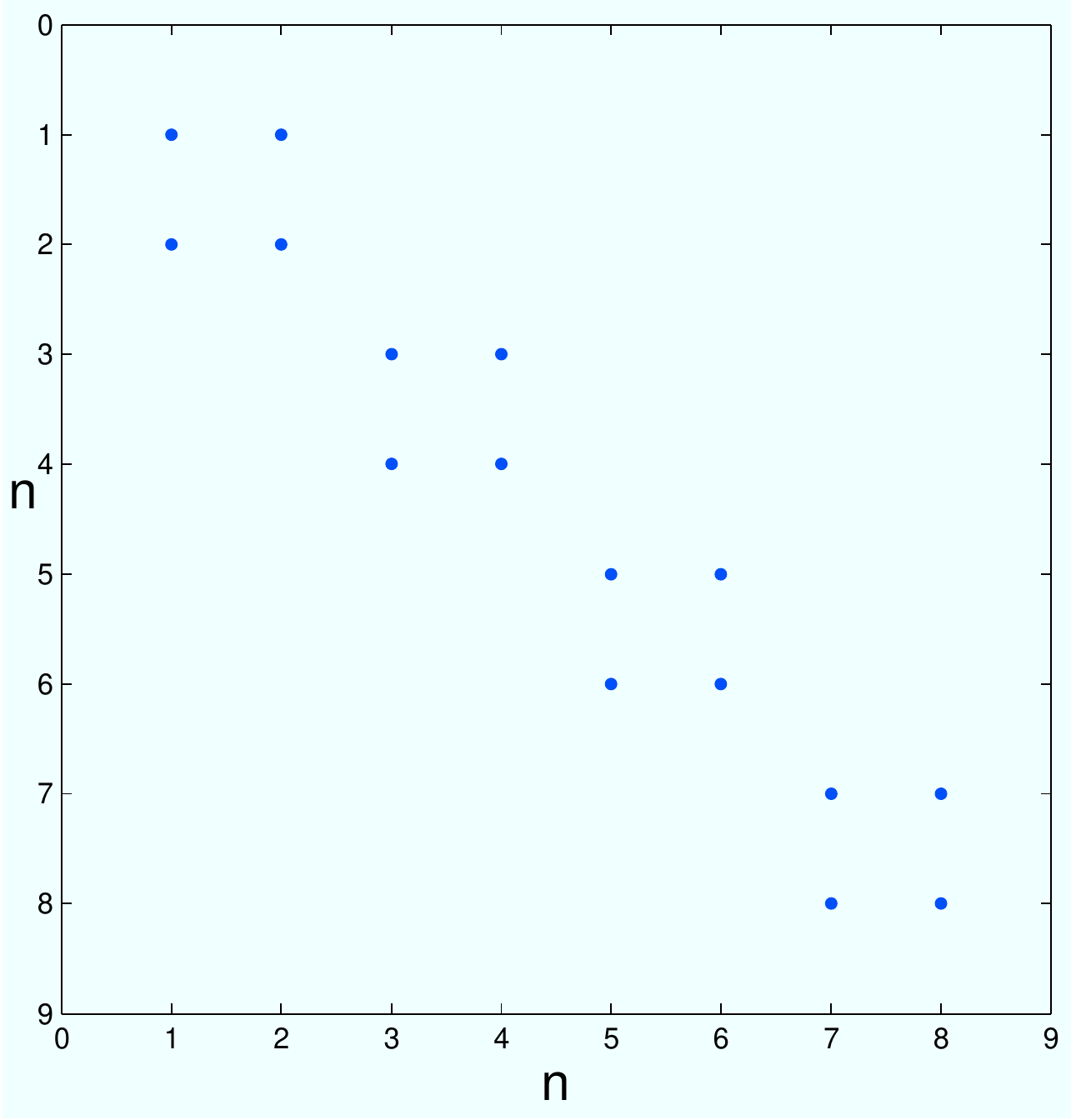}}		
         \quad
	\subfloat[$A^\intercal A$, $A=(P\tilde G_2 \tilde GP)\Sigma (G_2G)^\intercal $]{%
		\label{spy_AtA_2}%
		\includegraphics[scale=0.44]{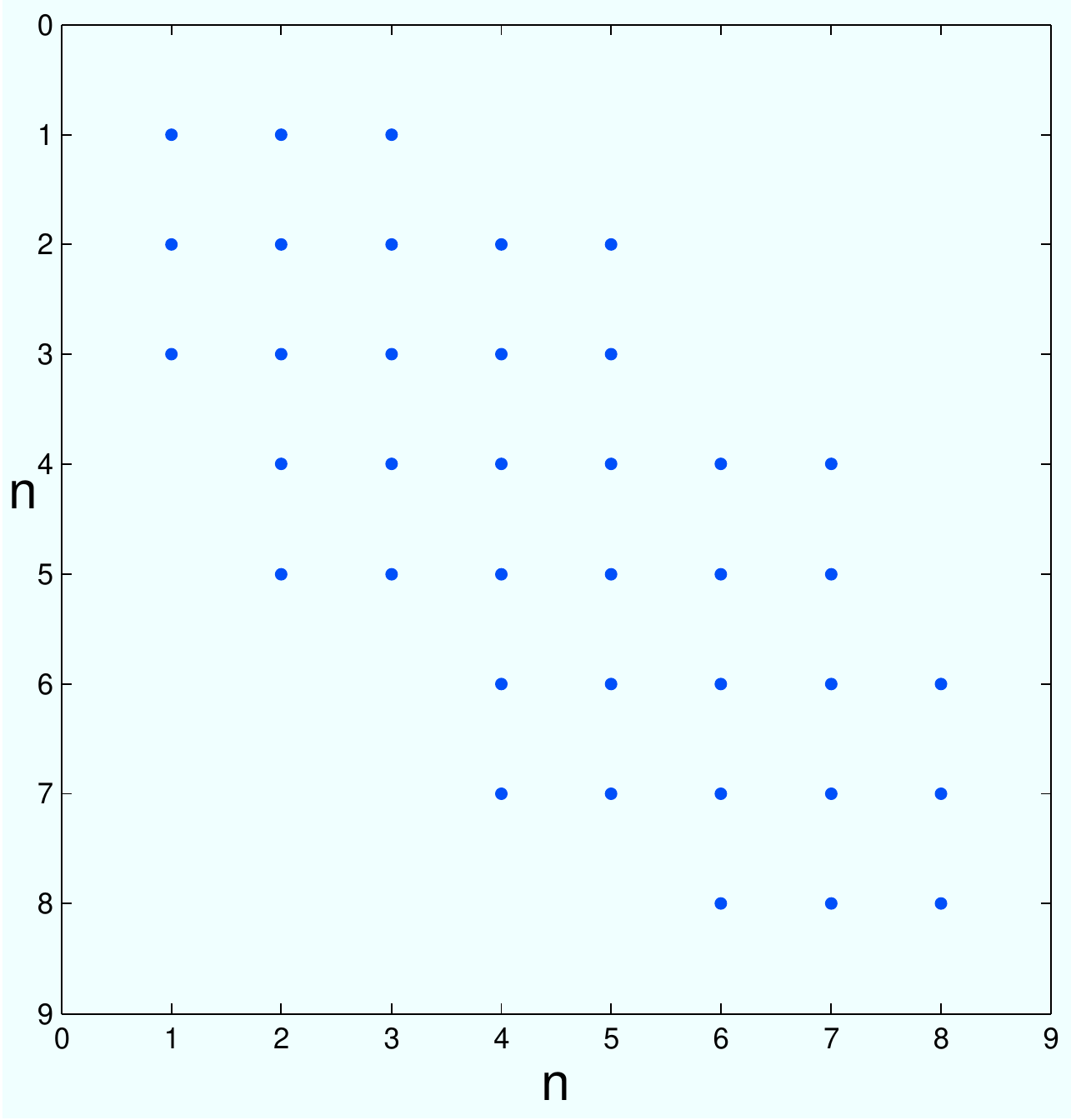}}
        \\
	\subfloat[$A^\intercal A$, $A=(P\tilde G \tilde G_2 \tilde GP) \Sigma (GG_2G)^\intercal $]{%
		\label{spy_AtA_3}%
		\includegraphics[scale=0.44]{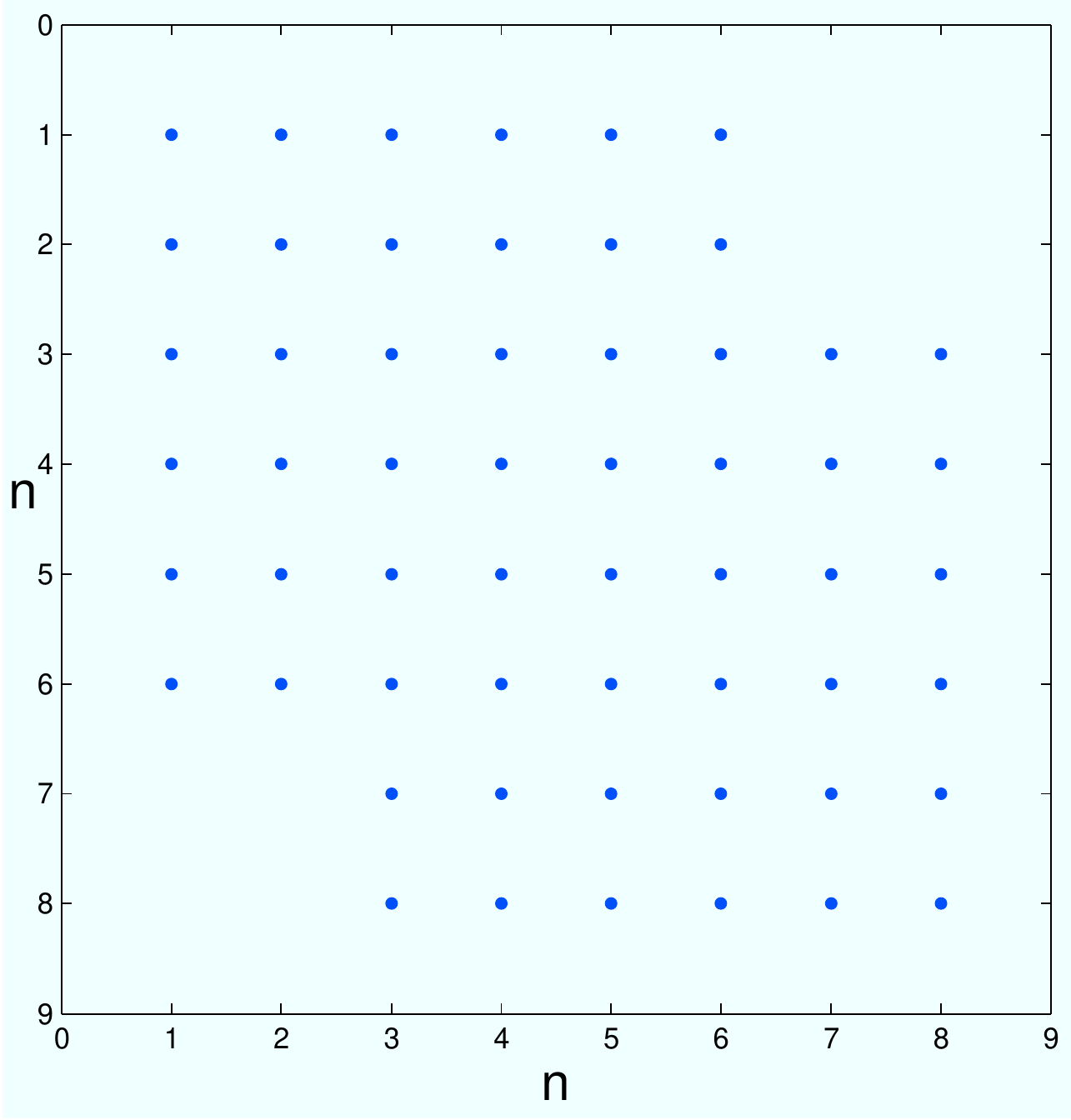}}
         \quad
	\subfloat[$A^\intercal A$, $A=(P\tilde G_2 \tilde G \tilde G_2 \tilde GP) \Sigma (G_2GG_2G)^\intercal $]{%
		\label{spy_AtA_4}%
		\includegraphics[scale=0.44]{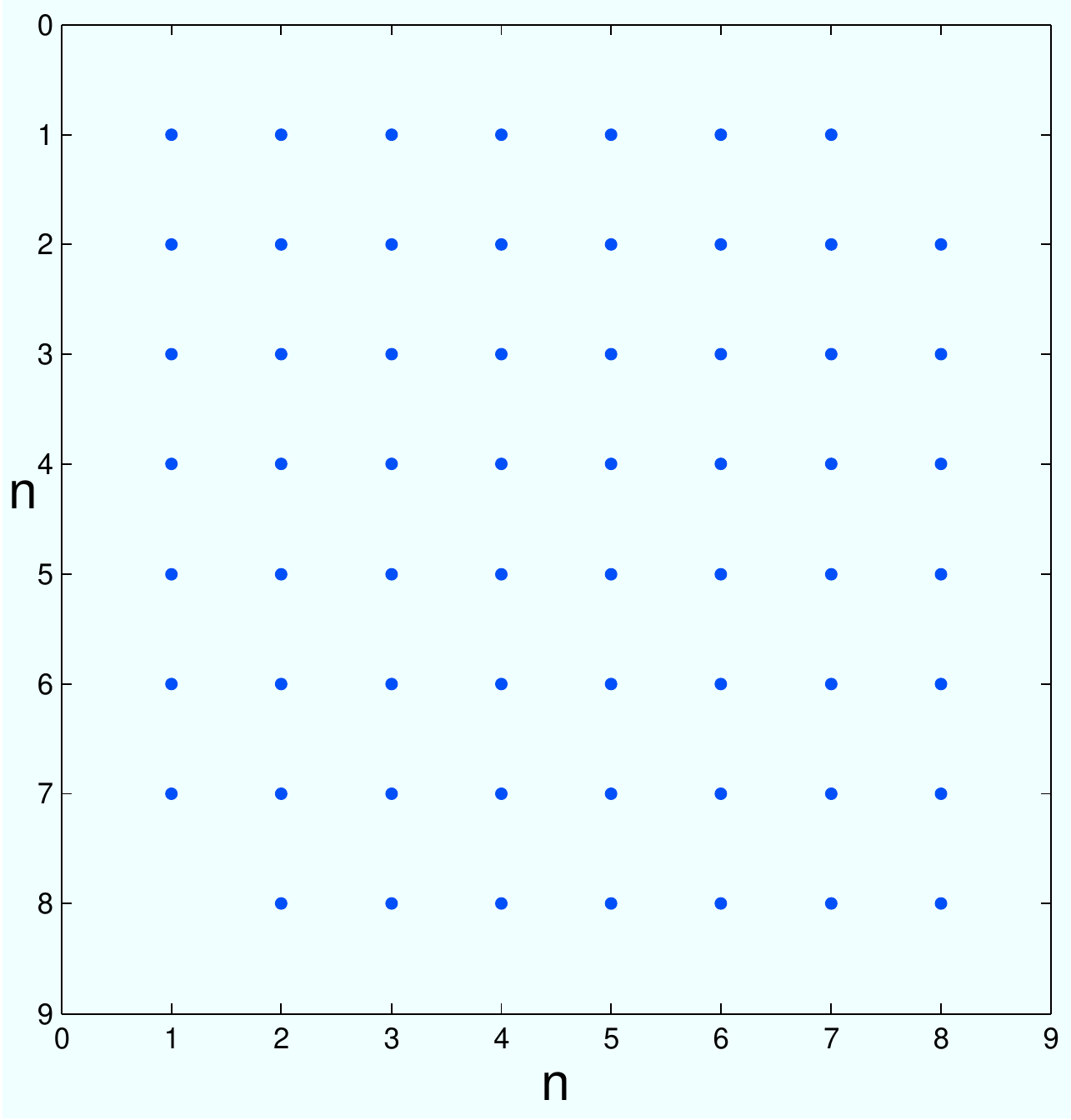}}
	\caption{Sparsity pattern of four examples of matrix $A^\intercal A$, where the Givens rotations $G$ and $G_2$ are explained in Subsections \ref{subsec:constA} and \ref{subsec:constA_2}}.
	\label{figGivensrotAtA}%
\end{figure}

In the example of Subsection \ref{subsec:constA}, two compositions of $n/2$ and $m/2$ Givens rotations, denoted by G and $\tilde{G}$, are applied
on an initial diagonal rectangular matrix $\Sigma$. If $n=2^3$ and $m=2n$ the sparsity pattern of the resulting matrix $A=(P\tilde{G}P)\Sigma G^\intercal$ is given in Subfigure \ref{spy_A_1}
and has $28$ nonzero elements, while the sparsity pattern of matrix $A^\intercal A$ is given in Subfigure \ref{spy_AtA_1} and has $16$ nonzero elements. 
Notice in this subfigure that the coordinates can be clustered in pairs of coordinates $(1,2)$, $(3,4)$, $(5,6)$ and $(7,8)$. One could apply 
another stage of Givens rotations. For example, one could construct matrix $A=(P\tilde{G}\tilde{G}_2P)\Sigma (G_2G)^\intercal$, where
$$
G_2 = G(i_1,j_1,\theta)G(i_2,j_2,\theta)\cdots ,G(i_k,j_k,\theta),\cdots,G(i_{n/2-1},j_{n/2-1},\theta)
$$
with 
$$
i_k = 2k, \quad j_k = 2k+1 \quad \mbox{for } \ k=1,2,3,\cdots,\frac{n}{2}-1. 
$$
and
$$
\tilde{G}_2 = \tilde{G}(l_1,p_1,\theta)\tilde{G}(l_2,p_2,\theta)\cdots ,\tilde{G}(l_k,p_k,\theta),\cdots,\tilde{G}(l_{m/2-1},p_{m/2-1},\theta)
$$
with 
$$
l_k = 2k, \quad p_k = 2k+1 \quad \mbox{for } \ k=1,2,3,\cdots,\frac{m}{2}-1. 
$$
Matrix $A=(P\tilde{G}\tilde{G}_2P)\Sigma (G_2G)^\intercal$ has $74$ nonzeros and it is shown in Subfigure \ref{spy_A_2}, while matrix $A^\intercal A$ has $38$
nonzeros and it is shown in Subfigure \ref{spy_AtA_2}.
By rotating again we obtain the matrix $A=(P\tilde{G}\tilde{G}_2\tilde{G}P)\Sigma (GG_2G)^\intercal$ in Subfigure \ref{spy_A_3} with $104$ nonzero
elements and matrix $A^\intercal A$ in Subfigure \ref{spy_AtA_3} with $56$ nonzero elements. 
Finally, the fourth Subfigures \ref{spy_A_4} and \ref{spy_AtA_4} show matrix $A=(P\tilde{G}_2\tilde{G}\tilde{G}_2\tilde{G}P)\Sigma (G_2GG_2G)^\intercal$ and $A^\intercal A$ with $122$ and $62$ nonzero elements, respectively. 

\section{Conditioning of the problem}\label{subset:cond}
Let us now precisely define how we measure the conditioning of problem \eqref{prob1}. 
For simplicity, throughout this section we assume that matrix $A$ has more rows than columns, $m\ge n$, and it is full-rank. Extension to the case of matrix $A$ with
more columns than rows is easy and we briefly discuss this at the end of this section.

We denote with $\mbox{span}(\cdot)$ the span of the columns of the input matrix.
Moreover, $S$ is defined in \eqref{bd13}, $S^c$ is its complement.

Two factors are considered that affect the conditioning of the problem. First, the usual condition number of the second-order derivative of $1/2\|Ax-b\|_2^2$ in \eqref{prob1}, which is simply
$\kappa(A^\intercal A) = \lambda_{1}(A^\intercal A)/\lambda_{n}(A^\intercal A)$, where $0<\lambda_n\le \lambda_{n-1} \le \cdots \le \lambda_1$ are the eigenvalues of matrix $A^\intercal A$. 
It is well-known that the larger $\kappa(A^\intercal A)$ is, the more difficult problem \eqref{prob1} becomes.

Second, the conditioning of the optimal solution $x^*$ of problem \eqref{prob1}. Let us explain what we mean by the conditioning of $x^*$.
We define a constant $\rho>0$ and the index set $\mathcal{I}_\rho:=\{i\in \{1,2,\cdots,n\} \ | \ \lambda_i(A^\intercal A) \ge \rho \}$.
Furthermore,
we define the projection $P_\rho = G_\rho G_\rho^\intercal$,
where $G_\rho\in\mathbb{R}^{n\times r}$, $r=|\mathcal{I}_{\rho}|$ and matrix $G_\rho$ has as columns the eigenvectors of matrix $A^\intercal A$
which correspond to eigenvalues with indices in $\mathcal{I}_\rho$. Then, the conditioning of $x^*$ is defined as
\begin{equation}\label{bd12}
\kappa_\rho(x^*) =
\begin{cases}
    \frac{\|x^*\|_2}{\|P_\rho x^*\|_2},& \text{if } P_\rho x^* \neq 0\\
    +\infty,& \mbox{otherwise}
\end{cases}
\end{equation}
For the case $P_\rho x^*\neq 0$, the denominator of \eqref{bd12} is the mass of $x^*$ which exists in the space spanned by eigenvectors of $A^\intercal A$ which correspond to eigenvalues that are larger than or equal to $\rho$.

Let us assume that there exists some $\rho$ which satisfies $\lambda_{n}(A^\intercal A) \le \rho \ll \lambda_{1}(A^\intercal A)$. 
If $\kappa_\rho(x^*) $ is large, i.e., $\|P_\rho x^*\|_2$ is close to zero, then the majority of the mass of $x^*$ is ``hidden" in 
the space spanned by eigenvectors which correspond to eigenvalues that are smaller than $\rho$, i.e., the orthogonal space of $\mbox{span}(G_\rho)$.
In Section \ref{sec:introducation} we referred to methods that do not incorporate information which correspond to small eigenvalues of $A^\intercal A$. 
Therefore, if the previous scenario holds, then we expect the performance of such methods to degrade.  
In Section \ref{sec:num} we empirically verify the previous arguments.

If matrix $A$ has more columns than rows then the previous definitions of conditioning of problem \eqref{prob1} are incorrect and need to be adjusted. 
Indeed, if $m<n$ and $\mbox{rank}(A)=\mbox{min}(m,n)=m$, then $A^\intercal A$ is a rank deficient matrix which has $m$ nonzero eigenvalues 
and $n-m$ zero eigenvalues. 
However, we can restrict the conditioning of the problem to a neighbourhood of the optimal solution of $x^*$. 
In particular, let us define a neighbourhood of $x^*$ so that all points in this neighbourhood have nonzeros at the same indices 
as $x^*$ and zeros elsewhere, i.e. $\mathcal{N}:=\{x\in\mathbb{R}^n \ | \ x_i\neq 0 \  \forall i\in S, \ x_i=0 \ \forall i \in S^c \}$.
In this case an important feature to determine the conditioning of the problem is the ratio of the largest and the smallest 
nonzero eigenvalues of $A_S^\intercal A_S$, where $A_S$ is a submatrix of $A$ built of columns of $A$ which belong to set $S$. 

\section{Construction of the optimal solution }\label{sec:contructsol}
Two different techniques are employed to generate the optimal solution $x^*$ for the experiments presented in Section \ref{sec:num}.
The first procedure suggests a simple random generation of $x^*$, see Procedure OsGen below.
\begin{algorithm}[H]
\floatname{algorithm}{Procedure}
\renewcommand{\thealgorithm}{}
\caption{Optimal solution Generator (OsGen)}
\begin{algorithmic}[1]
\STATE Given the required number $s\le \mbox{min}(m,n)$ of nonzeros in $x^*$ and a positive constant $\gamma>0$ do:
\STATE Choose a subset $S \subseteq \{1,2,\cdots,n\}$ with $|S|=s$.
\STATE $\forall i\in S$ choose $x_i^*$ uniformly at random in $[-\gamma,\gamma]$ and $\forall j \notin S$ set $x_j^*=0$. 
\end{algorithmic}
\end{algorithm} 

The second and more complicated approach is given in Procedure OsGen2. This procedure is applicable only 
in the case that $m\ge n$, however, it can be easily extended to the case of $m<n$. We focus in the former 
scenario since all experiments in Section \ref{sec:num} are generated by setting $m\ge n$.
\begin{algorithm}[H]
\floatname{algorithm}{Procedure}
\renewcommand{\thealgorithm}{}
\caption{Optimal solution Generator 2 (OsGen2)}
\begin{algorithmic}[1]
\STATE Given the required number $s\le \mbox{min}(m,n)$ of nonzeros in $x^*$, a positive constant $\gamma>0$, the right singular vectors $G$ and singular values $\Sigma$ of matrix $A$ do:
\STATE Solve approximately 
\begin{equation}\label{eq:130}
\begin{array}{cll}
x^*:= & \displaystyle\argmin_{x\in\mathbb{R}^{n}} & \|G^\intercal x - \gamma (\Sigma^\intercal \Sigma)^{-1} 1_n\|^2 \\
& \mbox{subject to:} & \|x\|_0 \le s, \\
\end{array}
\end{equation}
where $1_n\in\mathbb{R}^n$ 
is a vector of ones and $\|\cdot\|_0$ is the zero norm which returns the number of nonzero components of the input vector.
Problem \eqref{eq:130} can be solved approximately using an Orthogonal Matching Pursuit (OMP) \cite{cosamp} solver implemented in \cite{cosampImpl}.
\end{algorithmic}
\end{algorithm} 

The aim of Procedure OsGen2 is to find a sparse $x^*$ with $\kappa_\rho(x^*)$ arbitrarily large for some $\rho$ in the interval $\lambda_{n}(A^\intercal A) \le \rho \ll \lambda_{1}(A^\intercal A)$.
In particular, OsGen2 will return a sparse $x^*$ which can be expressed as $x^*=Gv$. The coefficients $v$ are close to the inverse of the eigenvalues of matrix $A^\intercal A$. 
Intuitively, this technique will create an $x^*$ which has strong dependence on subspaces which correspond to small eigenvalues of $A^\intercal A$.
The constant $\gamma$ is used in order to control the norm of $x^*$.  

The sparsity constraint in problem \eqref{eq:130}, i.e., $ \|x\|_0 \le s$, makes the approximate solution of this problem difficult when we use OMP, especially in the
case that $s$ and $n$ are large. To avoid this expensive task we can ignore the sparsity constraint in \eqref{eq:130}.
Then we can solve exactly and inexpensively the unconstrained problem and finally we can project the obtained solution in the feasible set defined by the sparsity constraint. 
Obviously, there is no guarantee that the projected solution is a good approximation to the one obtained in Step $2$ of Procedure OsGen2. However,
for all experiments in Section \ref{sec:num} that we applied this modification we obtained sufficiently large $\kappa_\rho(x^*)$. This means that 
our objective to produce ill-conditioned optimal solutions was met, while we kept the computational costs low. 
The modified version of Procedure OsGen2 is given in Procedure OsGen3.
\begin{algorithm}[H]
\floatname{algorithm}{Procedure}
\renewcommand{\thealgorithm}{}
\caption{Optimal solution Generator 3 (OsGen3)}
\begin{algorithmic}[1]
\STATE Given the required number $s\le \mbox{min}(m,n)$ of nonzeros in $x^*$, two non negative integers $s_1$ and $s_2$ such that $s_1+s_2 = s$, a positive constant $\gamma>0$, the right singular vectors $G$ and singular values $\Sigma$ of matrix $A$ do:
\STATE Solve exactly 
\begin{equation} \label{eq:131}
\begin{array}{cll}
x^*:= & \displaystyle\argmin_{x\in\mathbb{R}^{n}} & \|G^\intercal x - \gamma (\Sigma^\intercal \Sigma)^{-1} 1_n\|^2
\end{array}
\end{equation}
where $1_n\in\mathbb{R}^n$ 
is a vector of ones.
Problem \eqref{eq:131} can be solved exactly and inexpensively because $G^\intercal$ is an orthonormal matrix.
\STATE Maintain the positions and the values of the $s_1$ smallest and $s_2$ largest (in absolute values) components of ${x}^*$.
\STATE Set the remaining components of $x^*$ to zero. 
\end{algorithmic}
\end{algorithm} 

\section{Existing Problem Generators}\label{sec:existinggen}
So far in Section \ref{subsec:mgn} we have described in details our proposed problem generator. Moreover, in Section \ref{subsec:conA} we have described
how to construct matrices $A$ such that the proposed generator is scalable with respect to the number of unknown variables.  We now briefly describe 
existing problem generators and explain how our propositions add value to the existing approaches. 

Given a regularization parameter $\tau$ existing problem generators are looking for $A$, $b$ and $x^*$
such that the optimality conditions of problem \eqref{prob1}:
\begin{equation}\label{eq:system}
A^\intercal(Ax^* - b) \in -\tau \partial \|x^*\|_1,
\end{equation}
are satisfied. For example, 
in \cite{nesterovgen} the author fixes a vector of noise $e$ and an 
optimal solution $x^*$ and then finds $A$ and $b$ such that \eqref{eq:system} is satisfied. 
In particular, in \cite{nesterovgen} matrix $A=BD$ is used, where $B$ is a fixed matrix and $D$
is a scaling matrix such that the following holds.
\begin{equation*}
(BD)^\intercal e \in \tau \partial \|x^*\|_1,
\end{equation*}
Matrix $D$ is trivial to calculate, see Section $6$ in \cite{nesterovgen} for details. Then by setting $b = Ax^* + e$ \eqref{eq:system} is satisfied.
The advantage of this generator is that it allows control of the noise vector $e$, in comparison to our approach where the vector noise $e$ has 
to be determined by solving a linear system. On the other hand, one does not have direct control over the singular value decomposition of 
matrix $A$, since this depends on matrix $D$, which is determined based on the fixed vectors $e$ and $x^*$.

Another representative example is proposed in \cite{dirkLorenz}. This generator, which we discovered during the revision of our paper,
proposes the same setting as in our paper. In particular, given $A$, $x^*$ and $\tau$ one can construct a vector $b$ (or a noise vector $e$) such that \eqref{eq:system} is satisfied.
However, in \cite{dirkLorenz} the author suggests that $b$ can be found using a simple iterative procedure. Depending on matrix $A$ and how ill-conditioned it is,
this procedure might be slow. In this paper, we suggest that one can rely on numerical linear algebra tools, such as Givens rotation, in order to inexpensively
construct $b$ (or a noise vector $e$) using straightforwardly scalable operations.
Additionally, we show in Section \eqref{sec:num} that a simple construction of matrix $A$ is sufficient to extensively 
test the performance of methods. 
 
\section{Numerical Experiments}\label{sec:num}
In this section we study the performance of state-of-the-art first- and second-order methods as the conditioning and the dimensions of the problem increase.
The scripts that reproduce the experiments in this section as well as the problem generators that are described in Section \ref{sec:gen}
can be downloaded from: {\url{http://www.maths.ed.ac.uk/ERGO/trillion/}.}

\subsection{State-of-the-art Methods}
A number of efficient first- \cite{fista,changHsiehLin,HsiehChang,petermartin,peterbigdata,shwartzTewari,tsengblkcoo,nesterovhuge,tsengyun,wrightaccel,tonglange} 
and second-order \cite{NIPS2012_4523,proximalNewtonNocedal,2ndpaperstrongly,IEEEhowto:Jacekmf,IEEEhowto:boyd,proximalNewtonKatya,thesisschmidt} methods have been developed for the solution of problem \eqref{prob1}. 
In this section we examine the performance of the following state-of-the-art methods. Notice that the first three methods FISTA, PSSgb and PCDM 
do not perform smoothing of the $\ell_1$-norm, while pdNCG does.

\begin{itemize}
\item FISTA (Fast Iterative Shrinkage-Thresholding Algorithm) \cite{fista} is an optimal first-order method for problem \eqref{prob1}, which adheres to
the structure of GFrame. At a point $x$, FISTA builds a convex function:
$$
Q_\tau(y;x) := \tau \|y\|_1 + \frac{1}{2}\|Ax-b\|^2 + (A^\intercal(Ax-b))^\intercal (y-x) + \frac{L}{2} \|y-x\|_2^2,
$$
where $L$ is an upper bound of $\lambda_{max}(A^\intercal A)$, and solves subproblem \eqref{bd3} exactly using 
shringkage-thresholding \cite{ista,sparseMRI}.
An efficient implementation of this algorithm can be found as part 
of TFOCS (Templates for First-Order Conic Solvers) package \cite{convexTemplates} under the name N$83$. In this implementation the 
parameter $L$ is calculated dynamically.
\item PCDM (Parallel Coordinate Descent Method) \cite{peterbigdata} is a randomized parallel coordinate descent method. 
The parallel updates are performed asynchronously and the coordinates to be updated are chosen uniformly at random. 
Let $\varpi$ be the number of processors that are employed by PCDM. Then, at a point $x$, PCDM builds $\varpi$ convex 
approximations:
 $$
Q_\tau^i(y_i;x) := \tau |y_i| + \frac{1}{2}\|Ax-b\|^2 + (A_i^\intercal(Ax-b)) (y_i-x_i) + \frac{\beta L_i}{2} (y_i-x_i)^2,
$$
$\forall i=1,2,\cdots,\varpi$, where $A_i$ is the $i$th column of matrix $A$ and $L_i=(A^\intercal A)_{ii}$ is the $i$th 
diagonal element of matrix $A^\intercal A$ and $\beta$ is a positive constant which is defined in Subsection \ref{sec:paramtun}. 
The $Q_\tau^i$ functions are minimized exactly using shrinkage-thresholding. 
%
\item PSSgb (Projected Scaled Subgradient, Gafni-Bertsekas variant) \cite{thesisschmidt} is a  
second-order method. At each iteration of PSSgb the coordinates 
are separated into two sets, the working set $\mathcal{W}$ and the active set $\mathcal{A}$. The working 
set consists of all coordinates for which, the current point $x$ is nonzero. The active set is the complement of 
the working set $\mathcal{W}$. The following local quadratic model is build at each iteration 
$$
Q_\tau(y;x) := f_\tau(x) + (\tilde{\nabla} f_\tau(x))^\intercal(y-x) + \frac{1}{2}(y-x)^\intercal H (y-x),
$$
where $\tilde{\nabla} f_\tau(x)$ is a sub-gradient of $f_\tau$ at point $x$ with the minimum Euclidean norm, see Subsection {2.2.1} in \cite{thesisschmidt}
for details. Moreover, matrix $H$ 
is defined as:
\begin{equation*}
H = 
\begin{bmatrix}   H_{\mathcal{W}}    & 0 \\
                           0                            & H_{\mathcal{A}} 
       \end{bmatrix},
\end{equation*}
where $H_{\mathcal{W}}$ is an L-BFGS (Limited-memory Broyden-Fletcher-Goldfarb-Shanno) Hessian approximation with respect to the coordinates $\mathcal{W}$
and $H_{\mathcal{A}}$ is a positive diagonal matrix. The diagonal matrix $H_{\mathcal{A}}$ is a scaled identity matrix, where the Shanno-Phua/Barzilai-Borwein scaling
is used, see Subsection 2.3.1 in \cite{thesisschmidt} for details. 
The local model is minimized exactly since the inverse of matrix $H$ is known due
to properties of the L-BFGS Hessian approximation $H_{\mathcal{W}}$.
\item pdNCG (primal-dual Newton Conjugate Gradients) \cite{2ndpaperstrongly} is also a second-order method.
At every point $x$ pdNCG constructs a convex function $Q_\tau$ exactly as described for \eqref{bd4}. The 
subproblem \eqref{bd3} is solved \textit{inexactly} by reducing it to the linear system:
$$
\nabla^2 f_{\tau}^{\mu}(x) (y-x) = -\nabla f_{\tau}^{\mu}(x),
$$
which is solved approximately using preconditioned Conjugate Gradients (PCG).
A simple diagonal preconditioner is used for all experiments. The preconditioner is the inverse of the diagonal
of matrix ${\nabla}^2 f_\tau^\mu(x)$.
\end{itemize}


\subsection{Implementation details}
Solvers pdNCG, FISTA and PSSgb are implemented in MATLAB, while solver PCDM is a C++ implementation. 
We expect that the programming language will not be an obstacle for pdNCG, FISTA and PSSgb. 
This is because these methods rely only on basic linear algebra operations, such as the dot product, which are implemented in $C$++ in MATLAB by default. 
The experiments in Subsections \ref{subsec:increaseCond}, \ref{subsec:nontrivial}, \ref{subsec:incdim} 
were performed on a Dell PowerEdge R920 running Redhat Enterprise Linux with four Intel Xeon E7-4830 v2 2.2GHz processors, 20MB Cache, 7.2 GT/s QPI, Turbo (4x10Cores). 

The huge scale experiments in Subsection \ref{subsec:bigproblem} were performed on a Cray XC$30$ MPP supercomputer. 
This work made use of the resources provided by ARCHER (\url{http://www.archer.ac.uk/}), made available through the Edinburgh Compute and Data Facility (ECDF) 
(\url{http://www.ecdf.ed.ac.uk/}). According to the most recent list of commercial supercomputers, which is published in TOP$500$ list (\url{http://www.top500.org}), ARCHER is currently 
the $25th$ fastest supercomputer worldwide out of $500$ supercomputers. ARCHER has a total of $118, 080$ cores with performance $1,642.54$ TFlops/s on LINPACK benchmark and 
$2,550.53$ TFlops/s theoretical peak perfomance. The most computationally demanding experiments which are presented in Subsection \ref{subsec:bigproblem} required more than half of the cores 
of ARCHER, i.e., $65, 536$ cores out of $118, 080$.

\subsection{Parameter tuning}\label{sec:paramtun}

We describe the most important parameters for each solver, any other parameters are set to their default values. For pdNCG we set the smoothing parameter to $\mu = 10^{-5}$, this setting allows accurate solution 
of the original problem with an error of order $\mathcal{O}(\mu)$ \cite{IEEEhowto:Nesterov}. For pdNCG, PCG is terminated when the relative residual is less that $10^{-1}$
and the backtracking line-search is terminated if it exceeds $50$ iterations. Regarding FISTA the most important parameter is the calculation of the Lipschitz constant $L$, which is handled 
dynamically by TFOCS. For PCDM the coordinate Lipschitz constants $L_i$ $\forall i=1,2,\cdots,n$ are calculated exactly and parameter $\beta = 1 + (\omega-1)(\varpi-1)/(n-1)$,
where $\omega$ changes for every problem since it is the degree of partial separability of the fidelity function in \eqref{prob1}, which is easily calculated (see \cite{peterbigdata}), and  $\varpi=40$
is the number of cores that are used. For PSSgb we set the number of L-BFGS corrections to $10$.

We set the regularization parameter $\tau=1$, unless stated otherwise.
We run pdNCG for sufficient time such that the problems are adequately solved. Then, the rest of the methods are terminated when the objective function $ f_{\tau}$ in \eqref{prob1} 
is below the one obtained by pdNCG or when a predefined maximum number of iterations limit is reached. All comparisons are presented in figures which 
show the progress of the objective function against the wall clock time. This way the reader can compare the performance of the solvers for various levels of 
accuracy. We use logarithmic scales for the wall clock time and terminate runs which do not converge in about $10^5$ seconds, i.e., approximately $27$ hours.

\subsection{Increasing condition number of $A^\intercal A$}\label{subsec:increaseCond}
In this experiment we present the performance of FISTA, PCDM, PSSgb and pdNCG for increasing condition number of matrix $A^TA$ when
Procedure OsGen is used to construct the optimal solution $x^*$. We generate six matrices $A$ and two instances of $x^*$ for every matrix $A$; twelve instances in total.

The singular value decomposition of matrix $A$ is $A=\Sigma G^\intercal$, where $\Sigma$ is the matrix of singular values, the columns of matrices $I_m$ and $G$ are 
the left and right singular vectors, respectively,
see Subsection \ref{subsec:constA} for details about the construction of matrix $G$.
The singular values of matrices $A$ are chosen uniformly at random in the intervals $[0, 10^q]$, where $q=0,1,\cdots,5$, for each of the six matrices $A$.
Then, all singular values are shifted by $10^{-1}$.
The previous resulted in a condition number of matrix $A^\intercal A$ which varies from $10^2$ to $10^{12}$
with a step of times $10^2$. 
The rotation angle $\theta$ of matrix $G$ is set to $2{\pi}/{3}$ radians.
Matrices $A$ have $n=2^{22}$ columns, $m =2 n$ rows and rank $n$.
The optimal solutions $x^*$ have $s = n/2^7$ nonzero components for all twelve instances.

\begin{figure}[htbp]%
  \centering

  \subfloat[$\kappa(A^\intercal A) = 10^2$]{\label{fig1a}\includegraphics[scale=0.36]{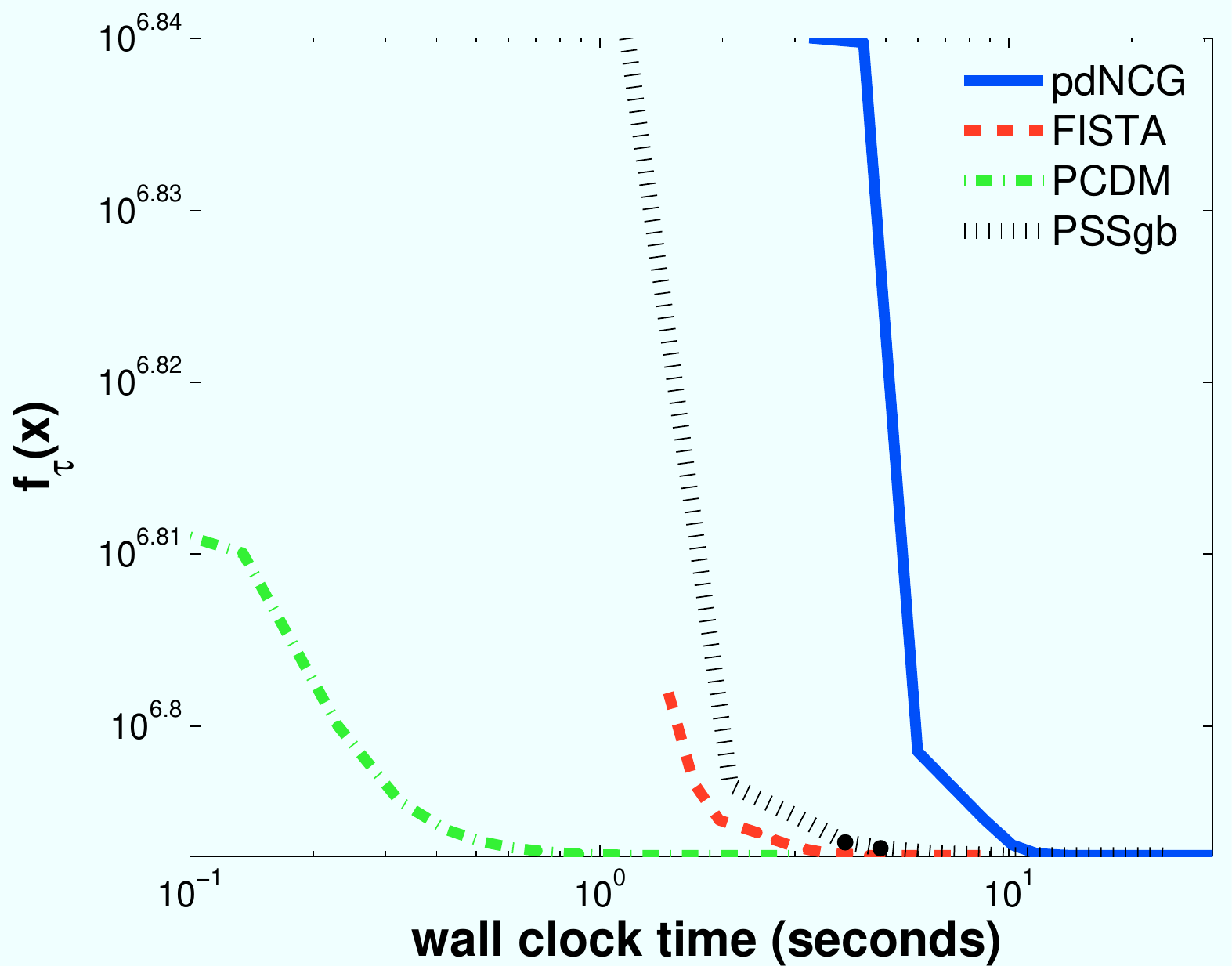}}
  \subfloat[$\kappa(A^\intercal A) = 10^4$]{\label{fig1b}\includegraphics[scale=0.36]{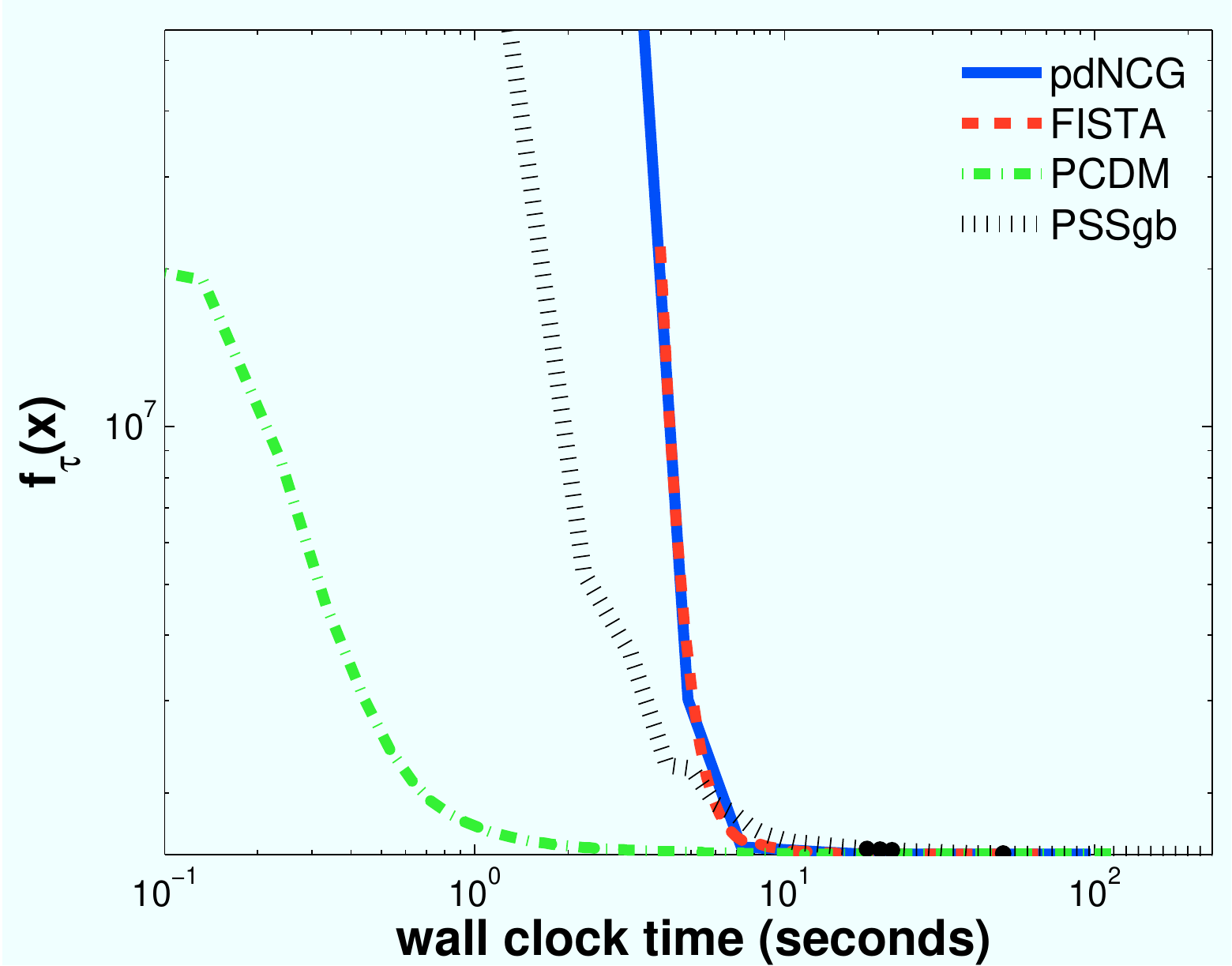}}
  \\
  \subfloat[$\kappa(A^\intercal A) = 10^6$]{\label{fig1c}\includegraphics[scale=0.36]{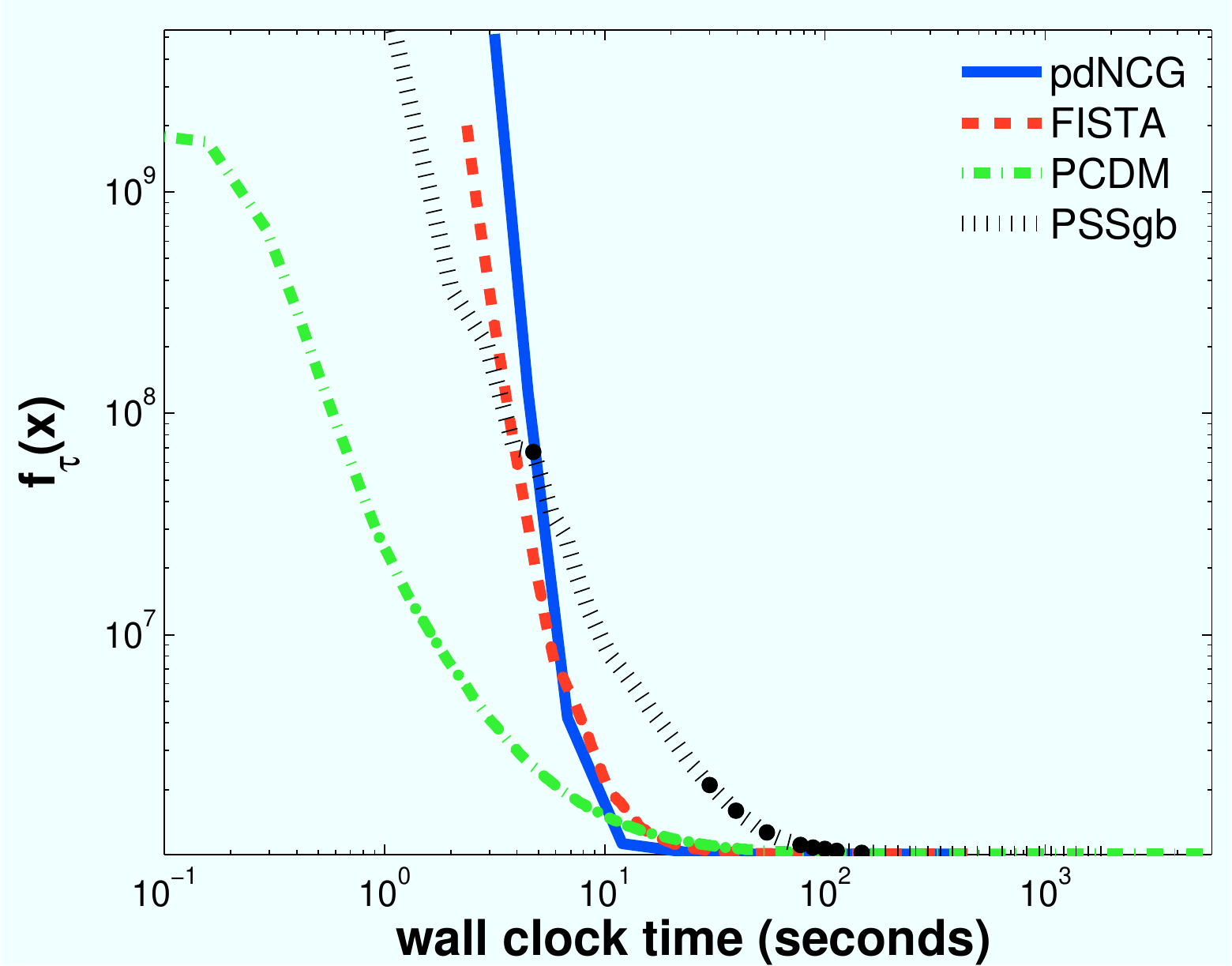}}
  \subfloat[$\kappa(A^\intercal A) = 10^8$]{\label{fig1d}\includegraphics[scale=0.36]{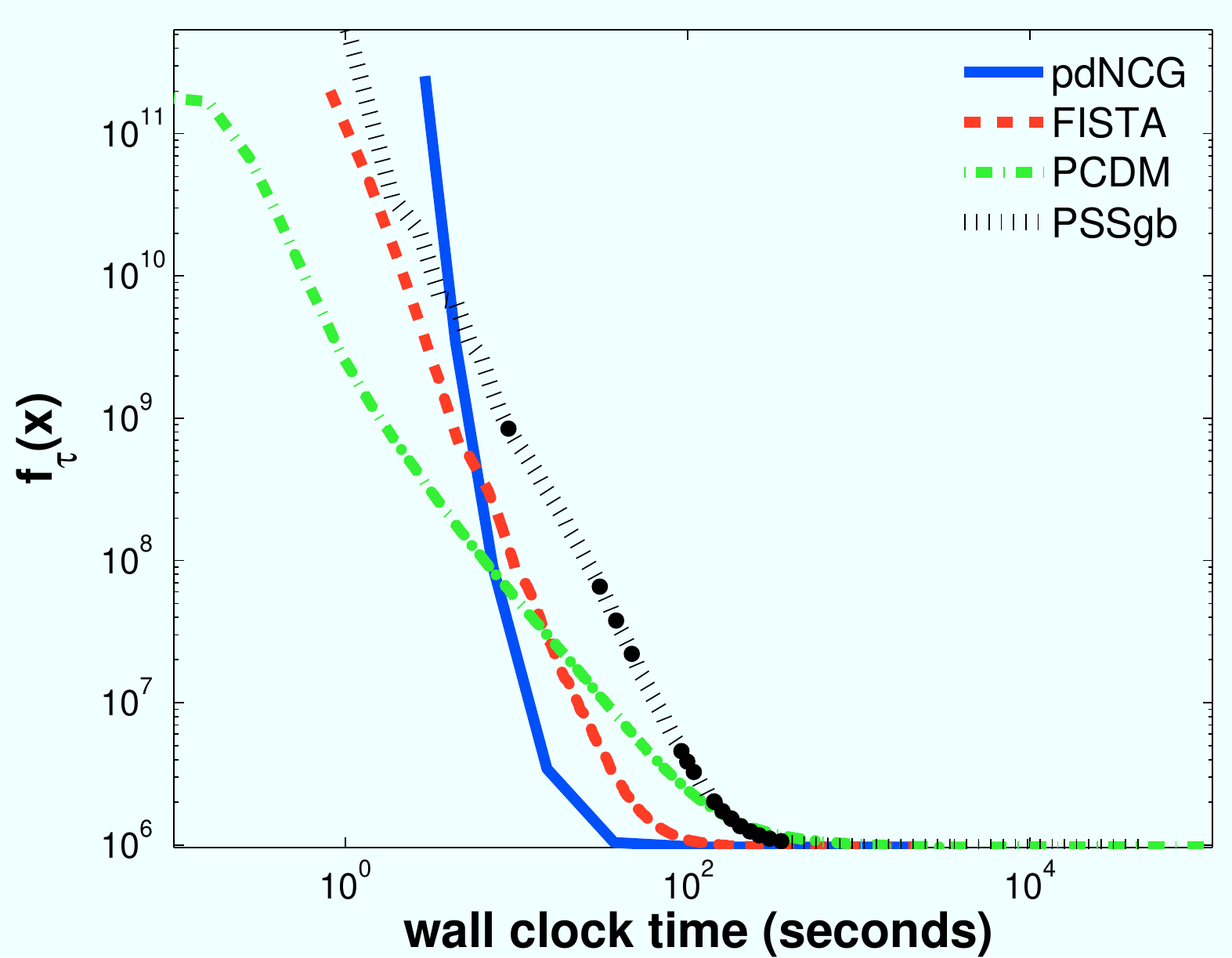}}
    \\
  \subfloat[$\kappa(A^\intercal A) = 10^{10}$]{\label{fig1e}\includegraphics[scale=0.36]{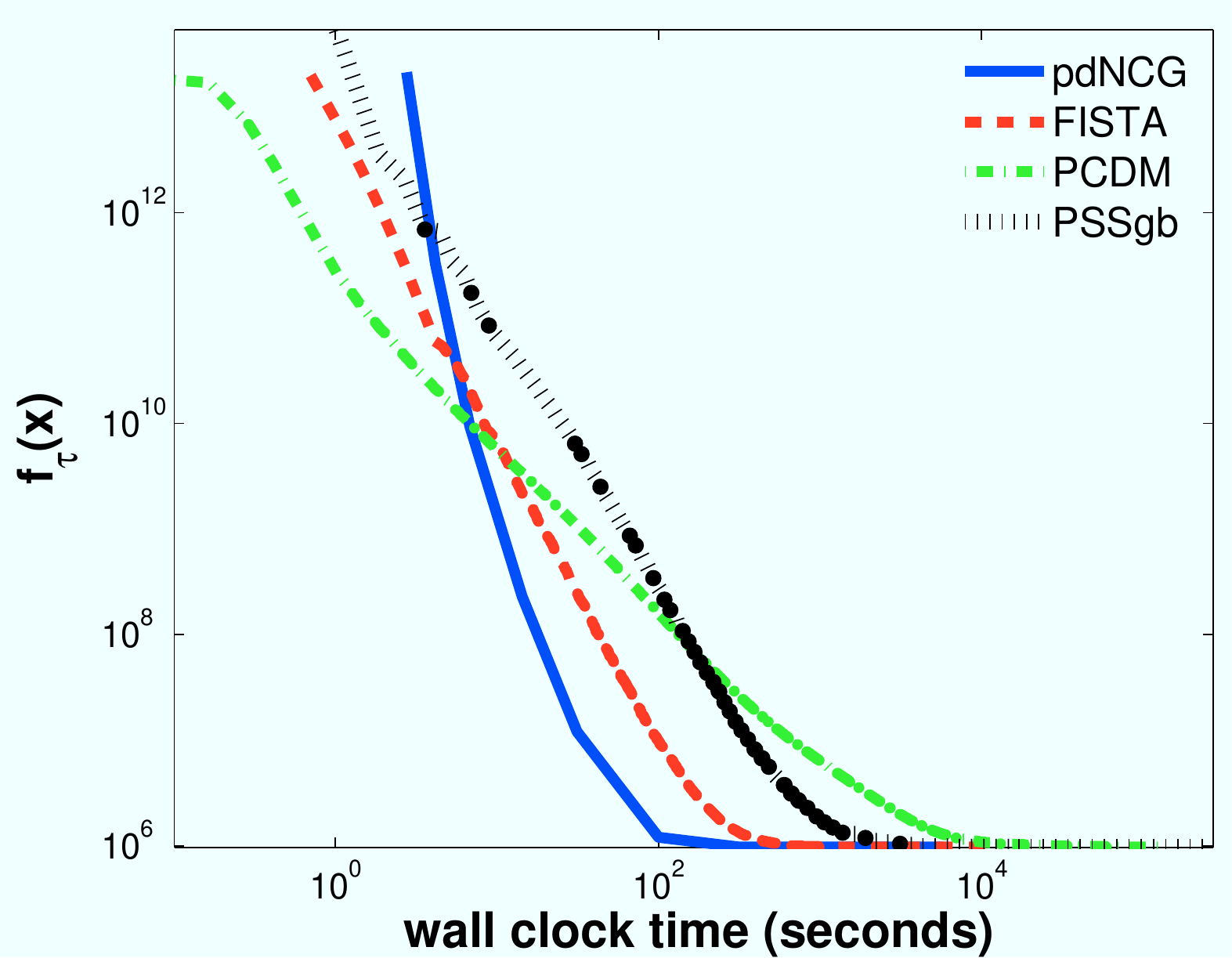}}
  \subfloat[$\kappa(A^\intercal A) = 10^{12}$]{\label{fig1f}\includegraphics[scale=0.36]{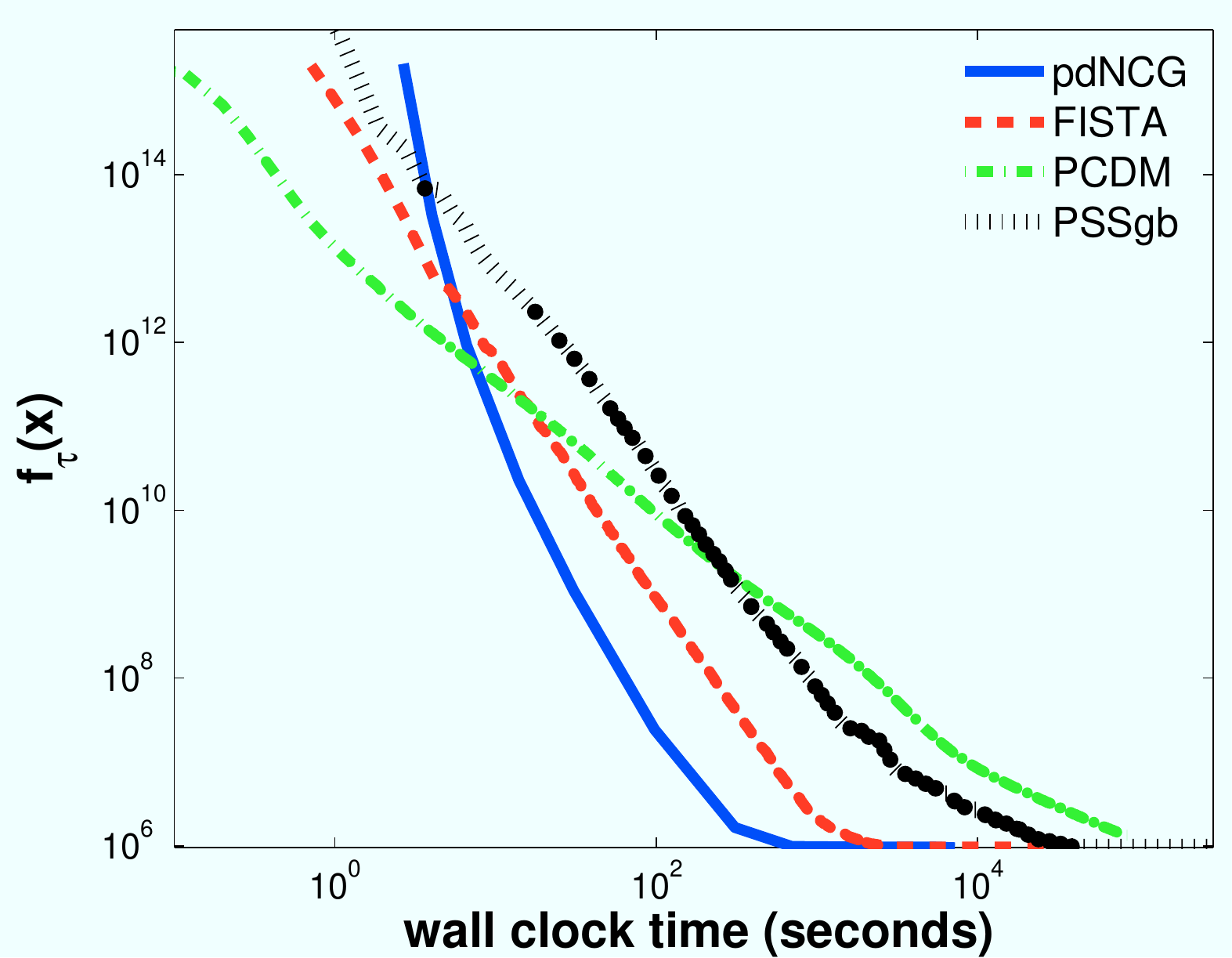}}
  
\caption{Performance of pdNCG, FISTA, PCDM and PSSgb on synthetic S-LS problems for increasing condition number of matrix $A^\intercal A$
and $\gamma=10$ in Procedure OsGen. 
The axes are in log-scale. In this figure $f_\tau$ denotes the objective value that was obtained by each solver.}
\label{fig1}%
\end{figure}

For the first set of six instances we set $\gamma=10$ in OsGen, which resulted in $\kappa_{0.1}(x^*) \approx 1$ for all experiments.
The results are presented in Figure \ref{fig1}. For these instances PCDM is clearly the fastest for $\kappa(A^\intercal A) \le 10^4$, while for $\kappa(A^\intercal A) \ge 10^6$ 
pdNCG is the most efficient. 

\begin{figure}[htbp]%
  \centering

  \subfloat[$\kappa(A^\intercal A) = 10^2$]{\label{fig2a}\includegraphics[scale=0.36]{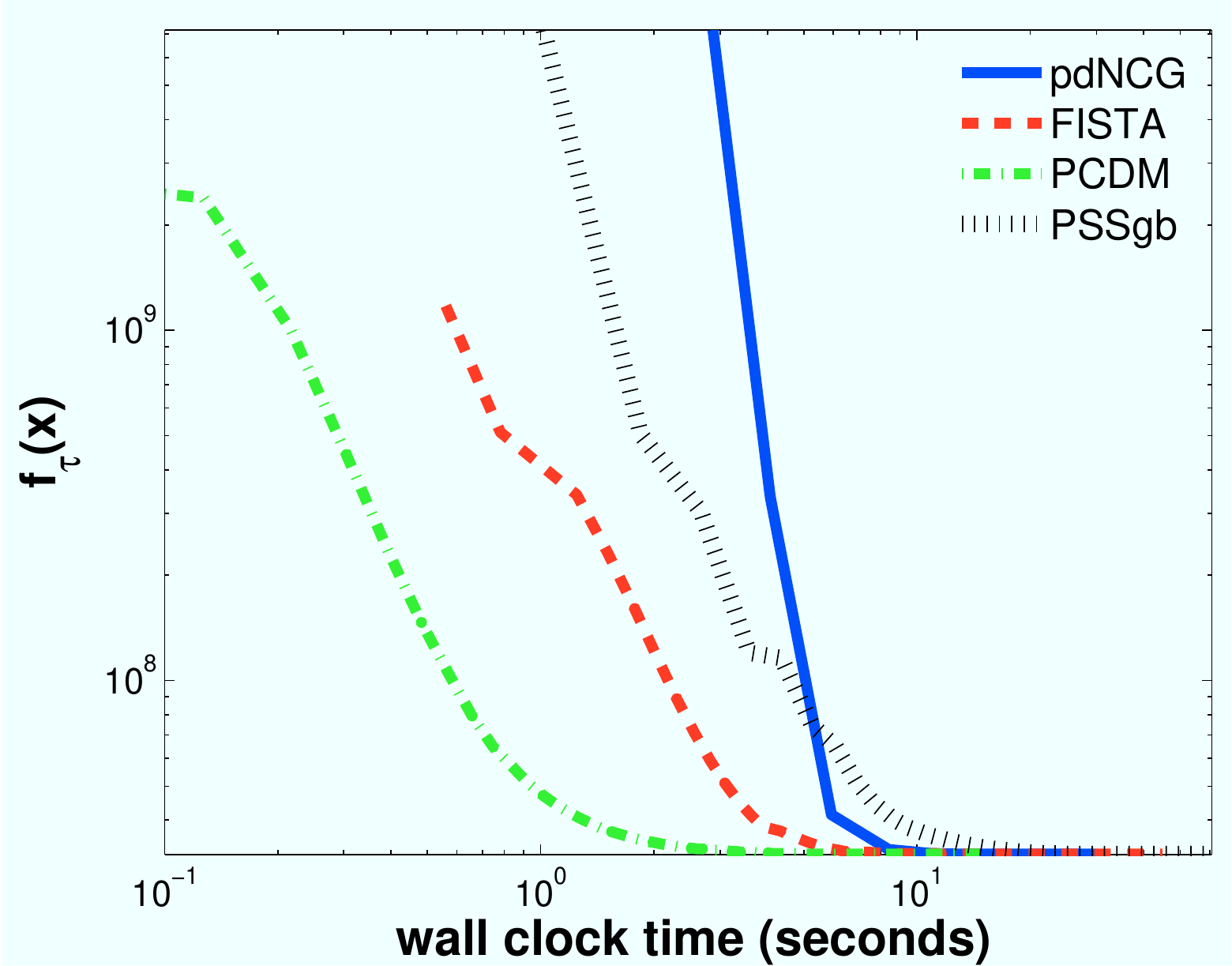}}
  \subfloat[$\kappa(A^\intercal A) = 10^4$]{\label{fig2b}\includegraphics[scale=0.36]{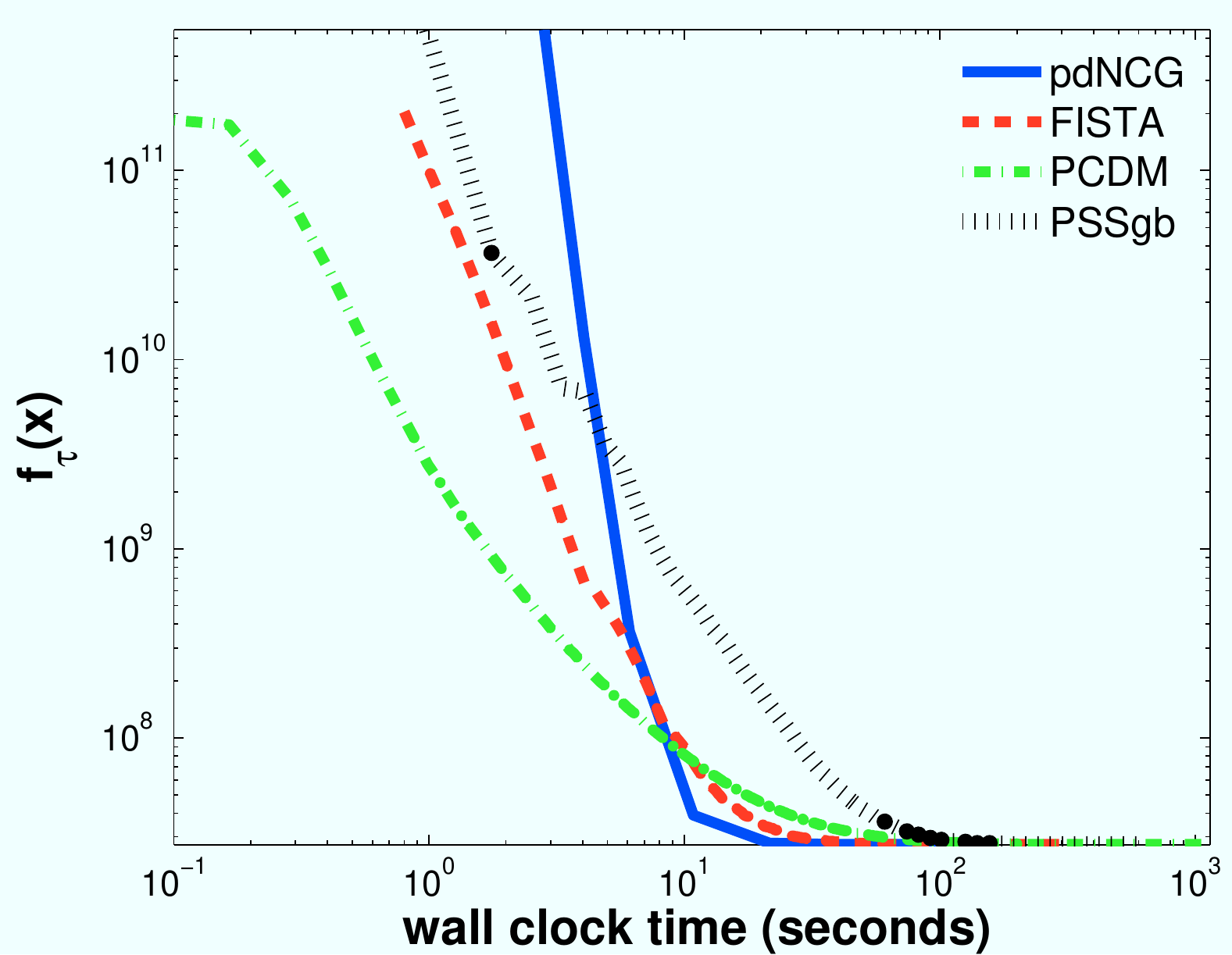}}
  \\
  \subfloat[$\kappa(A^\intercal A) = 10^6$]{\label{fig2c}\includegraphics[scale=0.36]{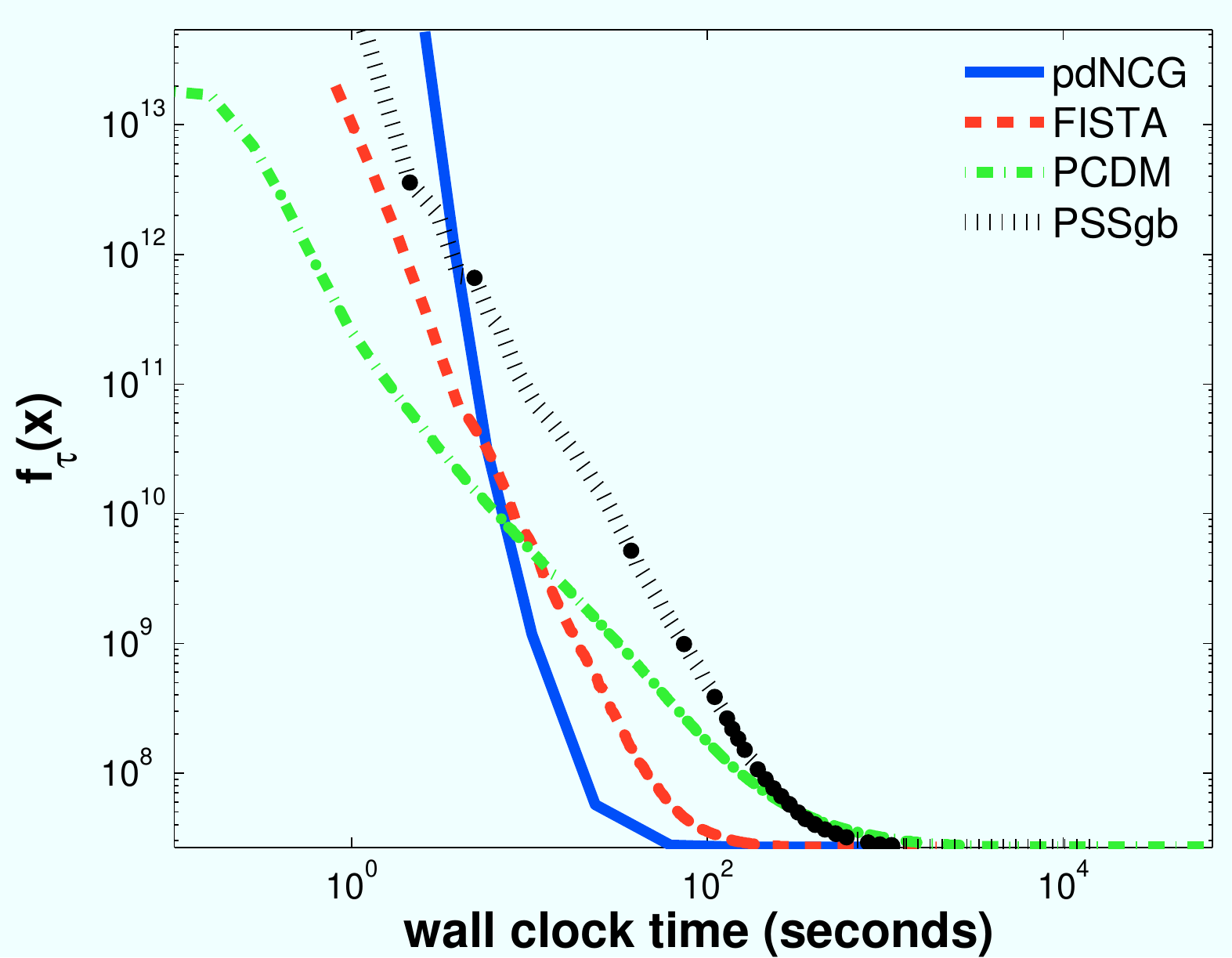}}
  \subfloat[$\kappa(A^\intercal A) = 10^8$]{\label{fig2d}\includegraphics[scale=0.36]{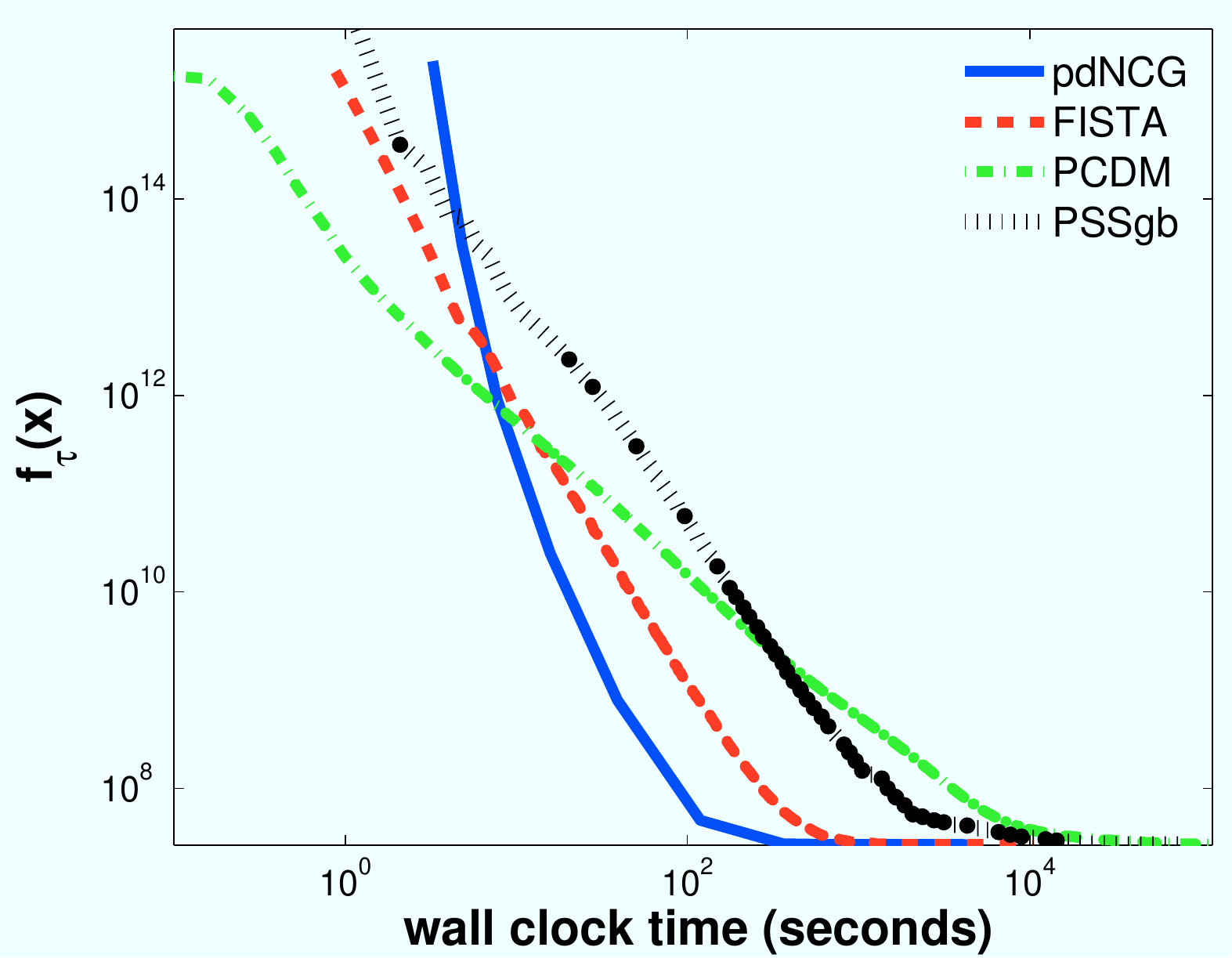}}
    \\
  \subfloat[$\kappa(A^\intercal A) = 10^{10}$]{\label{fig2e}\includegraphics[scale=0.36]{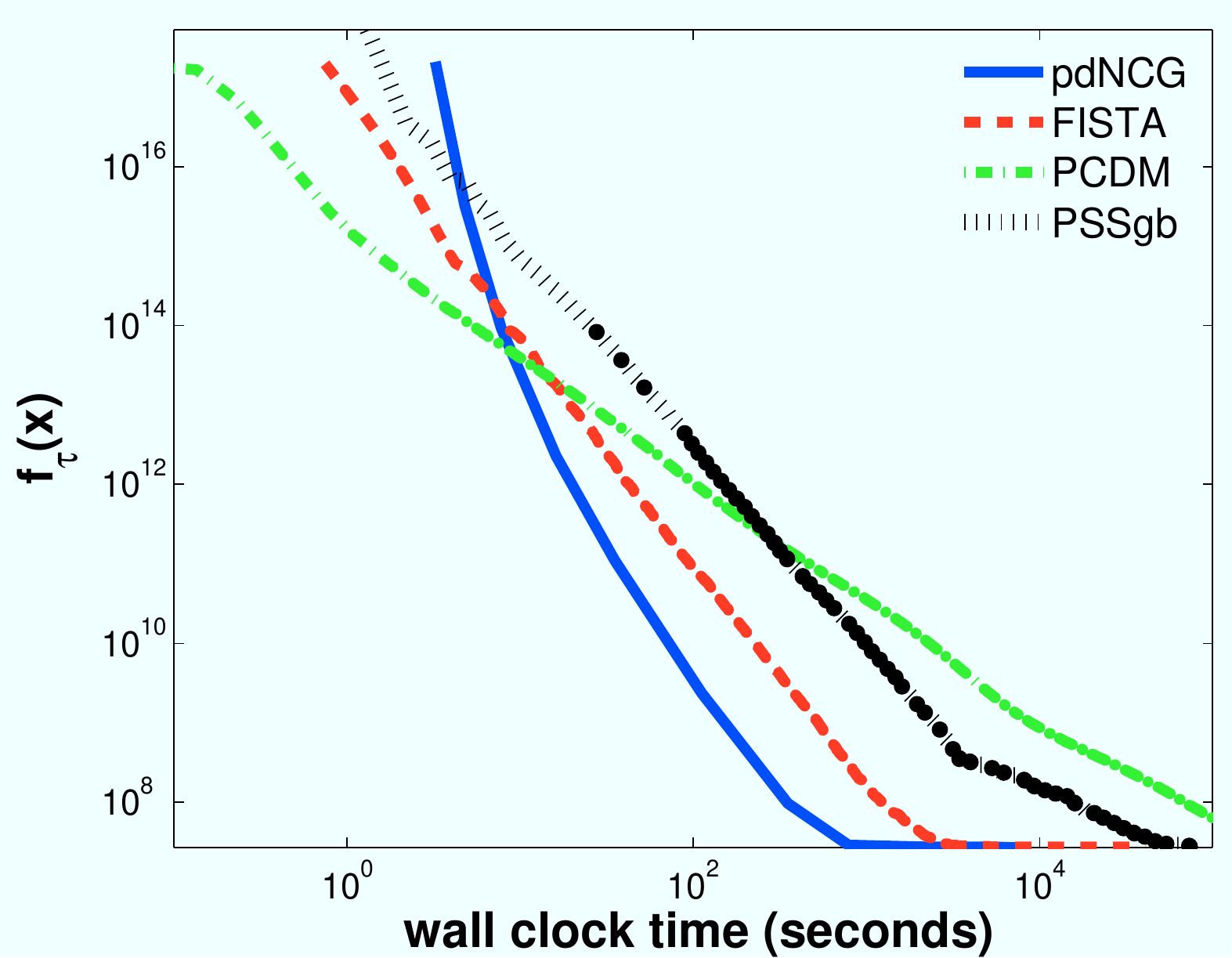}}
  \subfloat[$\kappa(A^\intercal A) = 10^{12}$]{\label{fig2f}\includegraphics[scale=0.36]{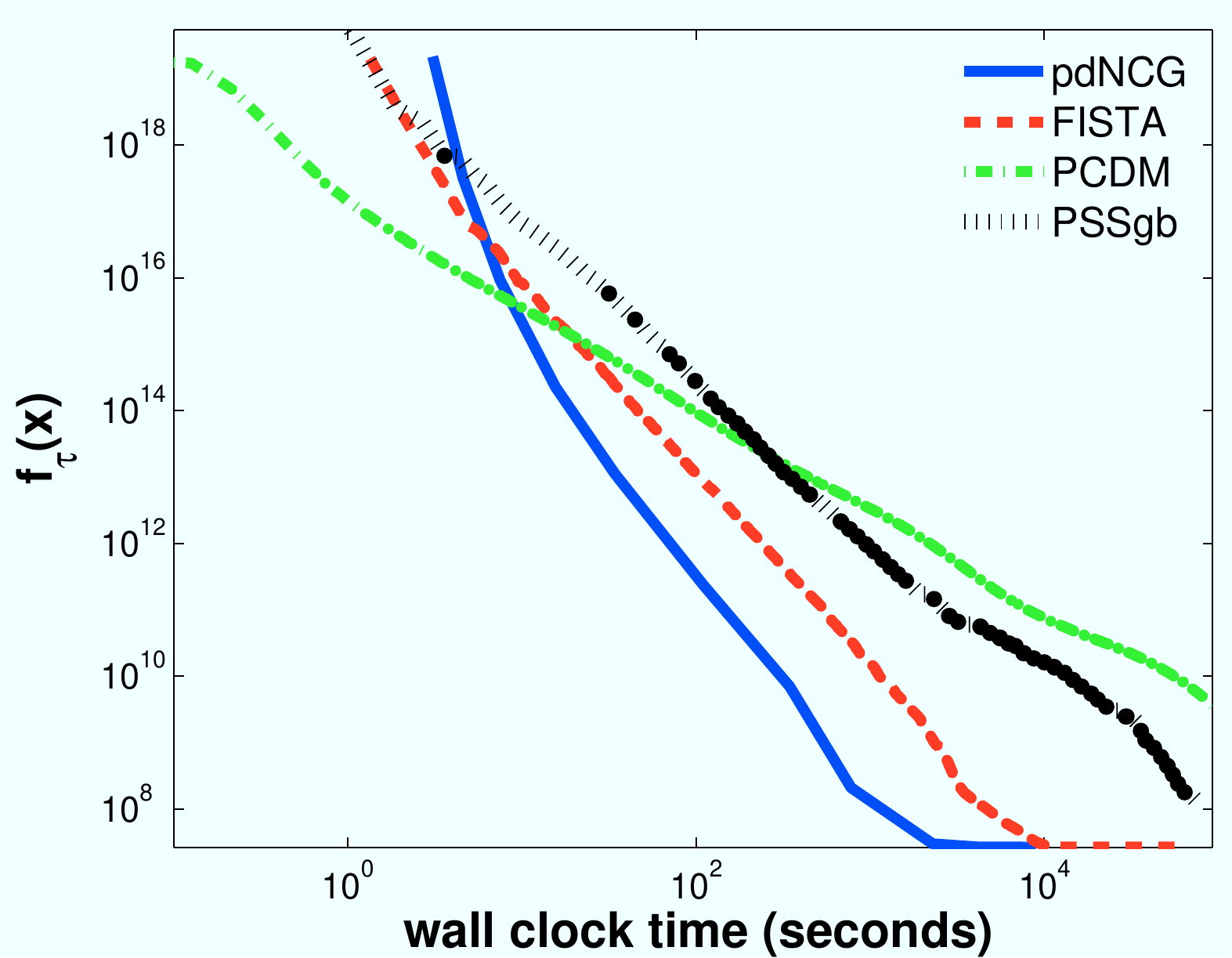}}
  
\caption{Performance of pdNCG, FISTA, PCDM and PSSgb on a synthetic S-LS problem for increasing condition number of matrix $A^\intercal A$
and $\gamma=10^3$ in Procedure OsGen. 
The axes are in log-scale.}
\label{fig2}%
\end{figure}

For the second set of six instances we set $\gamma=10^3$ in Procedure OsGen, which resulted in the same $\kappa_{0.1}(x^*)$ as before for every matrix $A$.
The results are presented in Figure \ref{fig2}. For these instances PCDM is the fastest for very well conditioned problems with $\kappa(A^\intercal A) \le 10^2$, while pdNCG is the fastest for $\kappa(A^\intercal A) \ge 10^4$.

We observed that
pdNCG required at most 30 iterations to converge for all experiments. For FISTA, PCDM and PSSgb the number of iterations 
was varying between thousands and tens of thousands iterations depending on the condition number of matrix $A^\intercal A$; the larger the condition 
number the more the iterations. However, the number of iterations is not a fair metric to compare solvers because every solver has different computational 
cost per iteration. In particular, FISTA, PCDM and PSSgb perform few inner products per iteration, which makes every iteration inexpensive, but the number of 
iterations is sensitive to the condition number of matrix $A^\intercal A$.  
On the other hand, for pdNCG the empirical iteration complexity is fairly stable, however, the number of inner products per iteration (mainly matrix-vector products with matrix $A$) 
may increase as the condition
number of matrix $A^\intercal A$ increases. 
Inner products are the major computational burden at every iteration for all solvers, therefore, the faster 
an algorithm converged in terms of wall-clock time the less inner products that are calculated. 
In Figures \ref{fig1} and \ref{fig2}
we display the objective evaluation against wall-clock time (log-scale) to facilitate the comparison of different algorithms. 

%

\subsection{Increasing condition number of $A^\intercal A$: non-trivial construction of $x^*$}\label{subsec:nontrivial}
In this experiment we examine the performance of the methods as the condition number of matrix $A^\intercal A$ increases,
while the optimal solution $x^*$ is generated using Procedure OsGen3 (instead of OsGen) with $\gamma = 100$ and $s_1 = s_2= s/2$.
Two classes of instances are generated, each class consists of four instances $(A,x^*)$ with $n=2^{22}$, $m =2 n$ and $s=n/2^7$. Matrix $A$ is constructed as in Subsection \ref{subsec:increaseCond}.
The singular values of matrices $A$ are chosen uniformly at random in the intervals $[0, 10^q]$, where $q=0,1,\cdots,3$, for all generated matrices $A$.
Then, all singular values are shifted by $10^{-1}$.
The previous resulted in a condition number of matrix $A^\intercal A$ which varies from $10^2$ to $10^{8}$
with a step of times $10^2$. 
The condition number of the generated optimal solutions was on average $\kappa_{0.1}(x^*)\approx 40$.

The two classes of experiments are distinguished based on the rotation angle $\theta$ that is used for the composition of Givens rotations $G$.
In particular, for the first class of experiments the angle is $\theta=2{\pi}/10$ radians,
while for the second class of experiments the rotation angle is $\theta=2{\pi}/10^{3}$ radians. 
The difference between the two classes is that the second class consists of matrices $A^\intercal A$ for which, a major part of their mass
is concentrated in the diagonal. This setting is beneficial for PCDM since it uses information only from the diagonal of matrices $A^\intercal A$. 
This setting is also beneficial for pdNCG since it uses a diagonal preconditioner
for the inexact solution of linear systems at every iteration. 
 
\begin{figure}[htbp]%
  \centering

  \subfloat[$\kappa(A^\intercal A) = 10^2$]{\label{fig3_2a}\includegraphics[scale=0.36]{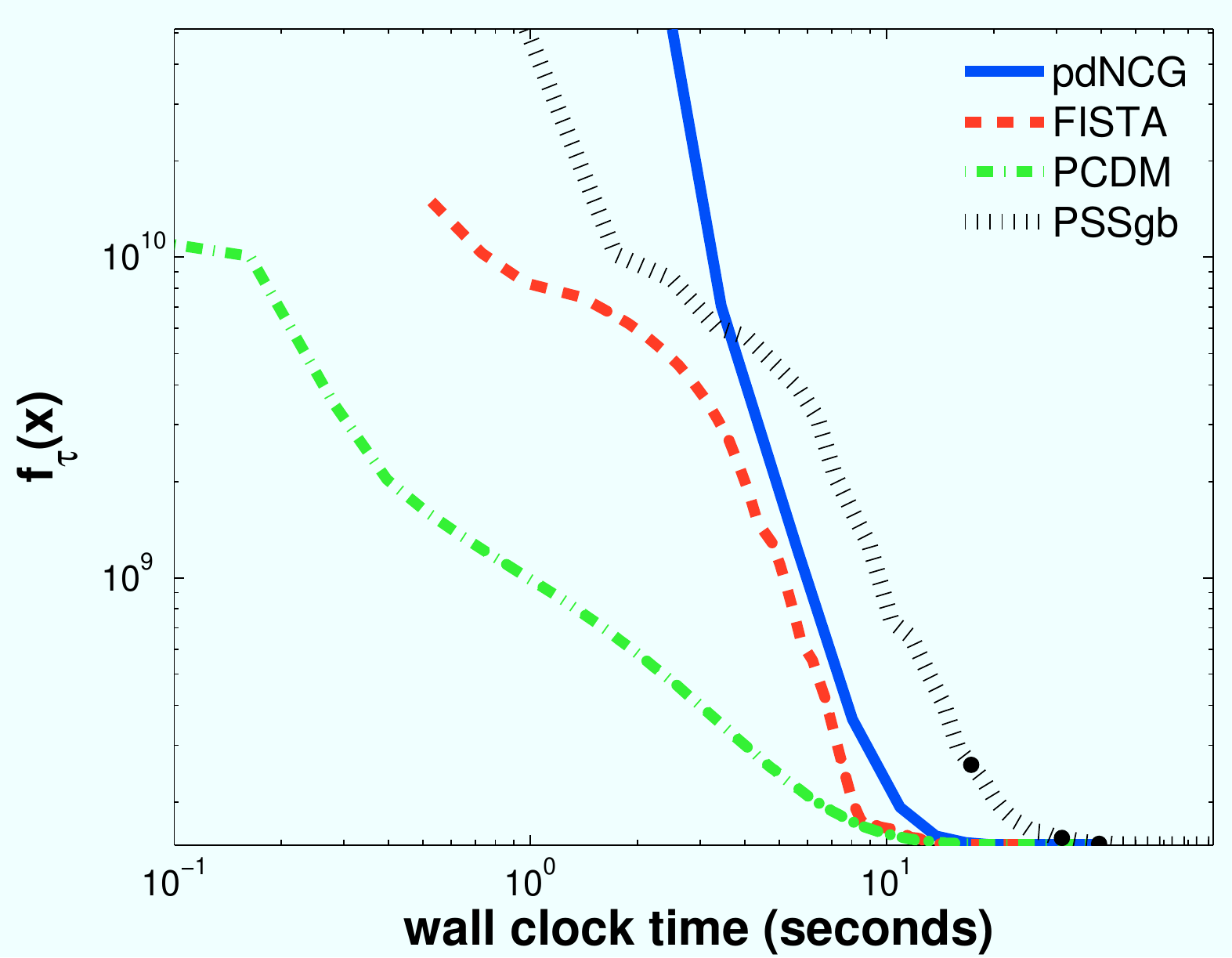}}
  \subfloat[$\kappa(A^\intercal A) = 10^4$]{\label{fig3_2b}\includegraphics[scale=0.36]{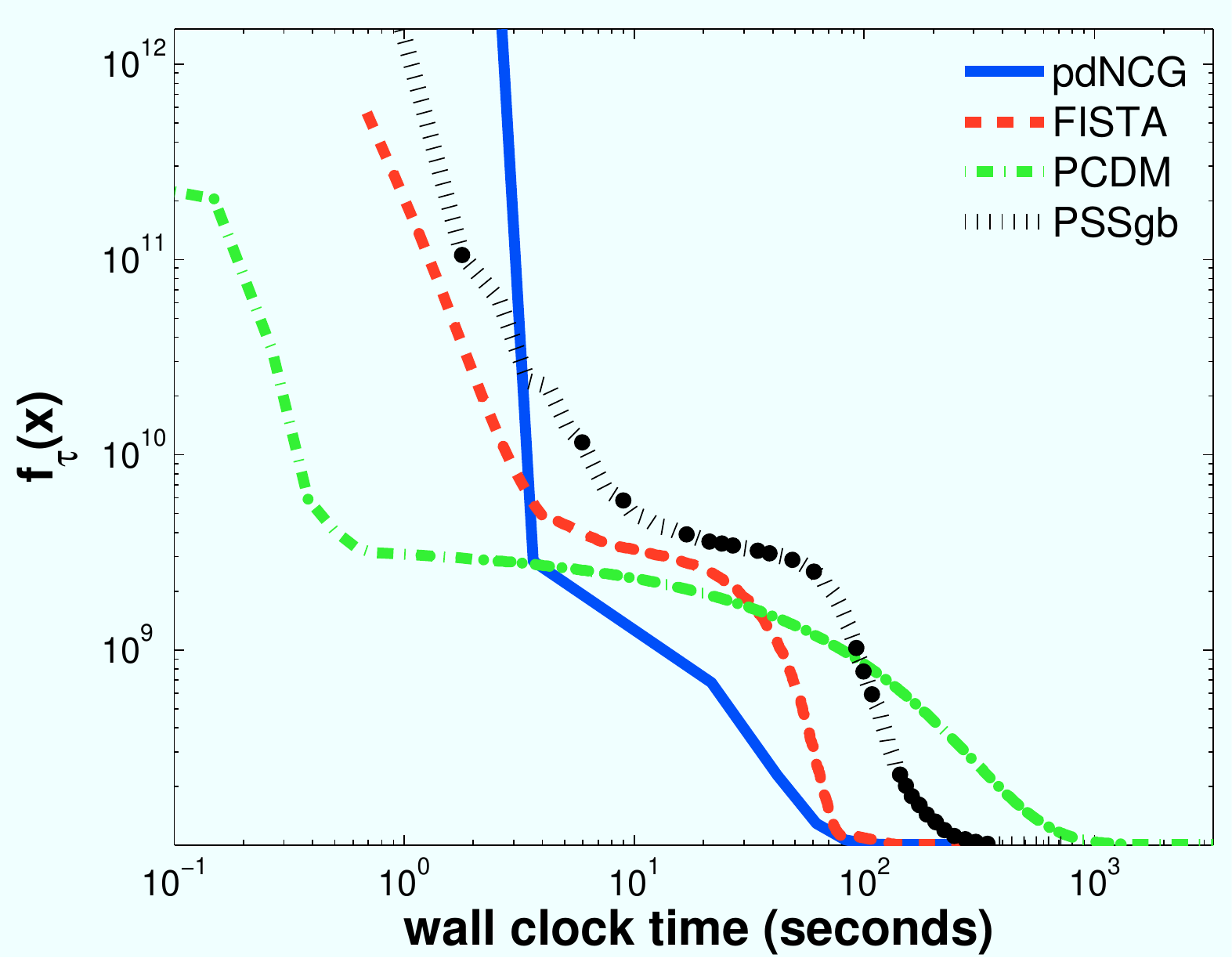}}
  \\
  \subfloat[$\kappa(A^\intercal A) = 10^6$]{\label{fig3_2c}\includegraphics[scale=0.36]{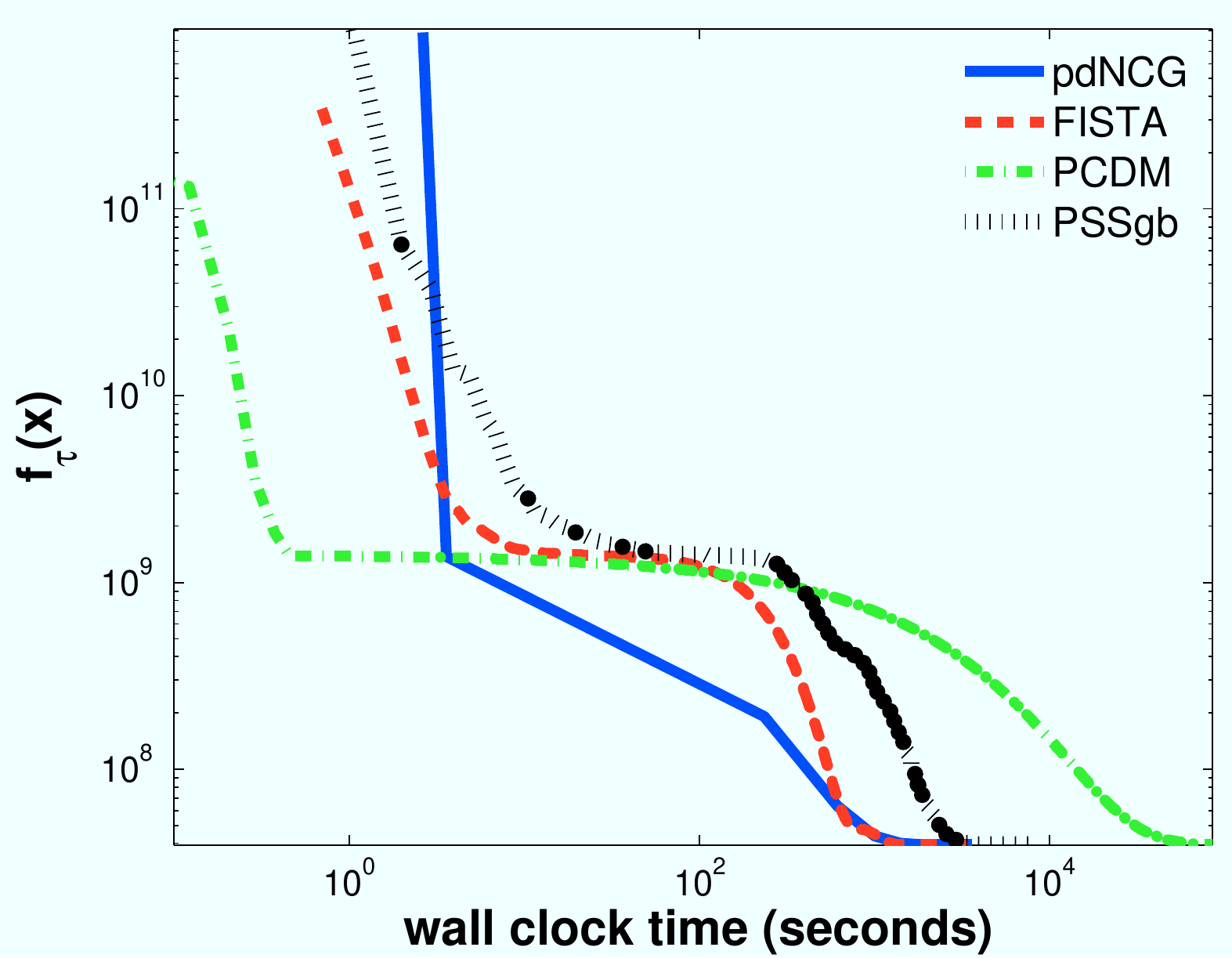}}
  \subfloat[$\kappa(A^\intercal A) = 10^8$]{\label{fig3_2d}\includegraphics[scale=0.36]{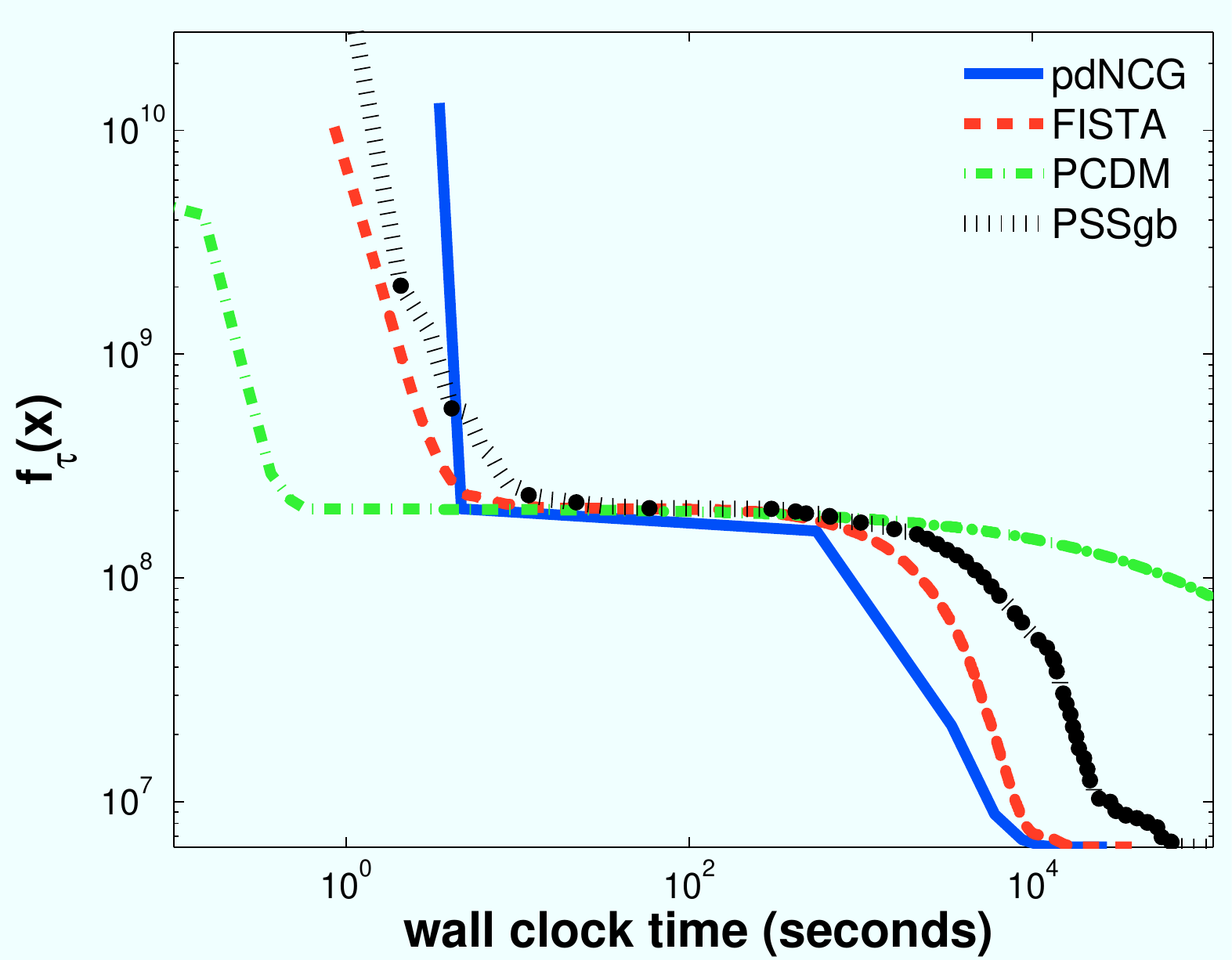}}
  
\caption{Performance of pdNCG, FISTA, PCDM and PSSgb on synthetic S-LS problems for increasing condition number of matrix $A^\intercal A$.
The optimal solutions have been generated using Procedure OsGen3 with $\gamma=100$ and $s_1 = s_2= s/2$.
The axes are in log-scale. The rotation angle $\theta$ in $G$ was $2\pi/10$. For condition number $\kappa(A^\intercal A) \ge 10^{6}$ PCDM was terminated after $1,000,000$ iterations, 
which corresponded to more than $27$ hours of wall-clock time.}
\label{fig3_2}%
\end{figure}

 The results for the first class of experiments are presented in Figure \ref{fig3_2}.
 For instances with $\kappa(A^\intercal A) \ge 10^{6}$ PCDM was terminated after $1,000,000$ iterations, 
which corresponded to more than $27$ hours of wall-clock time.

\begin{figure}[htbp]%
  \centering

  \subfloat[$\kappa(A^\intercal A) = 10^2$]{\label{fig3a}\includegraphics[scale=0.36]{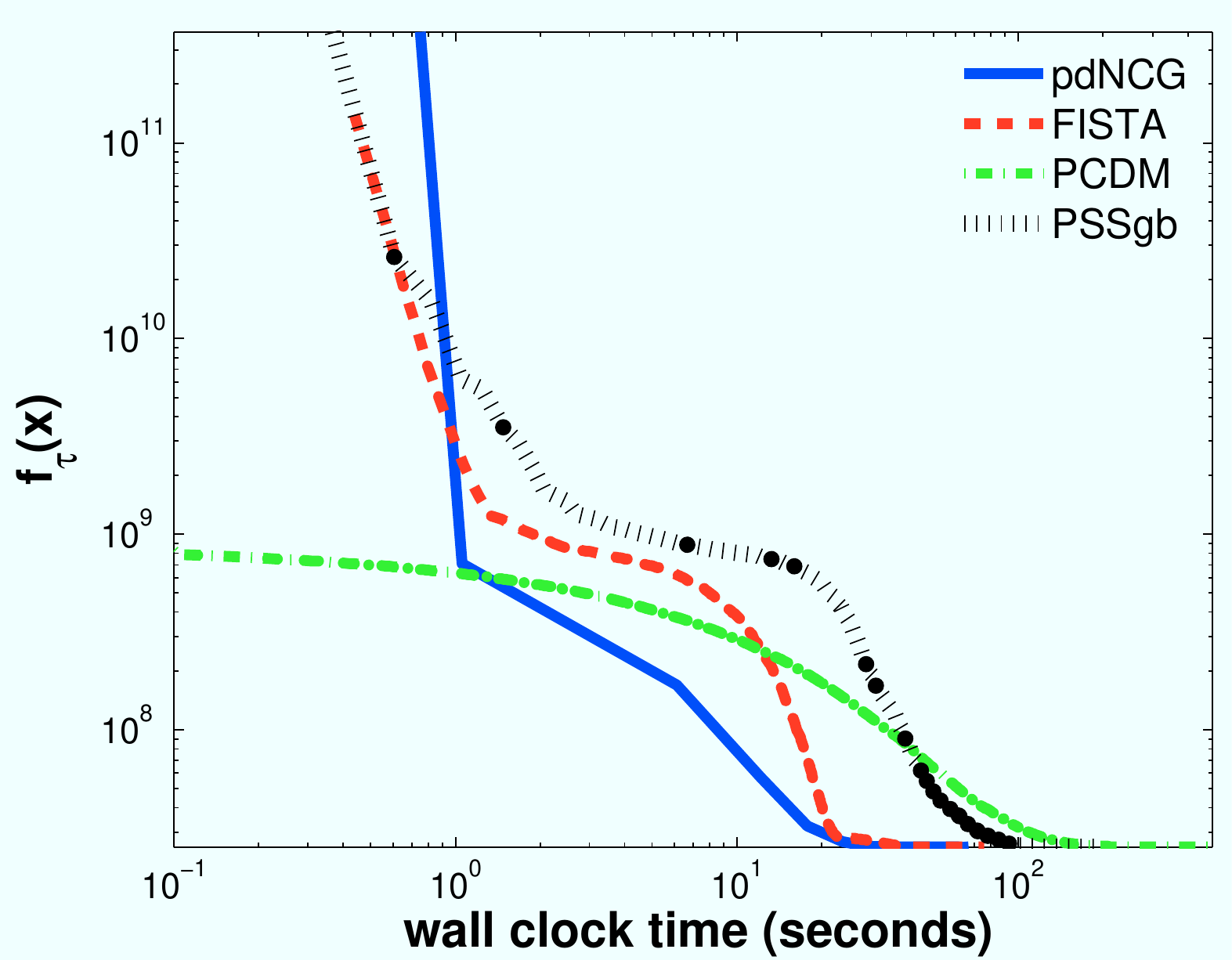}}
  \subfloat[$\kappa(A^\intercal A) = 10^4$]{\label{fig3b}\includegraphics[scale=0.36]{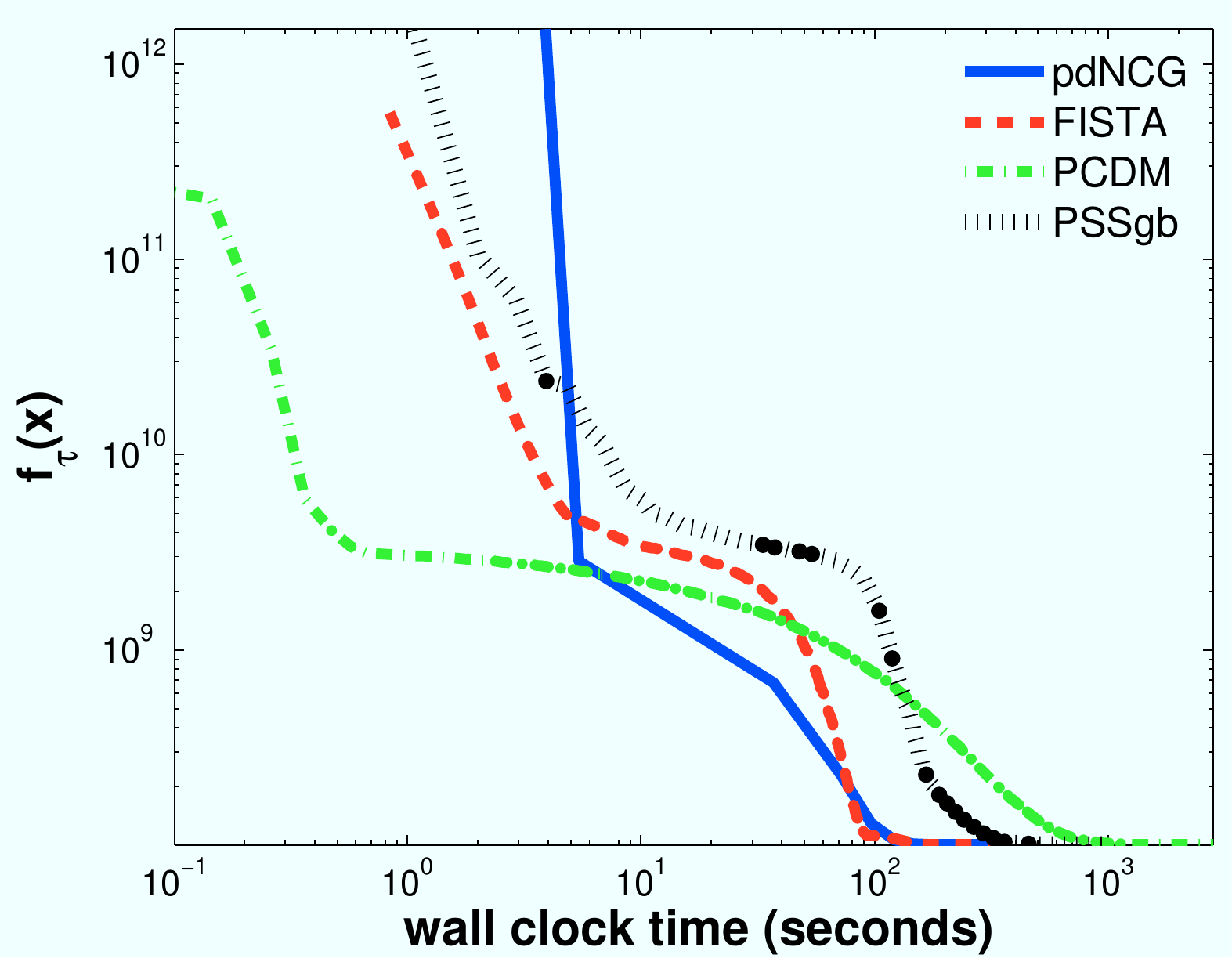}}
  \\
  \subfloat[$\kappa(A^\intercal A) = 10^6$]{\label{fig3c}\includegraphics[scale=0.36]{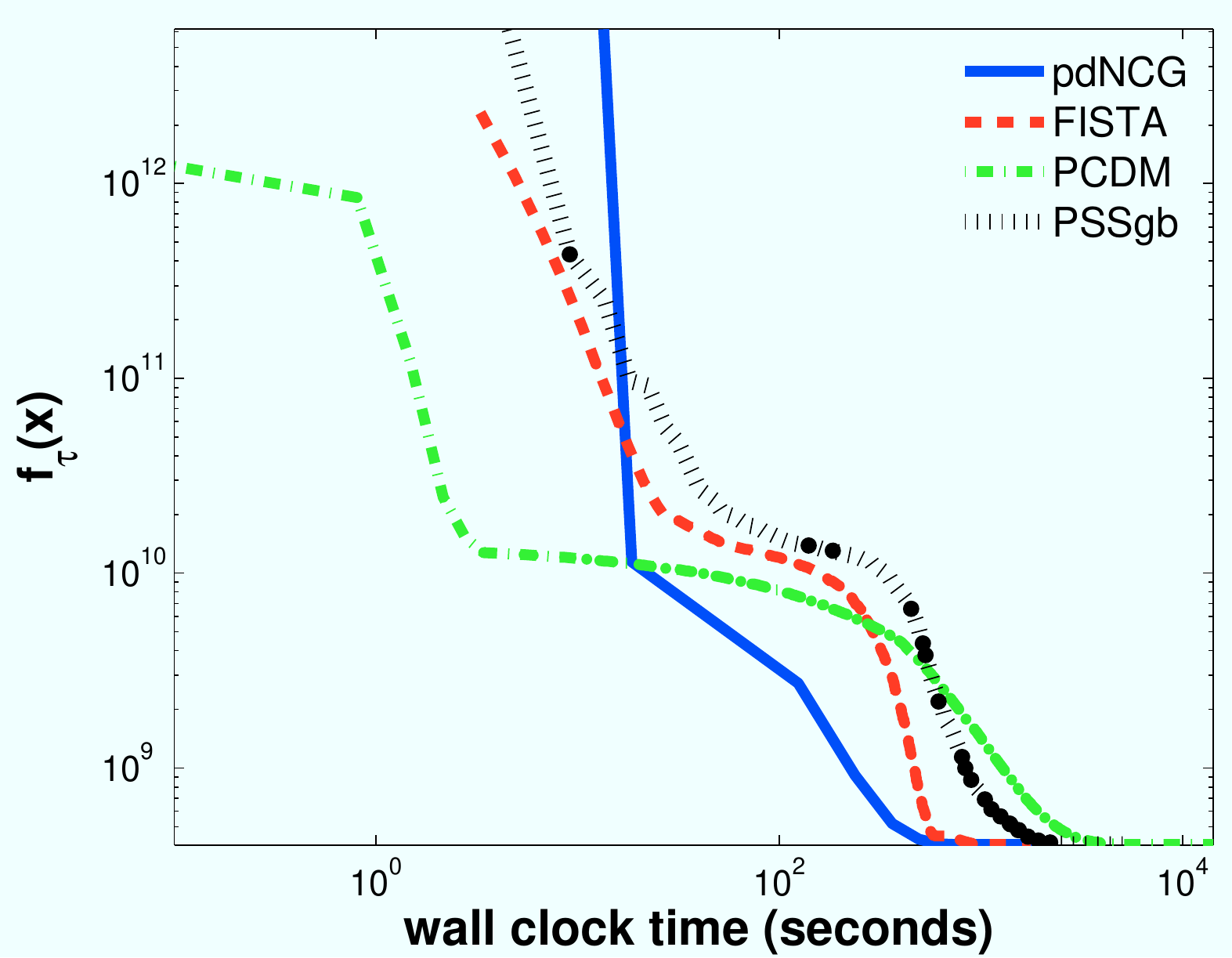}}
  \subfloat[$\kappa(A^\intercal A) = 10^8$]{\label{fig3d}\includegraphics[scale=0.36]{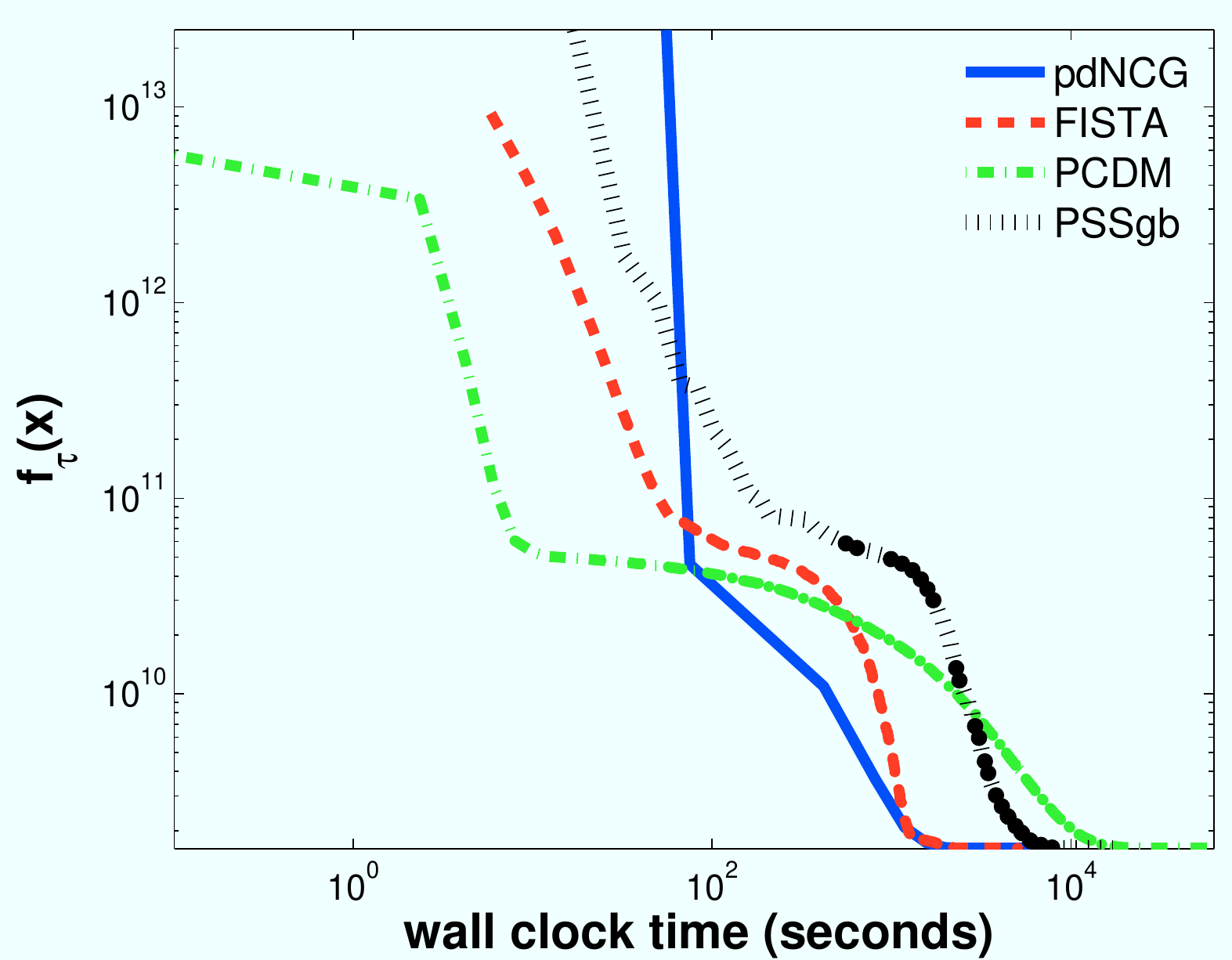}}
  
\caption{Performance of pdNCG, FISTA, PCDM and PSSgb on synthetic S-LS problems for increasing condition number of matrix $A^\intercal A$.
The optimal solutions have been generated by using Procedure OsGen3 with $\gamma=100$ and $s_1 = s_2= s/2$. The rotation angle $\theta$ in $G$ was $2\pi/10^{3}$.
The axes are in log-scale. }
\label{fig3}%
\end{figure}

The results for the second class of experiments are presented in Figure \ref{fig3}.
Notice in this figure that the objective function is only slightly reduced. This does not mean that the initial solution, which was the zero vector,
was nearly optimal. This is because 
noise with large norm, i.e., $\|Ax^* - b\|$ is large, was used in these experiments, therefore, changes in
the optimal solution did not have large affect on the objective function.  

\subsection{Increasing dimensions}\label{subsec:incdim}
In this experiment we present the performance of pdNCG, FISTA, PCDM and PSSgb as the number of variables $n$ increases. 
We generate four instances where the number of variables $n$ takes values $2^{20}$, $2^{22}$,
$2^{24}$ and $2^{26}$, respectively.
The singular value decomposition of matrix $A$ is $A=\Sigma G^\intercal$.
The singular values in matrix $\Sigma$ are chosen uniformly at random in the interval $[0, 10]$ and then are shifted by $10^{-1}$,
which resulted in $\kappa(A^\intercal A)\approx 10^{4}$.
The rotation angle $\theta$ of matrix $G$ is set to $2\pi/10$ radians.
Moreover, matrices $A$ have $m= 2 n$ rows and rank $n$.
The optimal solutions $x^*$ have $s = n/2^7$ nonzero components for each generated instance. 
For the construction of the optimal solutions $x^*$ we use Procedure OsGen3 with $\gamma = 100$ and $s_1 = s_2= s/2$, which resulted 
in $\kappa_{0.1}(x^*) \approx 3$ on average.

\begin{figure}[htbp]%
  \centering

  \subfloat[$n = 2^{20}$]{\label{fig4b}\includegraphics[scale=0.36]{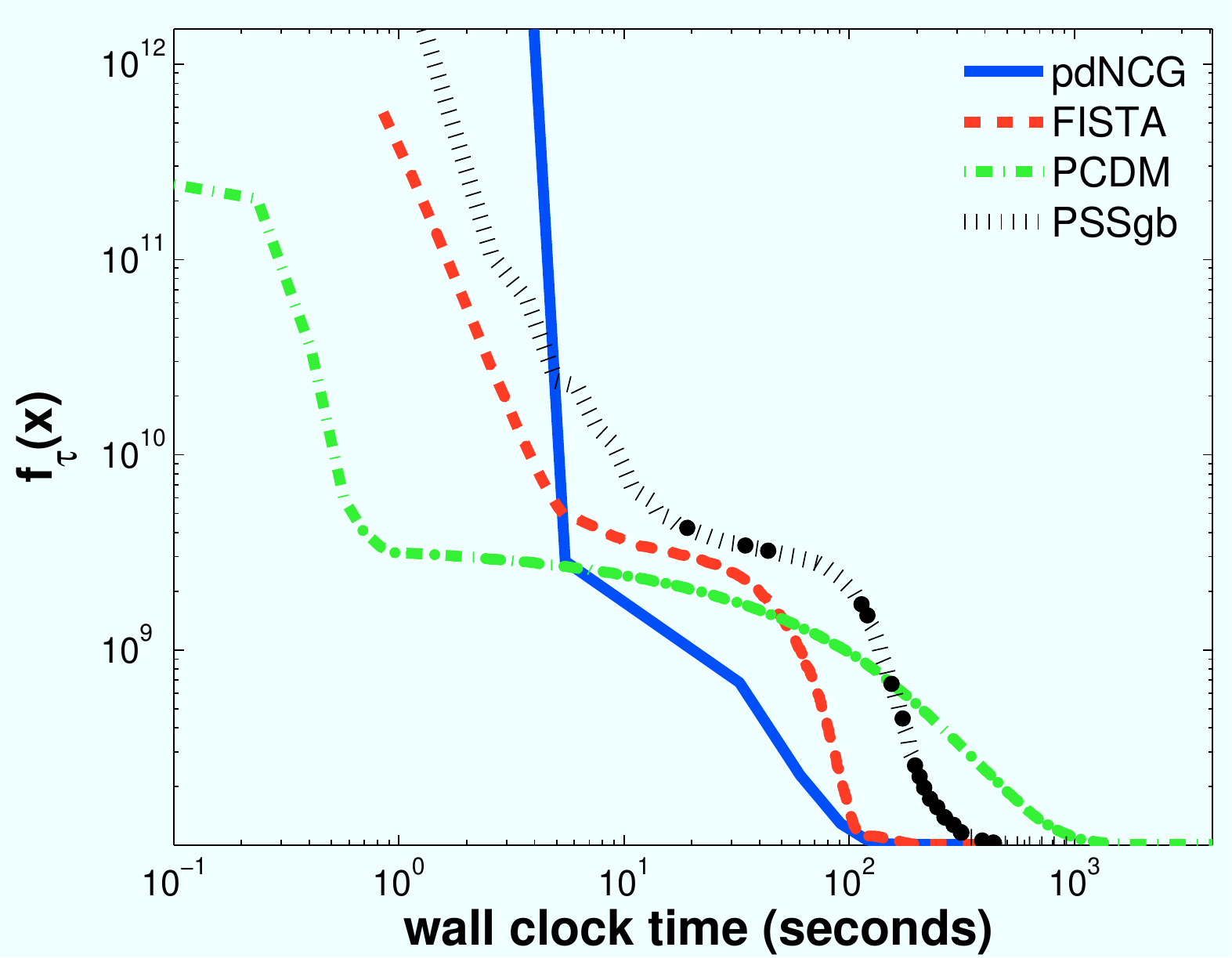}}
  \subfloat[$n = 2^{22}$]{\label{fig4d}\includegraphics[scale=0.36]{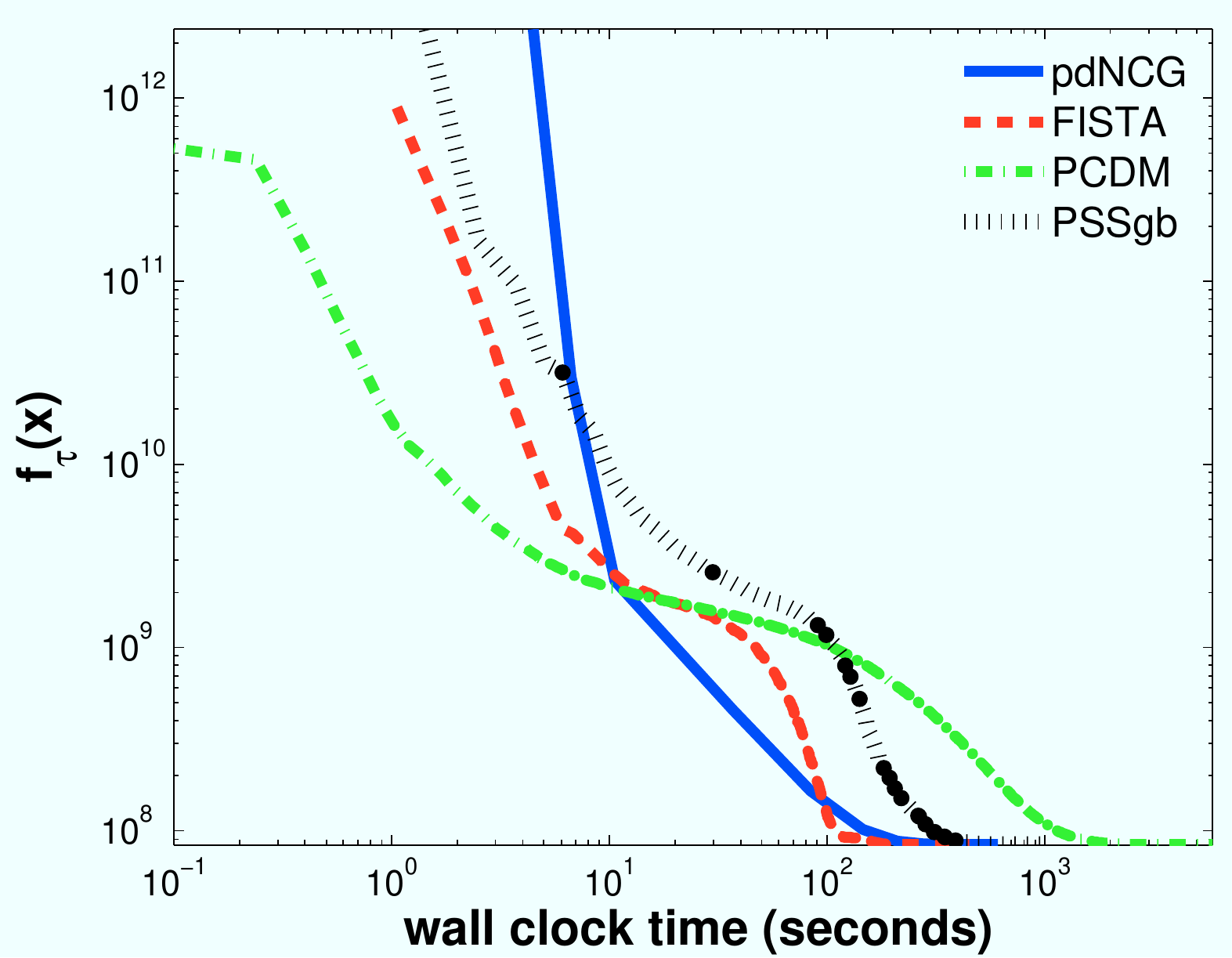}}
  \\
  \subfloat[$n = 2^{24}$]{\label{fig4f}\includegraphics[scale=0.36]{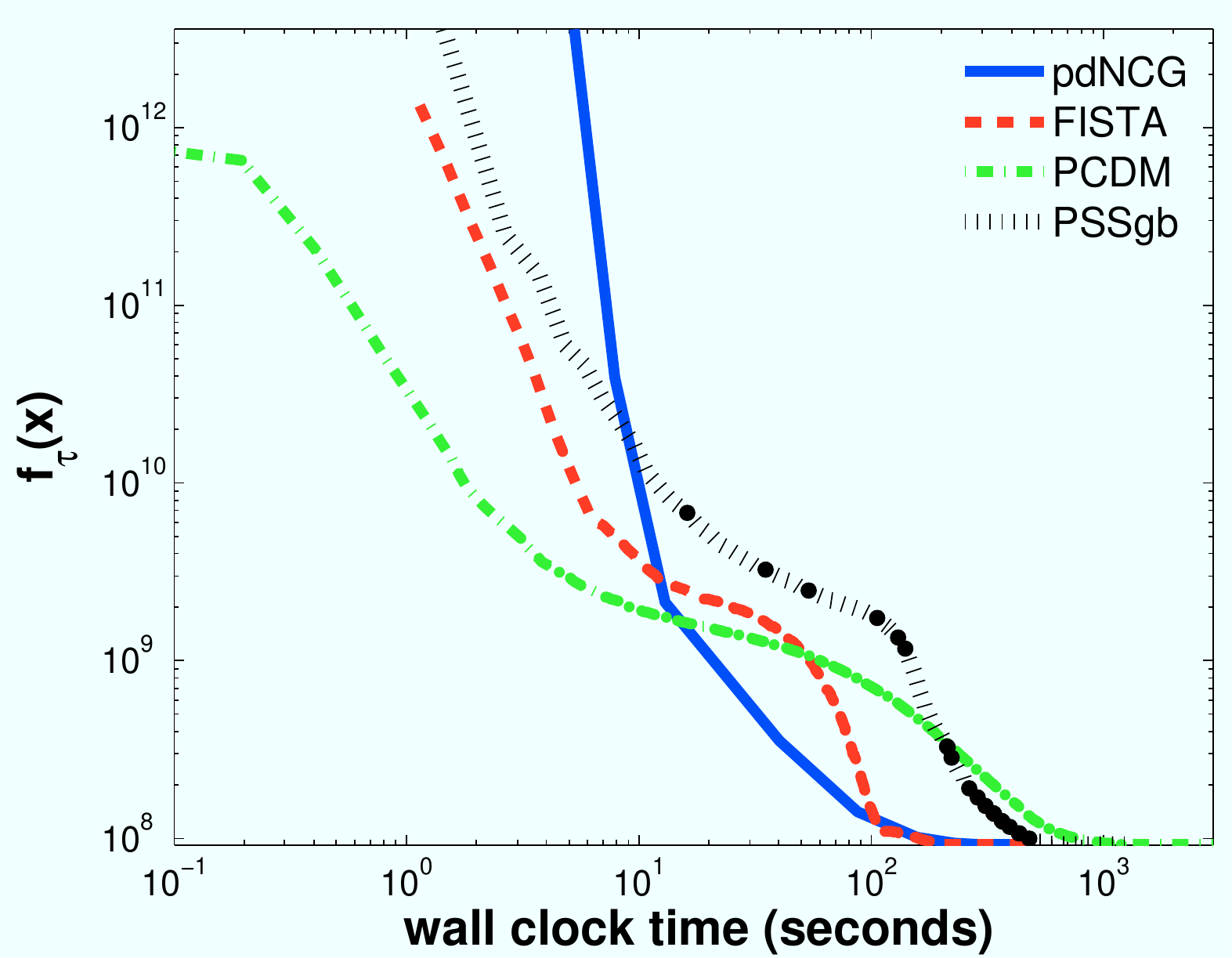}}
  \subfloat[$n = 2^{26}$]{\label{fig4h}\includegraphics[scale=0.36]{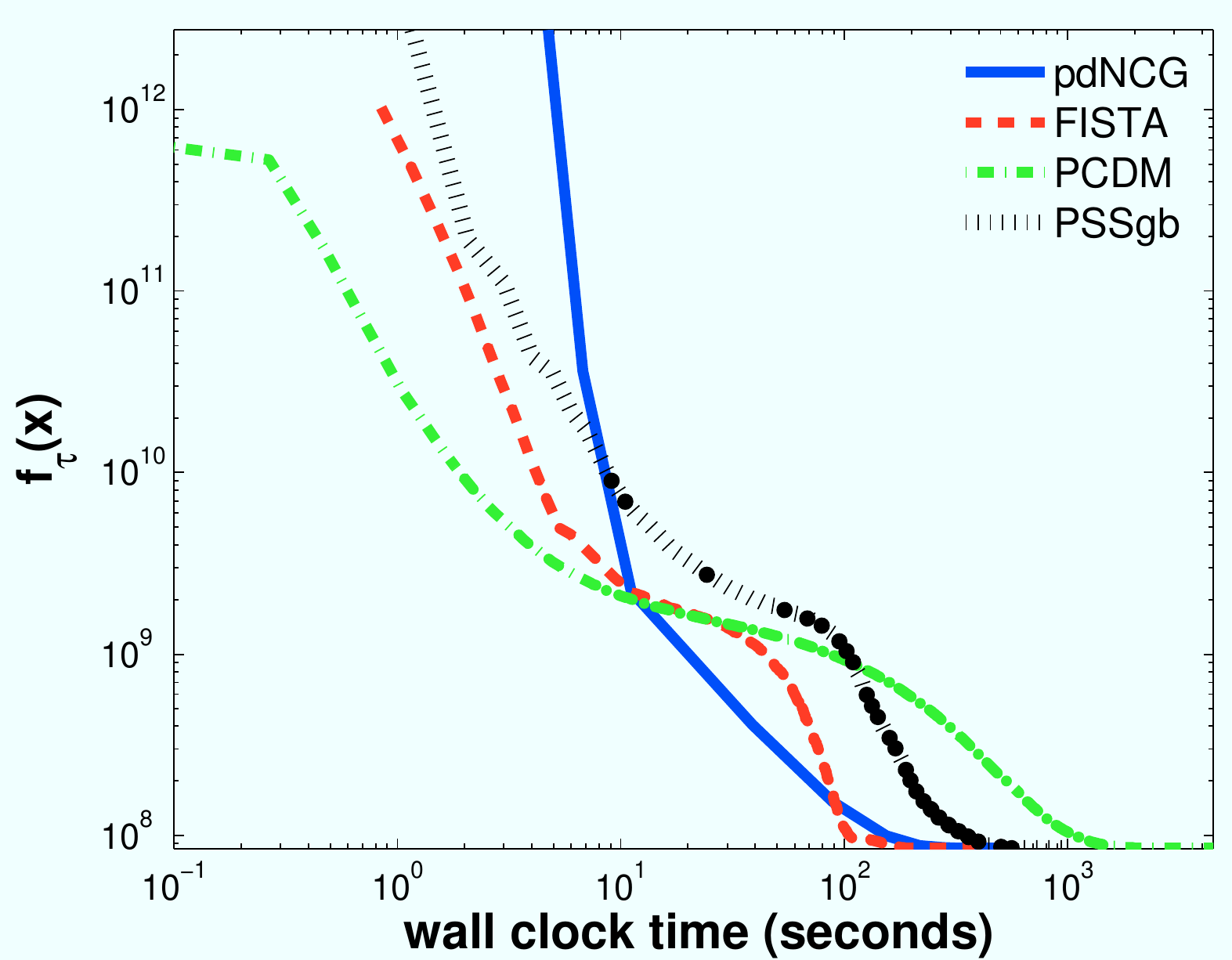}}
\caption{Performance of pdNCG, FISTA, PCDM and PSSgb on a synthetic S-LS problem for increasing number of variables $n$. 
The axes are in log-scale.}
\label{fig4}%
\end{figure}

The results of this experiment are presented in Figure \ref{fig4}. Notice that all methods have a linear-like scaling with respect to the size 
of the problem. 

\subsection{Increasing density of matrix $A^\intercal A$}
In this experiment we demonstrate the performance of pdNCG, FISTA, PCDM and PSSgb as the density of matrix $A^\intercal A$ 
increases. 
We generate four instances $(A,x^*)$. 
For the first experiment we generate matrix $A = \Sigma G^\intercal$, where $\Sigma$ is the matrix of singular values, the columns of matrices $I_m$ and $G$ are 
the left and right singular vectors, respectively.
For the second experiment we generate matrix $A=\Sigma (G_2G)^\intercal $, 
where the columns of matrices $I_m$ and $G_2G$ are the left and right singular vectors of matrix $A$, respectively; $G_2$ has been defined in Subsection \ref{subsec:constA_2}.
Finally, for the third and fourth experiments we have $A=\Sigma (GG_2G)^\intercal $ and $A=\Sigma (G_2GG_2G)^\intercal $, respectively.  
For each experiment the singular values of matrix $A$ are chosen uniformly at random in the interval $[0, 10]$ and then are shifted by $10^{-1}$,
which resulted in $\kappa(A^\intercal A)\approx10^4$.
The rotation angle $\theta$ of matrices $G$ and $G_2$ is set to $2\pi/10$ radians.
Matrices $A$ have $m= 2 n$ rows, rank $n$ and $n=2^{22}$.
The optimal solutions $x^*$ have $s = n/2^7$ nonzero components for each experiment.
Moreover, Procedure OsGen3 is used with $\gamma = 100$ and $s_1 = s_2= s/2$ for the construction of $x^*$ for each experiment, which resulted 
in $\kappa_{0.1}(x^*) \approx 2$ on average. 

\begin{figure}[htbp]%
  \centering

  \subfloat[$nnz(A^\intercal A)=2^{23}$]{\label{fig_sparsitya}\includegraphics[scale=0.36]{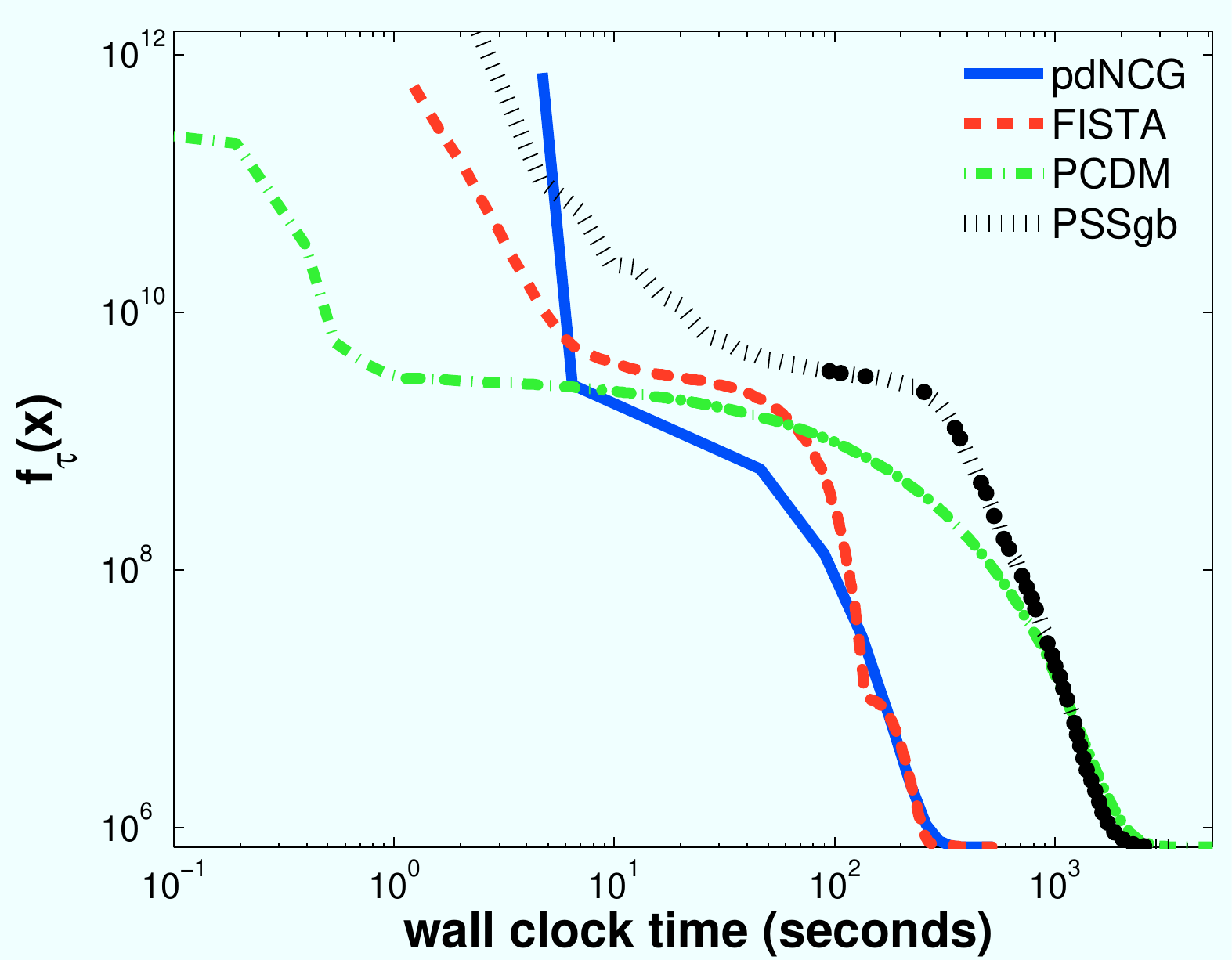}}
  \subfloat[$nnz(A^\intercal A)\approx 2^{24}$]{\label{fig_sparsityb}\includegraphics[scale=0.36]{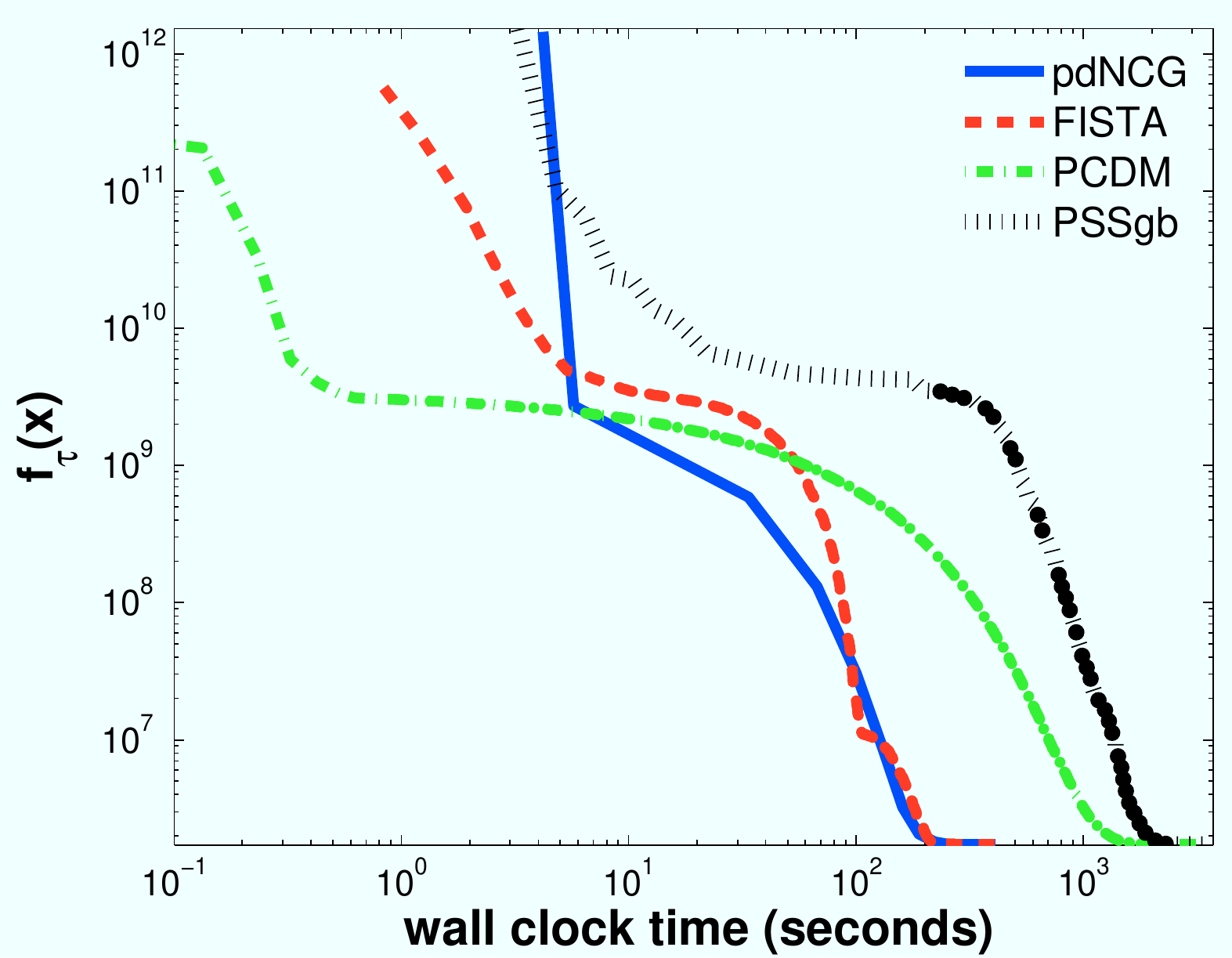}}
  \\
  \subfloat[$nnz(A^\intercal A)\approx 2^{24} + 2^{23}$]{\label{fig_sparsityc}\includegraphics[scale=0.36]{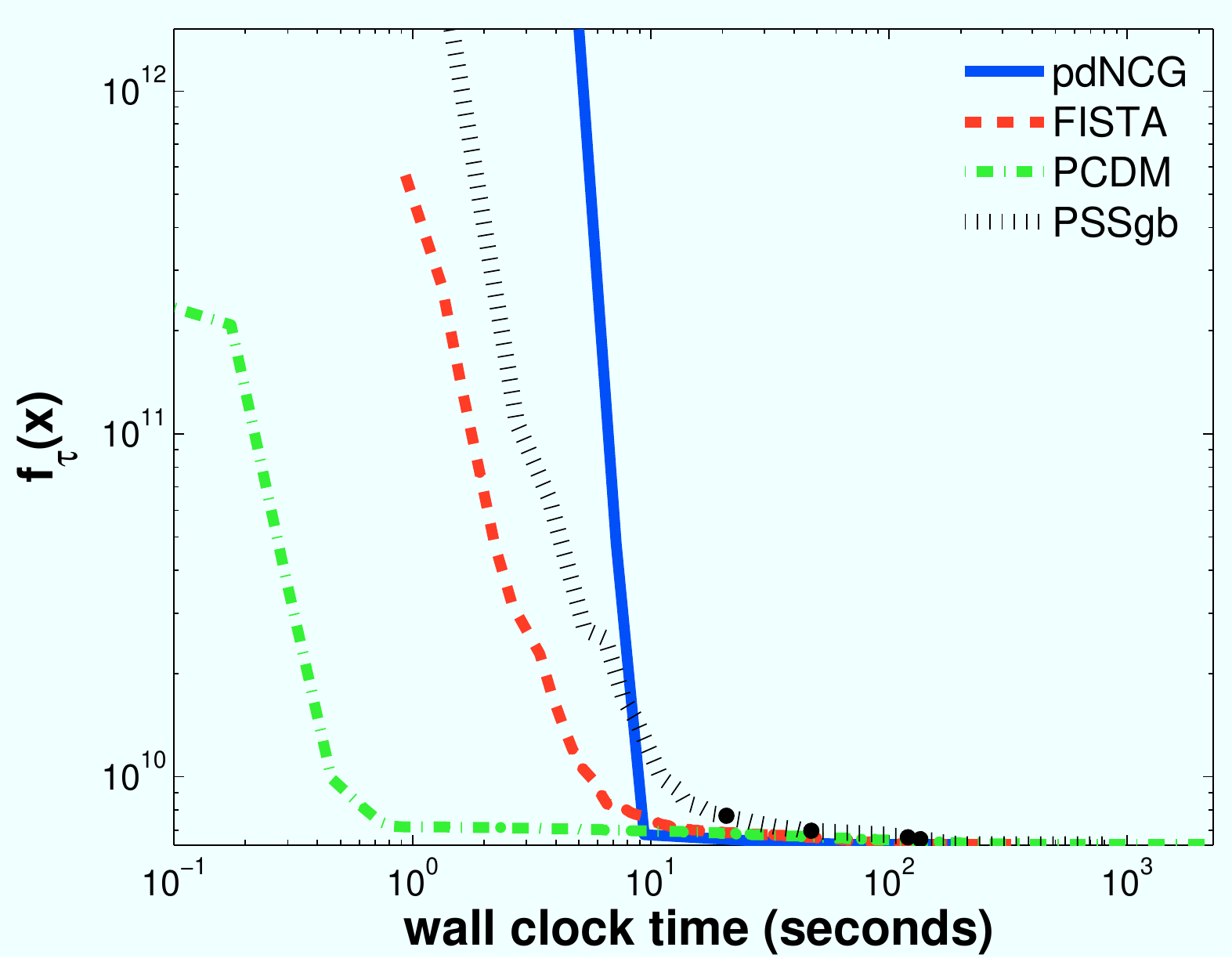}}
  \subfloat[$nnz(A^\intercal A)\approx 2^{25}$]{\label{fig_sparsityd}\includegraphics[scale=0.36]{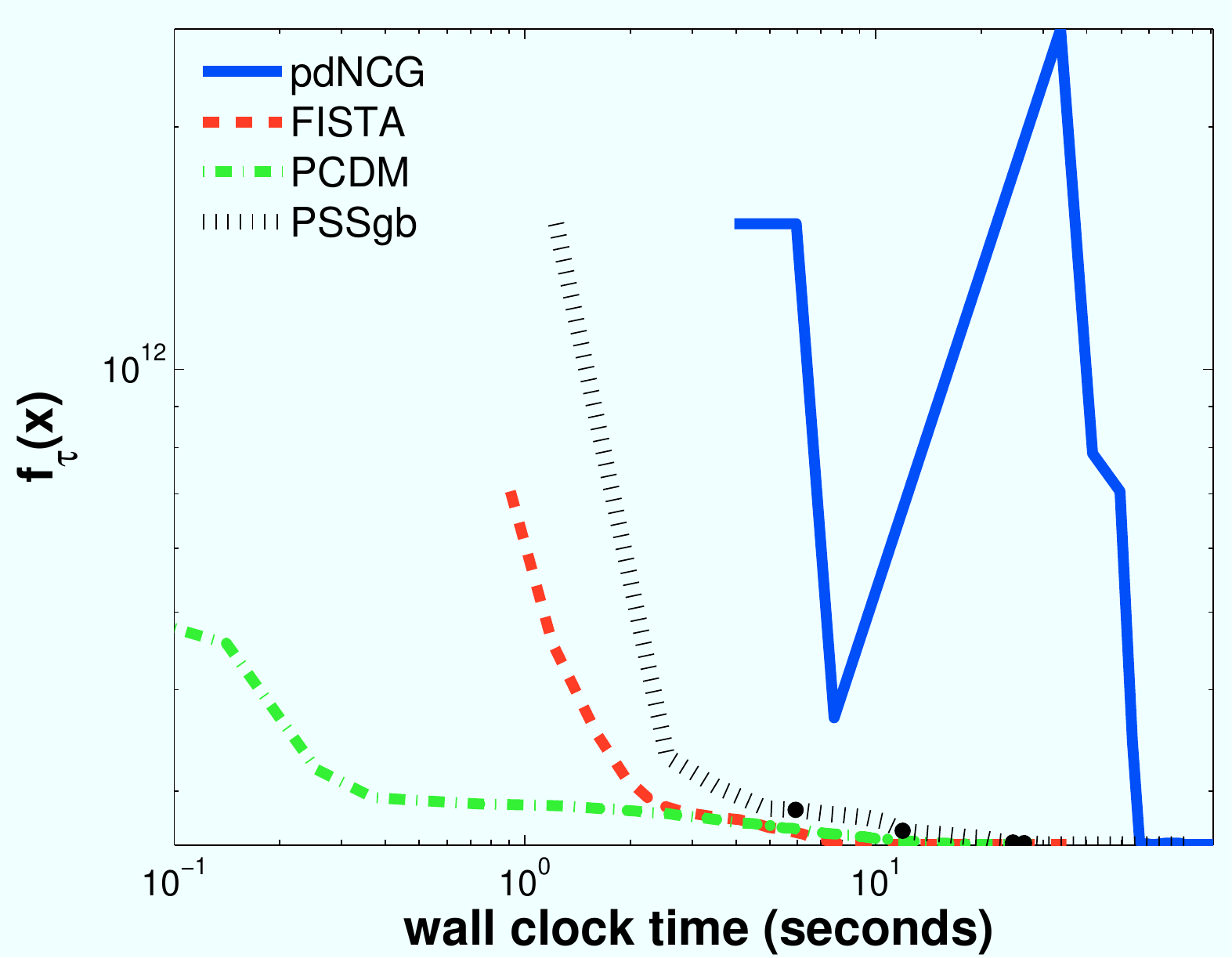}}
\caption{Performance of pdNCG, FISTA, PCDM and PSSgb on synthetic S-LS problems for increasing number of nonzeros of matrix $A$.
The axes are in log-scale.}
\label{fig_sparsity}%
\end{figure}

The results of this experiment are presented in Figure \ref{fig_sparsity}. Observe, that all methods had a robust performance 
with respect to the density of matrix $A$.

\subsection{Varying parameter $\tau$} 
In this experiment we present the performance of pdNCG, FISTA, PCDM and PSSgb as parameter $\tau$ varies from $10^{-4}$ to $10^4$ with a step of times $10^2$.
We generate four instances $(A,x^*)$, where matrix $A=\Sigma G^\intercal$ has $m= 2 n$ rows, rank $n$ and $n=2^{22}$.
The singular values of matrices $A$ are chosen uniformly at random in the interval $[0, 10]$ and then are shifted by $10^{-1}$,
which resulted in $\kappa(A^\intercal A)\approx 10^4$ for each experiment. 
The rotation angles $\theta$ for matrix $G$ in $A$ is set to $2\pi/10$ radians.
The optimal solution $x^*$ has $s = n/2^7$ nonzero components for all instances.
Moreover, the optimal solutions are generated using Procedure OsGen3 with $\gamma = 100$, which resulted 
in $\kappa_{0.1}(x^*) \approx 3$ for all four instances. 

\begin{figure}[htbp]%
  \centering

  \subfloat[$\tau=10^{-4}$]{\label{fig_taua}\includegraphics[scale=0.36]{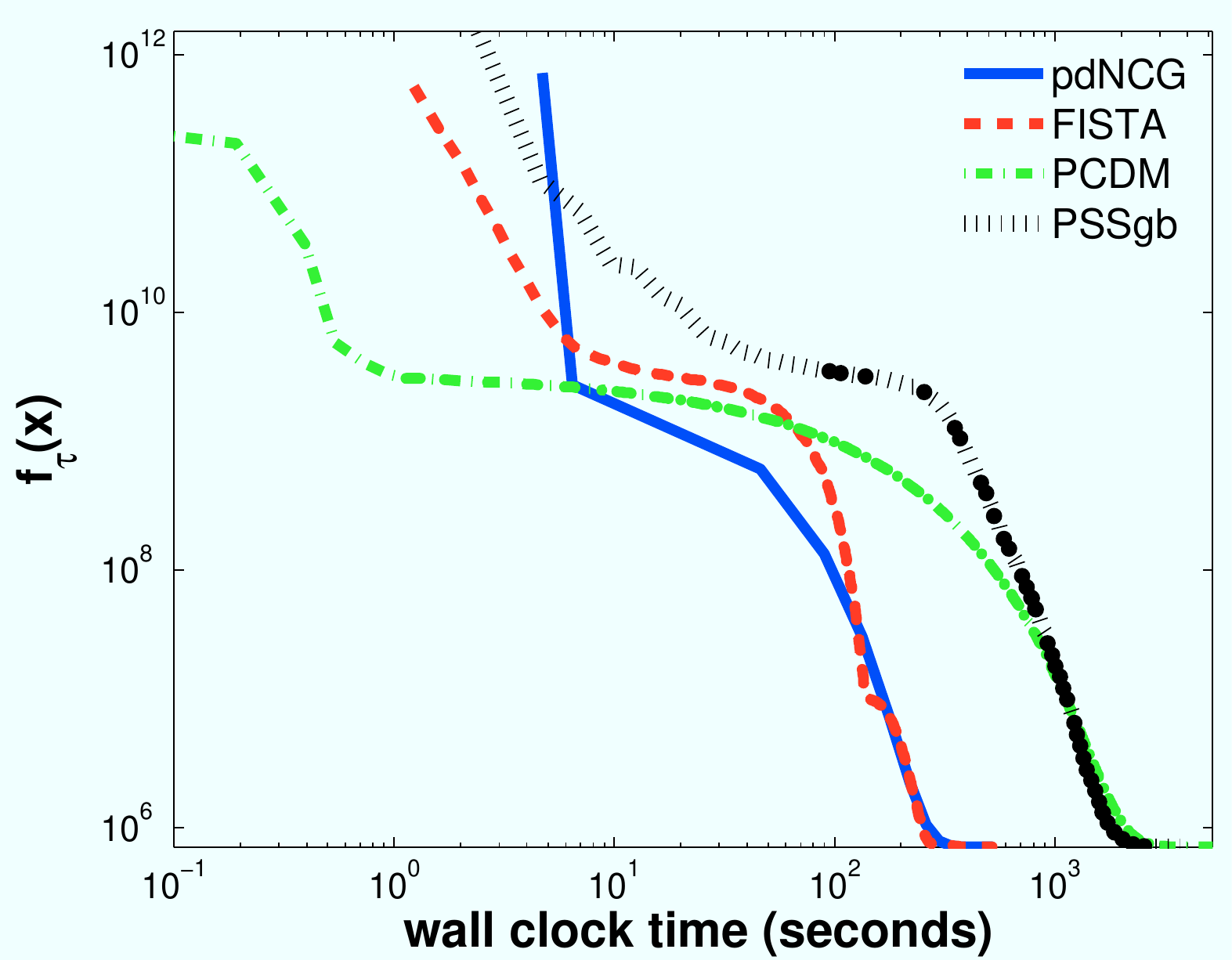}}
  \subfloat[$\tau=10^{-2}$]{\label{fig_taub}\includegraphics[scale=0.36]{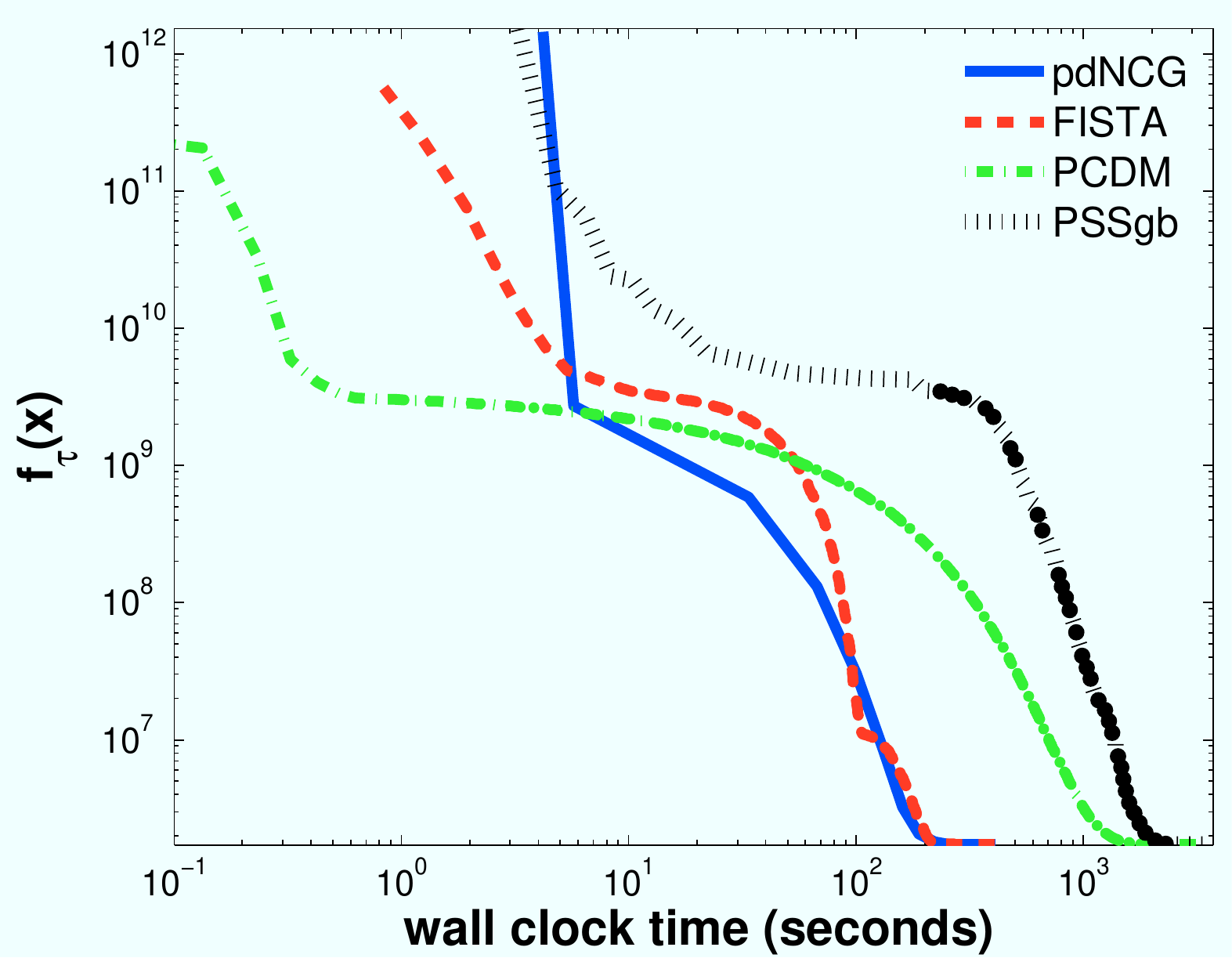}}
  \\
  \subfloat[$\tau=10^2$]{\label{fig_tauc}\includegraphics[scale=0.36]{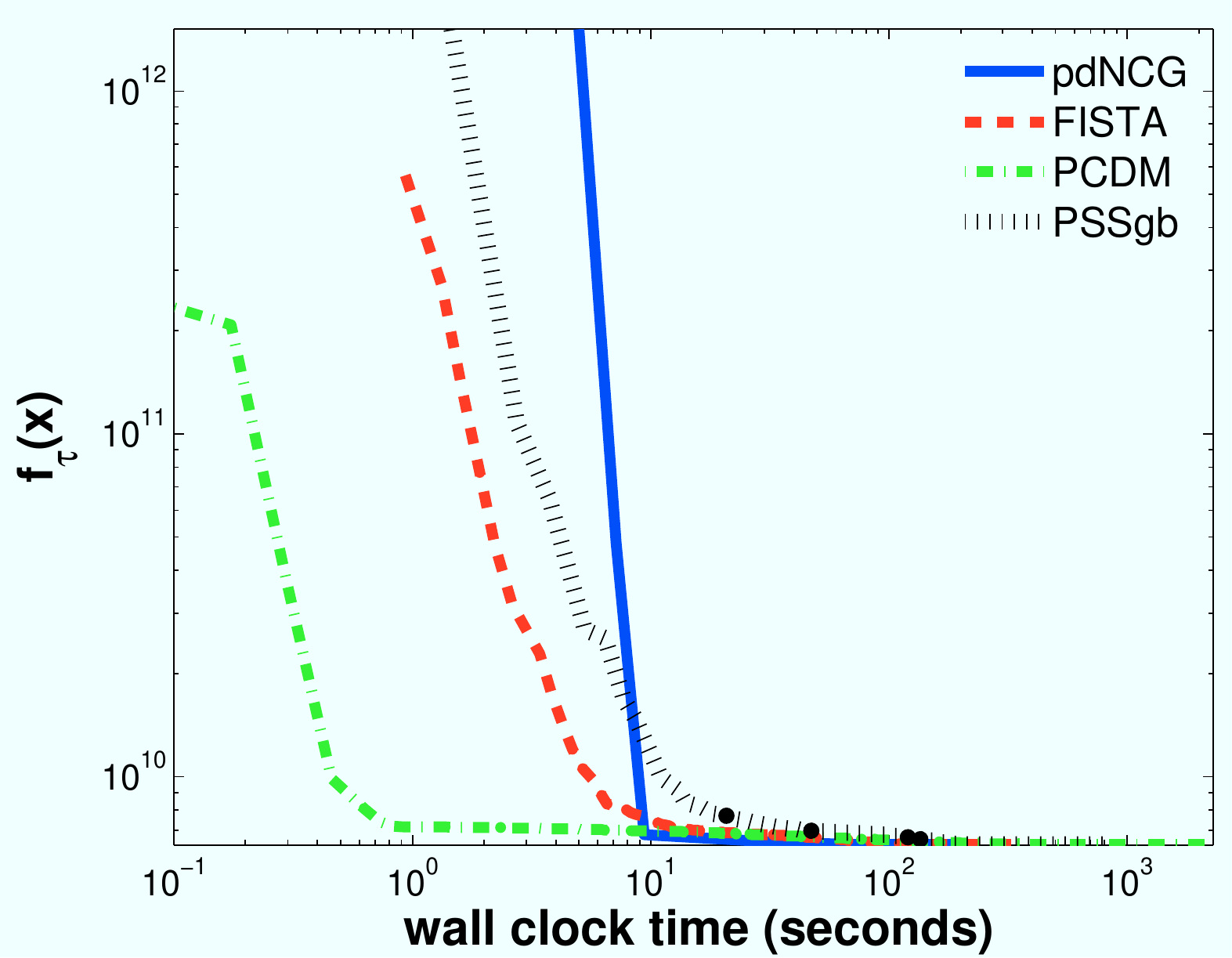}}
  \subfloat[$\tau=10^{4}$]{\label{fig_taud}\includegraphics[scale=0.36]{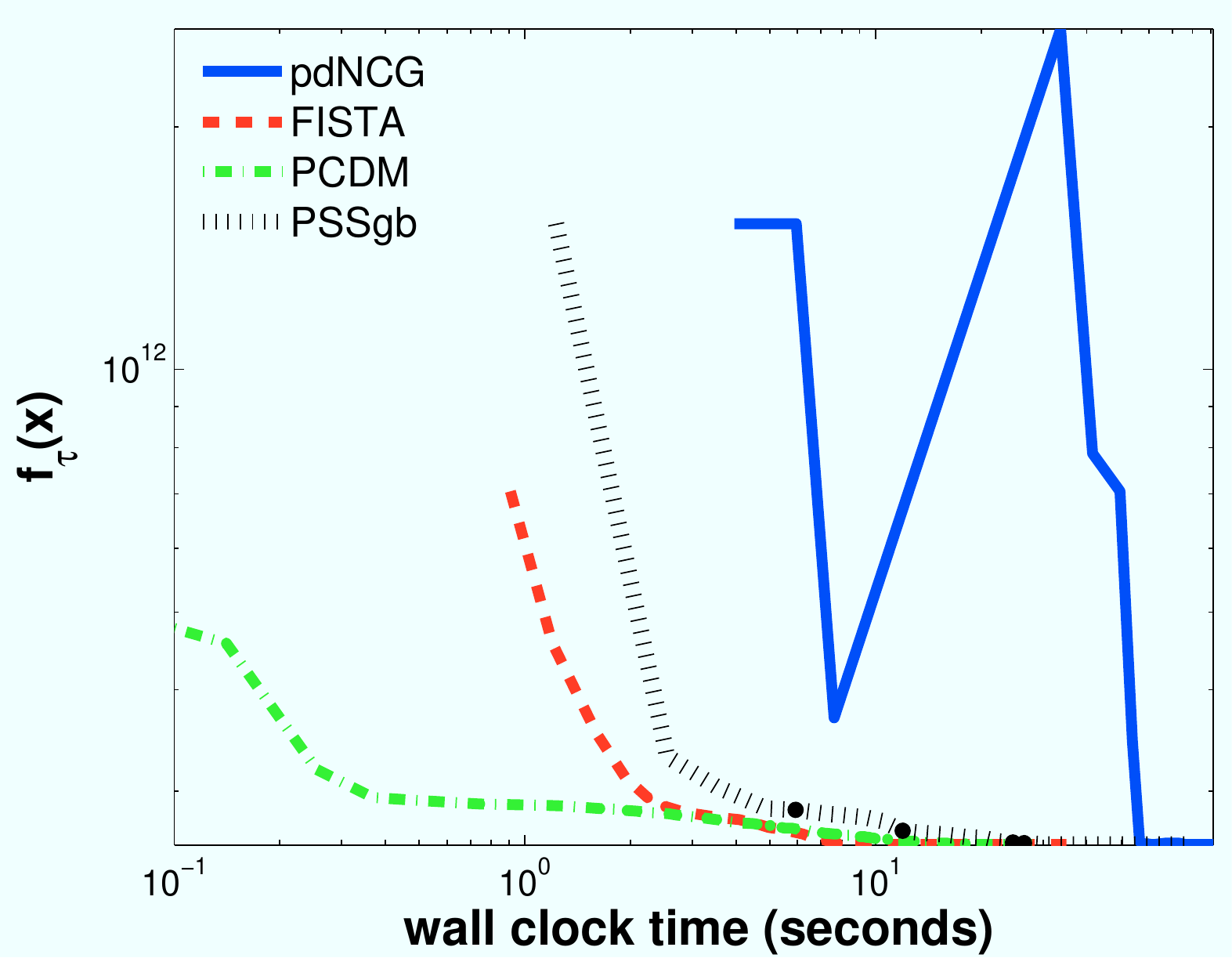}}
\caption{Performance of pdNCG, FISTA, PCDM and PSSgb on synthetic S-LS problems for various values of parameter $\tau$.
The axes are in log-scale. Observe in Subfigure \ref{fig_taud} that for pdNCG the objective function $f_\tau$
is not always decreasing monotonically. This is due to the 
backtracking line-search of pdNCG, which terminates in case that the maximum number of backtracking iterations 
is exceeded regardless if certain termination criteria are satisfied. }
\label{fig_tau}%
\end{figure}

The performance of the methods is presented in Figure \ref{fig_tau}. 
Notice in Subfigure \ref{fig_taud} that for pdNCG the objective function $f_\tau$
is not always decreasing monotonically. A possible explanation is that the 
backtracking line-search of pdNCG, which guarantees monotonic decrease of the objective function \cite{2ndpaperstrongly}, 
terminates in case that $50$ backtracking iterations  are exceeded, regardless if certain termination criteria are satisfied.

\subsection{Performance of a second-order method on huge scale problems}\label{subsec:bigproblem}
We now present the performance of pdNCG on synthetic huge scale (up to one trillion variables) S-LS problems as the number of variables and the number of processors 
increase. 

We generate six instances $(A,x^*)$, where the number of variables $n$ takes values $2^{30}$, $2^{32}$,
$2^{34}$, $2^{36}$, $2^{38}$ and $2^{40}$. Matrices $A=\Sigma G^\intercal $ have $m= 2 n$ rows and rank $n$.
The singular values $\sigma_i$ for $i=1,2,\cdots,n$ of matrices $A$ are set to $10^{-1}$ for odd $i$'s and $10^2$ for even $i$'s. 
The rotation angle $\theta$ of matrix $G$ is set to $2\pi/3$ radians.
The optimal solutions $x^*$ have $s = n/2^{10}$ nonzero components for each experiment. 
In order to simplify
the practical generation of this problem the optimal solutions $x^*$ are set to have $s/2$ components equal
to $-10^4$ and the rest of nonzero components are set equal to $10^{-1}$.

Details of the performance of pdNCG are given in Table \ref{table1}. Observe the nearly linear scaling of pdNCG with respect to the number of variables $n$
and the number of processors. For all experiments in Table \ref{table1} pdNCG required $8$ Newton steps to converge, $100$ PCG iterations per Newton step on average,
where every PCG iteration requires two matrix-vector products with matrix $A$.
\begin{table}  
\centering
\begin{tabular}{ccccc}
\toprule
 \multicolumn{1}{c}{$n$}    & Processors  & Memory (terabytes) & Time (seconds) \\ 
\midrule
 $2^{30}$ & $64$ & $0.192$  & 1,923  \\ 
 $2^{32}$ & $256$ & $0.768$ &  1,968 \\ 
 $2^{34}$ & $1024$ & $3.072$ &  1,986 \\ 
 $2^{36}$ & $4096$ & $12.288$ &  1,970 \\ 
 $2^{38}$ & $16384$ & $49.152$ &  1,990 \\ 
 $2^{40}$ & $65536$ & $196.608$ & 2,006 \\ 
\bottomrule
 \end{tabular}
 \caption{Performance of pdNCG for synthetic huge scale S-LS problems. All problems have been solved to a relative error of order $10^{-4}$ of the obtained solution}
\label{table1}
\end{table}


\section{Conclusion}\label{sec:con}
In this paper we developed an instance generator for $\ell_1$-regularized sparse least-squares problems. 
The generator is aimed for the construction of very large-scale instances. 
Therefore it scales well as the number of variables increases, both in terms of memory requirements and time.
Additionally, the generator allows control of the conditioning and the sparsity of the problem. 
Examples are provided on how to exploit the previous advantages of the proposed generator.
We believe that the optimization community needs such a generator to be able to perform fair 
assessment of new algorithms.

Using the proposed generator we constructed very large-scale sparse instances (up to one trillion variables), 
which vary from very well-conditioned to moderately ill-conditioned. We examined the performance of several representative 
first- and second-order optimization methods. The experiments revealed that \textit{regardless of the size
of the problem}, the performance of the methods crucially depends on the conditioning of the problem. 
In particular, the first-order methods PCDM and FISTA are faster for problems with small or moderate condition number, whilst,
the second-order method pdNCG is much more efficient for ill-conditioned problems. 

\begin{acknowledgements}
This work has made use of the resources provided by ARCHER (\url{http://www.archer.ac.uk/}), made available through the Edinburgh
Compute and Data Facility (ECDF) (http://www.ecdf.ed.ac.uk/).

The authors are grateful to Dr Kenton D' Mellow for providing guidance and helpful suggestions regarding the use of ARCHER and the solution of large scale problems.
\end{acknowledgements}


\bibliographystyle{plain}
\bibliography{KFandJG.bib}

\begin{thebibliography}{10}

\bibitem{SparsityInducing}
F.~Bach, R.~Jenatton, J.~Mairal, and G.~Obozinski.
\newblock Optimization with sparsity-inducing penalties.
\newblock {\em Journal Foundations and Trends in Machine Learning},
  4(1):1--106, 2012.

\bibitem{fista}
A.~Beck and M.~Teboulle.
\newblock A fast iterative shrinkage-thresholding algorithm for linear inverse
  problems.
\newblock {\em SIAM J. Imaging Sci.}, 2(1):183--202, 2009.

\bibitem{cosampImpl}
S.~Becker.
\newblock Co{S}a{MP} and {OMP} for sparse recovery.
\newblock
  \url{http://www.mathworks.co.uk/matlabcentral/fileexchange/32402-cosamp-and-omp-for-sparse-recovery},
  2012.

\bibitem{NIPS2012_4523}
S.~Becker and J.~Fadili.
\newblock A quasi-{N}ewton proximal splitting method.
\newblock In F.~Pereira, C.J.C. Burges, L.~Bottou, and K.Q. Weinberger,
  editors, {\em Advances in Neural Information Processing Systems 25}, pages
  2618--2626. Curran Associates, Inc., 2012.

\bibitem{IEEEhowto:Nesta}
S.~R. Becker, J.~Bobin, and E.~J. Cand\`{e}s.
\newblock Nesta: A fast and accurate first-order method for sparse recovery.
\newblock {\em SIAM J. Imaging Sciences}, 4(1):1--39, 2011.

\bibitem{convexTemplates}
S.~R. Becker, E.~J. Cand\'{e}s, and M.~C. Grant.
\newblock Templates for convex cone problems with applications to sparse signal
  recovery.
\newblock {\em Mathematical Programming Computation}, 3(3):165--218, 2011.
\newblock Software available at \url{http://tfocs.stanford.edu}.

\bibitem{inpainting}
M.~Bertalmio, G.~Sapiro, C.~Ballester, and V.~Caselles.
\newblock Image inpainting.
\newblock {\em Proceedings of the 27th annual conference on Computer graphics
  and interactive techniques (SIGGRAPH)}, pages 417--424, 2000.

\bibitem{coap-vol54}
M.~Bertero, V.~Ruggiero, and L.~Zanni.
\newblock Special issue: Imaging 2013.
\newblock {\em Computational Optimization and Applications}, 54:211--213, 2013.

\bibitem{distributedadmm}
S.~Boyd, N.~Parikh, E.~Chu, B.~Peleato, and J.~Eckstein.
\newblock Distributed optimization and statistical learning via the alternating
  direction method of multipliers.
\newblock {\em Journal Foundations and Trends in Machine Learning},
  3(1):1--122, 2011.

\bibitem{proximalNewtonNocedal}
R.~H. Byrd, J.~Nocedal, and F.~Oztoprak.
\newblock An inexact successive quadratic approximation method for convex l-1
  regularized optimization.
\newblock {\em Math. Program., Ser. B}, 2015.
\newblock {DOI}: 10.1007/s10107-015-0941-y.

\bibitem{introtv}
A.~Chambolle, V.~Caselles, D.~Cremers, M.~Novaga, and T.~Pock.
\newblock An introduction to total variation for image analysis.
\newblock {\em Radon Series Comp. Appl. Math}, 9:263--340, 2010.

\bibitem{ista}
A.~Chambolle, R.~A. DeVore, N.~Y. Lee, and B.~J. Lucier.
\newblock Nonlinear wavelet image processing: {V}ariational problems,
  compression, and noise removal through wavelet shrinkage.
\newblock {\em IEEE Trans. Image Process.}, 7(3):319--335, 1998.

\bibitem{changHsiehLin}
K.-W. Chang, C.-J. Hsieh, and C.-J. Lin.
\newblock Coordinate descent method for large-scale $\ell_2$-loss linear
  support vector machines.
\newblock {\em Journal of Machine Learning Research}, 9:1369--1398, 2008.

\bibitem{IEEEhowto:DonohoCompSens}
D.~L. Donoho.
\newblock Compressed sensing.
\newblock {\em IEEE Trans. Inf. Theory}, 52(4):1289--1306, 2006.

\bibitem{admmtutorial}
J.~Eckstein.
\newblock Augmented lagrangian and alternating direction methods for convex
  optimization: A tutorial and some illustrative computational results.
\newblock {\em RUTCOR Research Reports}, 2012.

\bibitem{2ndpaperstrongly}
K.~Fountoulakis and J.~Gondzio.
\newblock A second-order method for strongly convex $\ell_1$-regularization
  problems.
\newblock {\em Mathematical Programming (accepted)}, 2015.
\newblock {DOI}: 10.1007/s10107-015-0875-4, Software available at
  \url{http://www.maths.ed.ac.uk/ERGO/pdNCG/}.

\bibitem{gmlnet}
J.~Friedman, T.~Hastie, and R.~Tibshirani.
\newblock Regularization paths for generalized linear models via coordinate
  descent.
\newblock {\em Journal of Machine Learning Research}, 9:627--650, 2008.

\bibitem{goldsteinadmm}
T.~Goldstein, B.~O'Donoghue, and S.~Setzer.
\newblock Fast alternating direction optimization methods.
\newblock {\em Technical report, CAM Report 12-35, UCLA}, 2012.

\bibitem{IEEEhowto:Jacekmf}
J.~Gondzio.
\newblock Matrix-free interior point method.
\newblock {\em Computational Optimization and Applications}, 51(2):457--480,
  2012.

\bibitem{haleyin}
E.~T. Hale, W.~Yin, and Y.~Zhang.
\newblock Fixed-point continuation method for $\ell_1$-minimization:
  {M}ethodology and convergence.
\newblock {\em SIAM J. Optim.}, 19(3):1107--1130, 2008.

\bibitem{deblurringimages}
P.~C. Hansen, J.~G. Nagy, and D.~P. O'Leary.
\newblock {\em Deblurring Images: Matrices, Spectra and Filtering}.
\newblock SIAM, Philadelphia, PA., 2006.

\bibitem{datafitting}
P.~C. Hansen, V.~Pereyra, and G.~Scherer.
\newblock {\em Least Squares Data Fitting with Applications}.
\newblock JHU Press, 2012.

\bibitem{convergenceadmm}
B.~He and X.~Yuan.
\newblock On the $\mathcal{O}(1/t)$ convergence rate of alternating direction
  method.
\newblock {\em SIAM Journal on Numerical Analysis}, 50(2):700--709, 2012.

\bibitem{HsiehChang}
C.-J. Hsieh, K.-W. Chang, C.-J. Lin, S.~S. Keerthi, and S.~Sundararajan.
\newblock A dual coordinate descent method for large-scale linear {SVM}.
\newblock {\em Proceedings of the 25th international conference on Machine
  Learning, ICML 2008}, pages 408--415, 2008.

\bibitem{IEEEhowto:boyd}
S.-J. Kim, K.~Koh, M.~Lustig, S.~Boyd, and D.~Gorinevsky.
\newblock An interior-point method for large-scale $\ell_1$-regularized least
  squares.
\newblock {\em IEEE Journal Of Selected Topics In Signal Processing},
  1(4):606--617, 2007.

\bibitem{dirkLorenz}
D.~A. Lorenz.
\newblock Constructing test instances for basis pursuit denoising.
\newblock {\em IEEE Trans. Signal Process.}, 61(5):1210--1214, 2013.

\bibitem{sparseMRI}
M.~Lustig, D.~Donoho, and J.~M. Pauly.
\newblock Sparse {MRI}: {T}he application of compressed sensing for rapid {MR}
  imaging.
\newblock {\em Magnetic Resonance in Medicine}, 58(6):1182--1195, 2007.

\bibitem{cosamp}
D.~Needell and J.~A. Tropp.
\newblock Cosamp: Iterative signal recovery from incomplete and inaccurate
  samples.
\newblock {\em Applied and Computational Harmonic Analysis}, 26(3):301--321,
  2009.

\bibitem{booknesterov}
Y.~Nesterov.
\newblock {\em Introductory Lecture Notes On Convex Optimization. A Basic
  Course}.
\newblock Kluver, Boston, 2004.

\bibitem{nesterovSmooth}
Y.~Nesterov.
\newblock Smooth minimization of non-smooth functions.
\newblock {\em Mathematical Programming}, 103(1):127--152, 2005.

\bibitem{IEEEhowto:Nesterov}
Y.~Nesterov.
\newblock Smooth minimization of non-smooth functions.
\newblock {\em Math. Program.}, 103(1):127--152, 2005.

\bibitem{nesterovgen}
Yu. Nesterov.
\newblock Gradient methods for minimizing composite functions.
\newblock {\em Mathematical Programming}, 140(1):125--161, 2013.

\bibitem{proximalmethods}
N.~Parikh and S.~Boyd.
\newblock Proximal algorithms.
\newblock {\em Journal Foundations and Trends in Optimization}, 1(3):123--231,
  2013.

\bibitem{petermartin}
P.~Richt\'{a}rik and M.~Tak\'{a}\v{c}.
\newblock Iteration complexity of randomized block-coordinate descent methods
  for minimizing a composite function.
\newblock {\em Math. Program. Ser. A}, 144(1):1--38, 2014.

\bibitem{peterbigdata}
P.~Richt\'{a}rik and M.~Tak\'{a}\v{c}.
\newblock Parallel coordinate descent methods for big data optimization.
\newblock {\em Math. Program. Ser. A}, pages 1--52, 2015.
\newblock {DOI}: 10.1007/s10107-015-0901-6.

\bibitem{proximalNewtonKatya}
K.~Scheinberg and X.~Tang.
\newblock Practical inexact proximal quasi-{N}ewton method with global
  complexity analysis.
\newblock Technical report, March 2014.
\newblock arXiv:1311.6547 [cs.LG].

\bibitem{thesisschmidt}
M.~Schmidt.
\newblock Graphical model structure learning with l1-regularization.
\newblock {\em PhD thesis, University British Columbia}, 2010.

\bibitem{shwartzTewari}
S.~Shalev-Shwartz and A.~Tewari.
\newblock Stochastic methods for $\ell_1$-regularized loss minimization.
\newblock {\em Journal of Machine Learning Research}, 12(4):1865--1892, 2011.

\bibitem{tsengblkcoo}
P.~Tseng.
\newblock Convergence of a block coordinate descent method for
  nondifferentiable minimization.
\newblock {\em Journal of Optimization Theory and Applications},
  109(3):475--494, 2001.

\bibitem{nesterovhuge}
P.~Tseng.
\newblock Efficiency of coordinate descent methods on huge-scale optimization
  problems.
\newblock {\em SIAM J. Optim.}, 22:341--362, 2012.

\bibitem{tsengyun}
P.~Tseng and S.~Yun.
\newblock A coordinate gradient descent method for nonsmooth separable
  minimization.
\newblock {\em Math. Program., Ser. B}, 117:387--423, 2009.

\bibitem{gwacs}
S.~Vattikuti, J.~J. Lee, C.~C. Chang, S.~D. Hsu, and C.~C. Chow.
\newblock Applying compressed sensing to genome-wide association studies.
\newblock {\em GigaScience}, 3(10):1--17, 2014.

\bibitem{admmtv}
Y.~Wang, J.~Yang, W.~Yin, and Y.~Zhang.
\newblock A new alternating minimization algorithm for total variation image
  reconstruction.
\newblock {\em SIAM Journal on Imaging Sciences}, 1(3):248--272, 2008.

\bibitem{wrightaccel}
S.~J. Wright.
\newblock Accelerated block-coordinate relaxation for regularized optimization.
\newblock {\em SIAM Journal on Optimization}, 22(1):159--186, 2012.

\bibitem{tonglange}
T.~T. Wu and K.~Lange.
\newblock Coordinate descent algorithms for lasso penalized regression.
\newblock {\em The Annals of Applied Statistics}, 2(1):224--244, 2008.

\bibitem{svmcomparison}
G.-X. Yuan, K.-W. Chang, C.-J. Hsieh, and C.-J. Lin.
\newblock A comparison of optimization methods and software for large-scale
  $\ell_1$-regularized linear classification.
\newblock {\em Journal of Machine Learning Research}, 11:3183--3234, 2010.

\end{thebibliography}

\end{document}